\documentclass{amsart}
\usepackage{color}

\usepackage{amsmath}
\usepackage{amssymb}
\usepackage{mathrsfs}
\usepackage{enumitem}
\usepackage{moreenum}

\usepackage[dvipsnames]{xcolor}

\newcommand{\tieconcat}{%
	\mathbin{\mathpalette\dotieconcat\relax}%
}
\newcommand{\dotieconcat}[2]{
	\text{\raisebox{.8ex}{$\smallfrown$}}%
}

\newtheorem{theorem}{Theorem}[section]
\newtheorem{claim}[theorem]{Claim}

\newtheorem{conclusion}[theorem]{Conclusion}
\newtheorem{observation}[theorem]{Observation}

\newtheorem{ttt}[theorem]{Theorem}
\newtheorem{llll}[theorem]{Lemma}
\newtheorem{ccc}[theorem]{Claim}
\newtheorem{sccc}[theorem]{Subclaim}
\newtheorem{eee}[theorem]{Example}
\newtheorem{fff}[theorem]{Fact}
\newtheorem{rrr}[theorem]{Remark}
\newtheorem{sss}[theorem]{Statement}
\newtheorem{ddd}[theorem]{Definition}
\newtheorem{qqq}[theorem]{Question}
\newtheorem{cccc}[theorem]{Corollary}
\newtheorem{nnn}[theorem]{Notation}
\newtheorem{ppp}[theorem]{Problem}
\newtheorem{pppp}[theorem]{Proposition}
\newtheorem{ccccc}[theorem]{Conjecture}

\theoremstyle{definition}
\newtheorem{definition}[theorem]{Definition}

\newtheorem{convention}[theorem]{Convention}

\newtheorem{discussion}[theorem]{Discussion}

\theoremstyle{remark}
\newtheorem{remark}[theorem]{Remark}
\newtheorem{question}[theorem]{Question}

\newcommand{\beq}{\begin{equation} }
	
\newcommand{\bt}{\begin{ttt}}
\newcommand{\bl}{\begin{llll}}
\newcommand{\bcl}{\begin{ccc}}
	\newcommand{\bscl}{\begin{sccc}}
\newcommand{\bex}{\begin{eee}}
\newcommand{\bfact}{\begin{fff}}
\newcommand{\br}{\begin{rrr}\upshape}
\newcommand{\bst}{\begin{sss}}
\newcommand{\bdd}{\begin{ddd}\upshape}
										
\newcommand{\bqu}{\begin{qqq}}
\newcommand{\bnn}{\begin{nnn}}
\newcommand{\bpr}{\begin{ppp}}
\newcommand{\bprop}{\begin{pppp}}
\newcommand{\bcor}{\begin{cccc}}
\newcommand{\bcon}{\begin{ccccc}}

\newcommand{\et}{\end{ttt}}
\newcommand{\el}{\end{llll}}
\newcommand{\ecl}{\end{ccc}}
\newcommand{\escl}{\end{sccc}}
\newcommand{\eex}{\end{eee}}
\newcommand{\efact}{\end{fff}}
\newcommand{\er}{\end{rrr}}
\newcommand{\est}{\end{sss}}
\newcommand{\edd}{\end{ddd}}
\newcommand{\equ}{\end{qqq}}
\newcommand{\ecor}{\end{cccc}}
\newcommand{\econ}{\end{ccccc}}
\newcommand{\enn}{\end{nnn}}
\newcommand{\epr}{\end{ppp}}
\newcommand{\eprop}{\end{pppp}}

\newcommand{\otp}{{\rm otp}}

\newcommand{\feq}{{\rm feq}}
\newcommand{\ceq}{{\rm ceq}}
\newcommand{\Mod}{{\rm Mod}}

\newcommand{\Mos}{{\rm Mos}}
\newcommand{\Col}{{\rm Col}}
\newcommand{\NSOP}{{\rm NSOP}}

\newcommand{\card}{{\rm Card}}
\newcommand{\reg}{{\rm Reg}}
\newcommand{\ord}{{\rm ORD}}

\newcommand{\pr}{{\rm pr}}

\newcommand{\univ}{{\rm univ}}

\newcommand{\bfj}{{\mathbf j}}

\newcommand{\bfp}{{\mathbf p}}

\newcommand{\bfA}{{\mathbf A}}
\newcommand{\bfM}{{\mathbf M}}
\newcommand{\bfG}{{\mathbf G}}

\newcommand{\bfV}{{\mathbf V}}

\newcommand{\bfr}{{\mathbf r}}
\newcommand{\bfv}{{\mathbf v}}

\newcommand{\bfh}{{\mathbf h}}

\newcommand{\dom}{{\rm dom}}
\newcommand{\ran}{{\rm ran}}

\newcommand{\cf}{{\rm cf}}

\newcommand{\Rang}{{\rm Rang}}

\newcommand{\lh}{{\ell g}}

\newcommand{\wilog}{{\rm without loss of generality}}

\newcommand{\then}{{\underline{then}}}
\newcommand{\when}{{\underline{when}}}

\newcommand{\Then}{{\underline{Then}}}

\newcommand{\Iff}{{\underline{iff}}}
\newcommand{\mn}{{\medskip\noindent}}
\newcommand{\sn}{{\smallskip\noindent}}

\newcommand{\cA}{{\mathscr A}}

\newcommand{\varp}{{\varepsilon}}

\newcommand{\cH}{{\mathscr H}}

\newcommand{\cN}{{\mathscr N}}
\newcommand{\bbP}{{\mathbb P}}
\newcommand{\bbR}{{\mathbb R}}
\newcommand{\bbQ}{{\mathbb Q}}
 \newcommand{\cP}{{\mathscr P}}
  
  \newcommand{\crr}{\textrm{crit}}
  \newcommand{\um}{\upharpoonright}

\newcommand{\cS}{{\mathscr S}}

\newcommand{\beeq}{\begin{equation}}
\newcommand{\eeq}{\end{equation}}

\def\seq{\subseteq}
\def\se{\setminus}

\newcount\skewfactor
\def\mathunderaccent#1#2 {\let\theaccent#1\skewfactor#2
\mathpalette\putaccentunder}
\def\putaccentunder#1#2{\oalign{$#1#2$\crcr\hidewidth
\vbox to.2ex{\hbox{$#1\skew\skewfactor\theaccent{}$}\vss}\hidewidth}}
\def\name{\mathunderaccent\tilde-3 }

\newbox\noforkbox \newdimen\forklinewidth
\setbox0\hbox{$\textstyle\bigcup$}
\setbox1\hbox to \wd0{\hfil\vrule width \forklinewidth depth \dp0
                        height \ht0 \hfil}
\setbox\noforkbox\hbox{\box1\box0\relax}
\def\unionstick{\mathop{\copy\noforkbox}\limits}
\def\nonfork#1#2_#3{#1\unionstick_{\textstyle #3}#2}
\def\nonforkin#1#2_#3^#4{#1\unionstick_{\textstyle #3}^{\textstyle
    #4}#2}
\setbox0\hbox{$\textstyle\bigcup$}
\setbox1\hbox to \wd0{\hfil{\sl /\/}\hfil}
\setbox2\hbox to \wd0{\hfil\vrule height \ht0 depth \dp0 width
                                \forklinewidth\hfil}
\newbox\doesforkbox
\setbox\doesforkbox\hbox{\box1\box0\relax}
\def\nunionstick{\mathop{\copy\doesforkbox}\limits}

\def\fork#1#2_#3{#1\nunionstick_{\textstyle #3}#2}
\def\forkin#1#2_#3^#4{#1\nunionstick_{\textstyle #3}^{\textstyle
    #4}#2}

\newcommand{\stickT}{%
\setbox255=\hbox{\raise1ex\hbox{$\hspace{0.2pt}\,\bullet\,$}}
\mathord{\rlap{\hbox to\wd255{\hss\hbox{$|$}\hss}}
\box255}
}
\newcommand{\stickS}{%
\setbox255=\hbox{\raise0.6ex\hbox{$\scriptstyle\bullet$}}
\mathord{\rlap{\hbox to\wd255{\hss\hbox{$\scriptstyle|$}\hss}}
\box255}
}

\newenvironment{PROOF}[2][\proofname.]
   {\begin{proof}[#1]\renewcommand{\qedsymbol}{\textsquare$_{\rm #2}$}}
   {\end{proof}}

\usepackage{hyperref}
\begin{document}

\title {Universal graphs between a strong limit singular and its power}
\author{Márk Poór$^\dag$}
\address{Einstein Institute of Mathematics\\
	Edmond J. Safra Campus, Givat Ram\\
	The Hebrew University of Jerusalem\\
	Jerusalem, 9190401, Israel}
\author {Saharon Shelah$^\ast$}
\address{Einstein Institute of Mathematics\\
Edmond J. Safra Campus, Givat Ram\\
The Hebrew University of Jerusalem\\
Jerusalem, 9190401, Israel\\
 and \\
 Department of Mathematics\\
 Hill Center - Busch Campus \\
 Rutgers, The State University of New Jersey \\
 110 Frelinghuysen Road \\
 Piscataway, NJ 08854-8019 USA}
\email{shelah@math.huji.ac.il}
\urladdr{http://shelah.logic.at}
\thanks{$^\dag$The first author was supported by the Excellence Fellowship Program for International Postdoctoral Researchers of The Israel Academy of Sciences and Humanities, and by the National Research, Development and Innovation Office
	– NKFIH, grants no. 124749, 129211. \\
$^\ast$The second author was supported	by the Israel Science Foundation grant
	1838/19. Paper 1185 on Shelah's list. \\
	References like \cite[Th0.2=Ly5]{Sh:950} means the label
of Th.0.2 is y5.  The reader should note that the version in my
website is usually more updated than the one in the mathematical archive.}




\date {December 2021}

\begin{abstract}
The paper settles the problem of the consistency of the existence of a single universal graph between a strong limit singular  and its power. Assuming that in a model of $\mathbf{GCH}$ $\kappa$ is supercompact and the  cardinals $\theta < \kappa$, $\lambda > \kappa$ are regular,   as an application of a more general method we obtain a forcing extension in which  $\cf(\kappa) = \theta$, the Singular Cardinal Hypothesis fails at $\kappa$ and there exists a universal graph in cardinality $\lambda \in (\kappa,2^\kappa)$.
\end{abstract}

\maketitle
\numberwithin{equation}{section}
\setcounter{section}{-1}
\newpage

\centerline {Annotated Content}
\bigskip

\noindent
\S0 \quad Introduction, pg.\pageref{0}
\bigskip

\noindent
\S1 \quad The Frame and Deducing the Consistency Results (label h),
pg.\pageref{1}
\bigskip

\noindent
\S2 \quad Proving Known Forcings Fit the Frame (label b), pg.\pageref{2}
\bigskip

\noindent
\S3 \quad The Preparatory Forcing (label d), pg.\pageref{3}
\newpage

\section {Introduction} \label{0}
\bigskip

\subsection {Background} \label{0A}\
\bigskip

The existence of universal graphs in infinite cardinalities  has been widely investigated (where we mean that  the graph $G$ is universal in cardinality $|G|$ if each graph of the same cardinality is isomorphic to some induced subgraph of $G$). 
By the classical result \cite{random}, the so called countable random graph is a universal graph on $\aleph_0$ (which is also unique up to isomorphism). A classical result (which is now a standard induction argument) yields that there is a $\kappa^+$-saturated graph on $2^\kappa$\cite{CK73}, and so a graph on $2^\kappa$ into which each graph on $\kappa^+$ embeds (and we can replace $\kappa^+$, $2^\kappa$, $\kappa^+$-saturated by $\kappa$, $2^{<\kappa}$, $\kappa$-special).
Therefore under $\mathbf{GCH}$ in every uncountable cardinality there is a universal graph.
(However, concerning certain proper class of graphs the situation is more intricate, even for the countable case, see \cite{FuKo}, \cite{Ko}, \cite{Sh:492}, \cite{Sh:1033}, \cite{Sh:1161}.) For the problem of universal objects in more complex  theories (i.e. than that of the graphs), and the relevance of the present work in model theory 
see also the survey \cite{Sh:1151}, or earlier \cite{Dj05}, see lately
\cite{Sh:1162}, \cite{Sh:1164}.

On the other hand, without assuming $\mathbf{GCH}$ it is in general much more difficult to construct universal objects, while  there are certainly no universal graphs after adding enough Cohen subsets see \cite{Sh:409}.

As for positive results, for regular cardinals $\kappa < \lambda$ consistently there is a universal graph on $\lambda$, while $2^\kappa > \lambda$ \cite{Sh:175a}. While the argument  in \cite{Sh:175a} also gives a universal $\omega$-edge colored graph on $\omega_1$ with $\neg \mathbf{CH}$ (which feature  will utilized in this paper), recently \cite{Sh:1088} proved that assuming $\neg \mathbf{CH}$ a universal graph on $\omega_1$ does not imply that there is a universal $\omega$-edge colored graph.
 (Again we remark that, if we restrict ourselves to  specific classes of graphs there are both negative \cite {Kj}, and positive results \cite{Mk90}, for weak universal families (see \ref{(*)} below) in the absence of $\mathbf{GCH}$ see \cite{Sh:457}, \cite{Sh:614}. In all the above the case $\lambda = \kappa^+$ was considerably easier.)

In the present paper we investigate universal graphs in the interval between a strong limit singular cardinal  and its power. The question is also motivated by the following. Recall that for $\mu = \aleph_0$ its power $2^{\aleph_0}$ may be large, moreover a relevant forcing axiom (e.g.\ $\mathbf{MA}$) possibly holds. Similarly for $\mu = \aleph_1 = 2^{\aleph_0}$, $2^{\mu}$ large, or $\mu = \mu^{<\mu}$ parallel results hold for forcing notions which are e.g. $<\mu$-complete, satisfying a strong form of $\mu^+$-cc (the strong form is necessary, see \cite{Sh:1036}). On the other hand for $\mu$ strong limit singular we know much less, therefore the existence of universals also serves as a central test problem regarding the consistency of forcing axioms at $\mu$.  

More directly we continue the work of D{\v z}amonja-Shelah in \cite{Sh:659}, which proved for  the case $\cf(\mu) = \aleph_0$ (assuming a supercompact) the consistency of
\mn
\begin{enumerate}[label = $(*)$, ref = $(*)$]
\item \label{(*)}
\begin{enumerate}
\item[(a)]  $\mu$ is strong limit singular  and $\mu^{++} < 2^\mu$,
\sn
\item[(b)]  there is a graph $G_*$ of cardinality $\mu^{++}$ which is
  universal for graphs of cardinality $\mu^+$ (equivalently there is a
  sequence $\bar G = \langle G_\alpha:\alpha < \mu^{++}\rangle$ of
graphs each of cardinality $\mu^+$, universal for the family of such graphs).

\end{enumerate}
\end{enumerate}
\mn
see  \cite{Sh:659} for the case $\cf(\mu) = \aleph_0$, and Cummings-D{\v z}amonja-Magidor-Morgan-Shelah  prove this for arbitrary cofinality in \cite{Sh:963}.
Earlier Mekler-Shelah \cite{Sh:274} had proved such consistency results
replacing (b) by uniformization results; also starting naturally with a
supercompact cardinal. Later, \ref{(*)} was proved to be consistent for small singular $\mu$'s, see \cite{Sma}, \cite{Unio}.

Our aim is to solve the problem arising naturally there: First replacing weak universal by universal. Second, replacing $\lambda = \mu^+$ by $\lambda \in (\mu, 2^\mu)$, so formulating the following assertion the following assertion:
\mn
\begin{enumerate}[label = $(*)^+$, ref = $(*)^+$]
\item \label{(*)+}
\begin{enumerate}
\item[(a)]  $\mu$ is strong limit singular and $\mu^{++} < 2^\mu$,
\sn
\item[(b)]  there is a universal graph $G_*$ in $\mu^+$, i.e. universal for
  graphs of cardinality $\mu^+$, $G_*$ itself of cardinality $\mu^+$,
\item[(b)$^+$] as (b), but replacing $\mu^+$ for some cardinal in $(\mu,2^\mu)$.
\end{enumerate}
\end{enumerate}
\mn

Our proof starts with a supercompact cardinal $\kappa$, and we show (as part of a more general axiomatic frame) that a stronger version of a universal on $\lambda > \kappa$ (e.g.\ $\lambda = \kappa^+$) is sufficient for the existence of a universal graph on $\lambda$ even after forcing with some $\bbP$ satisfying the axiomatic requirements. Then we first build a general frame for the preparation, and then construct the strong universal as in \cite{Sh:175a} suited to the present frame. (Here we remark that some large cardinal hypotheses are essential, as the failure of the Singular Cardinal Hypothesis itself implies that there is an inner model with the Mitchell order $o(\kappa) = \kappa^{++}$ for a measurable cardinal $\kappa$ (in fact these are equiconsistent).)

The paper is organized as follows. 
In \S\ref{1} we introduce the concept of $(\lambda,\kappa)-i$ $(i = 0,1$) systems, and in Claim $\ref{h8}$ we prove that extending a ground model already admitting some strong version of universal using such a $(\lambda,\kappa)-i$ system results in a model with the desired universal object. 
In \S\ref{2} we prove that Prikry forcing, Magidor forcing and Radin forcing give rise to a $(\lambda,\kappa)-1$ system provided the relevant filters satisfy some reasonable assumptions.
In \S\ref{3a} we prepare the ground, in Claim $\ref{c10}$ build the frame to force $(\lambda,\kappa)-1$ systems using a supercompact cardinal.
In \S\ref{3B} we construct a forcing for obtaining the strong universal fitting in the frame in Claim $\ref{c10}$.

In works in preparation we intend to replace graphs by more general classes; much of our for is not specific to graphs. Also for consistency of \ref{(*)+} for a small singular $\mu$, e.g. $\mu = \aleph_\omega = \beth_\omega$.

\bigskip

\subsection {Preliminaries}
\bigskip

		We are interested in universal objects in the class of graphs, i.e.\ models of the first order language admitting no functions, only a single symmetric, nonreflexive binary relation.
	Under ordinals we always mean von Neumann ordinals, and for a set $X$ the symbol $|X|$ always refers to the smallest ordinal with the same cardinality. If $f$ is a mapping with $\dom(f) \supseteq X$, then $f``X = \{ f(x): \ x \in X\}$, i.e.\ the pointwise image of $X$. For a set $X$ the symbol $\cP(X)$ denotes the power set of $X$, while if $\kappa$ is an ordinal we use the standard notation $[X]^\kappa$ for $\{Y \in \cP(X): \ |Y| = \kappa\}$, similarly for $[X]^{<\kappa}$, $[X]^{< \kappa}$, etc.\ By a sequence we mean a function on an ordinal, where for a sequence $ \overline{s}= \langle s_\alpha: \ \alpha < \dom(\overline{s}) \rangle$ the length of $\overline{s}$ (in symbols $\lh(\overline{s})$) denotes $\dom(\overline{s})$. Moreover, for sequences $\overline{s}$, $\overline{t}$ let $\overline{s} \tieconcat \overline{t}$ denote the natural concatenation (of length $\lh(\overline{s}) + \lh(\overline{t})$).
	For a set $X$, and ordinal $\alpha$ we use	$^\alpha X = \{ \overline{s}: \ \lh(\overline{s}) =  \alpha, \ \ran(\overline{s}) \subseteq X\}$, and for cardinals $\lambda$, $\kappa$ we use the symbol  $\lambda^\kappa = |^\kappa \lambda|$ (that is, the least ordinal equivalent to it).

	Regarding iterated forcing and quotient forcing we will mostly use the terminology of the survey \cite{Ba}.
	However we adhere to the following conventions.
\begin{convention}
	Regarding forcing we follow the convention that ``$p \leq q$" means that $q$ is the stronger, i.e.\ giving more information.
\end{convention}
\begin{convention}
	A notion of forcing $\bbP$ is $<$$\mu$-directed closed ($<$$\mu$-closed, resp.), if for any directed (increasing, resp.) system $\{p_\alpha: \ \alpha < \nu < \mu\}$  there exists a common upper bound $p_*$ in $\bbP$.
\end{convention}
A filter $\mathcal{F} \subseteq \cP(X)$ is $\kappa$-complete, if for each $\{F_\alpha: \ \alpha < \nu < \kappa\} \subseteq \mathcal{F}$ we have $\bigcap_{\alpha < \nu} F_\alpha \in \mathcal{F}$. A partial order $T$ is $\mu$-directed, if for each $\{t_\alpha: \ \alpha < \nu < \kappa\} \subseteq T$ there exists a common upper bound $t_* \in T$.

\section {The frame and deducing the consistency results} \label{1}
\bigskip

\subsection {What We Do} \label{0B}\
\bigskip

In the present paper we introduce a more general framework and apply it for the class of graphs.

We shall start with a large cardinal, like a Laver indestructible
supercompact, or with forcing a relative of it.  We then have a two
step forcing.

\underline{First}, a forcing $\bbP$ with the following three properties:
\mn
\begin{enumerate}
\item[(a)]  preserving the largeness of $\kappa$,
\sn
\item[(b)]  moreover, in $\bfV^{\bbP}$
\newline
  there is a normal
$\kappa$-complete filter on $\kappa$ such that $(D, {}^*\supseteq)$ is
$\lambda^+$-directed for a suitable cardinal $\lambda < 2^\kappa$,
\item[(c)]  preparing the ground for the results we like to have on
  $\lambda$, e.g. has a strong version of ``there is a universal graph
  in $\lambda,\lambda < 2^\kappa$".
\end{enumerate}
\mn
\underline{Second}, a forcing $\bbQ$ (in $\bfV^{\bbP}$) such that: 
\mn
\begin{enumerate}
\item[(d)]  makes $\kappa$ singular,
\sn
\item[(e)]  preserves $\kappa$ is strong limit and $2^\kappa$ large,
\sn
\item[(f)]  but to get the desired property of $\lambda$, we use
  $\bbQ$ that fits the frame in Definition \ref{h2} below,
\sn
\item[(g)]  then prove the existence of a universal object using the frame
\end{enumerate}
(or instead of $(\bfV^{\bbP})^{\bbQ}$ use
$\bfV^{\bbP}[\name X]$ for a $\bbQ$-name $\name X$).

Now in \S1, Definition \ref{h2} define the family of
$(\lambda,\kappa)$-systems fitting (f), then we deduce the
existence of universal graphs in $\lambda$ (a case of (g)).

In \S2 we shall
prove that classical forcings for making $\kappa$ singular fit our
frame, i.e. satisfy (d)-(g).

In \S3 we shall deal with finding $\bbP$ as in (a),(b),(c), so have to
combine the specific forcing (say forcing a universal graph in
$\lambda$, i.e. clause (c)) and guaranteeing the existence of e.g. a
normal ultrafilter of which is $\lambda^+$-complete in a suitable
sense (i.e. clause (b)).
\bigskip

\begin{definition}\label{sydf}{\ } \\
\label{h2}
1) We say $\bfr$ is a $(\lambda,\kappa)-1$-system \when \, $\bfr = (\bbR, \name X, \le_{\pr}, \cS)=(\bbR_\bfr, \name X_\bfr, \le_{\bfr, \pr}, \cS_\bfr)$ satisfies the following
\mn
\begin{enumerate}[label = $(\alph*)$, ref = $\alph*$]
\item  $\kappa$ is inaccessible,
\sn
\item  $\lambda \in [\kappa^+,2^\kappa)$,
\sn
\item  $\bbR$ is a forcing notion preserving ``$\kappa$ is strong limit",
\sn
\item  $\name X$ is an $\bbR$-name of a subset of $\kappa$,
\sn
\item  $\le_{\pr} \subseteq \le_{\bbR}$ is a quasi-order,
\sn
\item  for each $p \in \bbR$ we have $\cS_p \subseteq \{\bar q \in {}^\kappa
  \bbR:p \le_{\pr} q_\varp$ for every $\varp < \kappa\}$,
\sn
\item \label{syg} 
whenever $\gamma < \kappa$, $p \in \bbR$, $\name\tau$ are such that  $p \Vdash "\name\tau$ is an $\bbR$-name of a
  subset of $\gamma$ or just of $\cH^\bfV(\gamma)$" 
  \then \,:
\sn
\begin{enumerate}
\item[$(*)$]  there are $\bar q  \in \cS_p,\bar\gamma = \langle \gamma_\varp: \ \varp < \kappa \rangle \in {}^\kappa \kappa$ and
  $\bar c = \langle   c_\varp:\varp < \kappa\rangle$, where
\sn
\begin{enumerate}
\item[$\bullet_1$]   each $c_\varp \in \bfV$ is a code for  $|\cH^\bfV(\gamma_\varp)|$-Borel function $F_{c_\varp}$ from
$\cP(\cH^\bfV(\gamma_\varp))$ into $\cP(\cH^\bfV(\gamma))$,
\sn
\item[$\bullet_2$]   $q_\varp \Vdash ``\name\tau =
 F_{c_\varp}(\name X \cap \cH(\gamma_\varp))"$, and so it belongs to
$\cH(\kappa)$;
\end{enumerate}

\end{enumerate}
\sn
\item[(h)]  if $\bar q_\alpha \in \cS_p$ for $\alpha < \lambda$, \then
  \, for some $q_* \in \bbR$ for every $\alpha < \lambda$ there is $\varp_\alpha
  < \kappa$ such that $q_{\alpha,\varp_\alpha} \le_{\bbR} q_*$.
\end{enumerate}
\mn
2) We say $\bfr$ is a $(\lambda,\kappa)-2$-system \when \, above in clause
(g) we restrict ourselves to $\name\tau$'s such that $\Vdash ``\name\tau
\in \bfV[\name X]"$;,


\noindent
2A) We may omit the 1 in ``1-system".

\noindent
3) We say $\bfr$ is nice \when \, the
forcing $\bbR_{\bfr}$ collapses no cardinal.

\noindent
4) It is enough to require in clause $(g)$ that $(*)$ holds for every $p$, $\name \tau$ such that $p \Vdash ``\name \tau \in \{0,1\}"$, as this formally weaker assumption easily implies clause $(g)$. 
\end{definition}

\begin{discussion}{\ } \\
\label{h5}
\noindent 1) Here we only deal with the question ``when is there a universal
graph in the cardinal $\lambda$?".

\noindent
2) Of course, in Definition \ref{h2}, we are interested in the case
$\Vdash_{\bbR_{\bfr}} ``\kappa$ is singular".

\noindent
3) There are such $\bfr$'s: Prikry forcing, Magidor forcing, cases of
Radin forcing, see \ref{b2} on.

\noindent
name in a derived forcing.
\end{discussion}

\begin{claim}
	For a fixed $\iota \in \{0,1 \}$ and cardinals $\kappa$, $\lambda$
\label{h8} \begin{enumerate}[label = $\arabic*)$, ref = $\arabic*)$]
	\item   In $\bfV_\iota$ there is a universal graph of cardinality
$\lambda$ \when \,:
\mn
\begin{enumerate}
\item[(a)] $\bfr \in \bfV$ is $(\lambda,\kappa)-\iota$-system, we define $\bfV_\iota = \bfV^{\bbR_\bfr}$ if $\iota =1$, and  $\bfV_\iota = \bfV[X_\bfr]$ in case of $\iota = 2$,
 \sn
\item[(b)]  $\kappa < \lambda < 2^\kappa$ (e.g. $\lambda = \kappa^+$),
\sn
\item[(c)]  in $\bfV$, there is a universal member of
  $(K_{\kappa})_\lambda$, see below (Definition $\ref{h11}$).
\end{enumerate}
\mn
\item Moreover, $\univ(K_\lambda) \le \chi$ (where $K_\lambda$ denotes the class of graphs on $\lambda$-many vertices) is true in $\bfV_\iota$, when:

\mn
\begin{enumerate}
\item[(a),(b)]  as above
\sn
\item[(c)]  $\univ((K_{\kappa})_\lambda) \le \chi$ in $\bfV$.
\end{enumerate}
\end{enumerate}
\end{claim}

\begin{definition}
\label{h11}
$(K_{\kappa})_\lambda$ is the class of edge colored graphs with
$\kappa$ colors, equivalently $M \in (K_{\kappa})_\lambda$ \Iff \,
\mn
\begin{enumerate}
\item[(a)]  $M = (|M|,R^M_\varp)_{\varp < \kappa}$,
\sn
\item[(b)]  $\|M\| = \lambda$,
\sn
\item[(c)]  $R^M_\varp$ is a symmetric irreflexive two-place relation
  on $|M|$,
\sn
\item[(d)]  $\langle R^M_\varp:\varp < \kappa\rangle$ is a partition
  of $\{(a,b):a \ne b \in M\}$.
\end{enumerate}
\end{definition}

\begin{PROOF}{Claim \ref{h8}}(Claim \ref{h8})
 Let (in $\bfV$) $\langle (c_\varp, \gamma_\varp):\ \varp < \kappa\rangle$ list
 \[ \{ (c, \gamma): \ c: \cP(\cH^\bfV(\gamma)) \to \{0,1\}  \text{ is a code for an } |\cH^\bfV(\gamma)|-\text{Borel function}  \}, \]

Assume that
\begin{enumerate}
	\item[$(*)_1$] there is a sequence $ \langle M^*_\delta: \ \delta < \chi \rangle$ in $(K_{\kappa})_\lambda$ forming a universal sequence for $(K_{\kappa})_\lambda$ (in
the universe $\bfV$, of course; so $\chi = 1$ means that $M^*_0 \in (K_\kappa)_\lambda$ is universal ),
\end{enumerate}
where $M^*_\delta = (\lambda,\ldots,R^{M^*_\delta}_\varp,\ldots)_{\varp < \kappa}$,
it is enough to prove that in $\bfV_\iota$ we have $\univ(K_\lambda) \leq \chi$.

Now we define the sequence of $\bbR_{\bfr}$-names $\name M_\delta$ ($\delta < \chi$) for graphs as follows.
\mn
\begin{enumerate}
\item[$(*)_2$]
\begin{enumerate}
\item[(a)]  the set of nodes of $\name M_\delta$ is $\{\alpha: \ \alpha < \lambda\}$,
\sn
\item[(b)]  for $\alpha \ne \beta < \mu$ let $(\alpha,\beta) \in
  R^{\name M_\delta}$ \Iff \, for the unique $\varp < \kappa$ with $(\alpha,\beta) \in
  R^{M^*_\delta}_\varp$ we have $F_{c_\varp}(\name X \cap \cH(\gamma_\varp)) = 1$.
\end{enumerate}
\end{enumerate}
\mn
So clearly
\mn
\begin{enumerate}
\item[$(*)_3$] for each $\delta < \chi$ $\name M_\delta$ is an $\bbR_{\bfr}$-name of a graph with set
  of nodes $\lambda$.
\end{enumerate}
\mn
Hence it suffices to prove:
\mn
\begin{enumerate}
\item[$(*)_4$]
$\Vdash_{\bbR} \ `` \bfV_\iota \models \langle \name M_\delta: \ \delta < \chi \rangle$ is a universal sequence in $K_\lambda$".
\end{enumerate}
\mn
So why does $(*)_4$ hold?  Assume
\mn
\begin{enumerate}
\item[$(*)_{4.1}$]   $p \Vdash ``\name N \in \bfV_\iota$ is a graph with set of nodes
$\lambda"$.
\end{enumerate}
\mn
Let $\langle (\alpha_\gamma,\beta_\gamma):\gamma < \lambda \rangle \in \bfV$ list
the set of pairs $(\alpha,\beta)$ such that $\alpha < \beta <
\lambda$.   For each $\gamma < \lambda$ (considering the $\bbR_{\bfr}$-names $\name\tau_\gamma$ for the truth value of
$(\alpha_\gamma,\beta_\gamma) \in R^{\name N}$) clause (g) of
Definition \ref{h2} 1) gives  $\bar q_{\gamma} = \langle q_{\gamma,\varp}: \ \varp < \kappa \rangle \in
\cS_p$, $\bar\zeta_{\gamma} = \langle \zeta_{\gamma, \varp}: \ \varp < \kappa \rangle \in {}^\kappa \kappa$ and $\bar c_{\gamma}
= \langle c_{\gamma,\varp}:\varp < \kappa\rangle$ such that for each $\gamma < \lambda$
\mn
\begin{enumerate}
\item[$\bullet_1$]  $c_{\gamma,\varp}$ is a  code for a $|\cH(\zeta_{\gamma,\varp})|$-Borel function 
$\cP(\cH(\zeta_{\gamma,\varp})) \rightarrow \{0,1\}$ ($\varp < \kappa$)
\sn
\item[$\bullet_2$]  $q_{\gamma, \varp} \Vdash_{\bbR}
  ``(\alpha_\gamma,\beta_\gamma) \in R^{\name M} \Leftrightarrow
F_{c_{\gamma,\varp}}(\name X \cap \cH(\zeta_{\gamma,\varp}))=1$".
\end{enumerate}

\mn
Now by clause (h) of Definition \ref{h2} 1), there are $q$ and $\langle
\varp_\gamma = \varp(\gamma):\gamma < \lambda\rangle \in
{}^\lambda\kappa$ such that:
\mn
\begin{enumerate}
\item[$\bullet_3$]  $q$ is above $q_{\gamma,\varp(\gamma)}$ for every
$\gamma < \lambda$.
\end{enumerate}
\mn
Now we define $N_*$ as follows:
\mn
\begin{enumerate}
\item[$(*)_{4.3}$]
\begin{enumerate}
\item[(a)]  $N_* = (\lambda,(R^{N_*}_\varp)_{\varp < \kappa})$, where
\sn
\item[(b)]  $R^{N_*}_\varp =
  \{(\alpha_\gamma,\beta_\gamma):\gamma < \lambda$ and $\varp_\gamma 
= \varp\}$ ($\forall \varp \in \kappa$).
\end{enumerate}
\end{enumerate}
\mn
Clearly
\mn
\begin{enumerate}
\item[$(*)_{4.4}$]  $N_* \in (K_{\kappa})_\lambda$ (with set of nodes
  $\lambda$) belongs to $\bfV$.
\end{enumerate}
\mn
Now choose a suitable $\delta < \chi$ and a function $f$ so that:
\mn
\begin{enumerate}
\item[$(*)_{4.5}$]  $f: N^* \to M^*_{\delta}$ is an embedding, $f \in \bfV$ 
\end{enumerate}
\mn
[which exists by $(*)_1$.]
\mn
Finally it is straightforward to check that
\begin{enumerate}
\item[$(*)_{4.6}$]  $q \Vdash ``f$ is an embedding of $\name N$ into
 $\name M_\delta"$.
\end{enumerate}
\mn
[Recall how we defined $N_*$ from $N$ $(*)_{4.3}$, $M_\delta$ from $M_\delta^*$ $(*)_2$, and the choice of $f$ $(*)_{4.5}$.]

\end{PROOF}

\noindent
Naturally we ask
\begin{question}{\ } \\
\label{h14}
1) What about $(K_{\kappa})_\lambda$?

\noindent
2) More seriously about the theory of triangle free graphs, or of
$T_{\feq}$ (equivalently $T_\ceq$, see \cite{Sh:1164}). \underline{Note}, on $T_\textrm{feq}$ see \cite{Sh:457}, or \cite{Sh:614}, and on the non-existence in case of $T_{\ceq}$ see \cite{Sh:1064}.

\noindent
3) Moreover, $(\Mod_T,\prec),T$ simple?  Or even $\NSOP_2$?  (of
cardinality $< \kappa$).  We have to be more careful because of,
e.g. function symbols.
\end{question}

\noindent
A work in preparation deals with \ref{h14} 2), 3).
Concerning \ref{h14} 1) (note that this does not reflect on Claim $\ref{h8}$):
\begin{claim}
\label{h17}
Assume $\kappa$ is strong limit singular and $\kappa < \lambda <
2^\kappa$.  \Then \, in $(K_{\kappa})_\lambda$ there is no universal member.
\end{claim}
\begin{remark}
	It suffices to have $\beth_\omega(\cf(\kappa))< \kappa$, and $(\alpha < \kappa \ \to \ |\alpha|^{\cf(\kappa)} < \kappa)$.
\end{remark}
\begin{PROOF}{\ref{h17}}
By  \cite[Thm 1.13 and 1.14 (2) on RGCH]{Sh:829}

\begin{enumerate}
	\item[$(*)_0$]
 there\footnote{In fact, instead
  ``$\kappa$ strong limit singular $\sigma$ as above", it suffices to assume less.}
  is a regular $\sigma \in (\cf(\kappa),\kappa)$ such that $\lambda^{[\sigma, \kappa]} = \lambda$, i.e.\
  there is $\cP' \subseteq \{u \subseteq \lambda:|u| \leq
  \kappa\}$ of cardinality $\lambda$ such that every $u \subseteq
  \lambda$ of cardinality $\leq \kappa$ is the union $< \sigma$ members of $\cP'$.
  
\end{enumerate}
 Therefore, as $\sigma = \cf(\sigma) > \cf(\kappa)$
\mn
\begin{enumerate}
\item[$(*)_1$]  there is $\cP \subseteq \{u \subseteq \lambda:|u| <
\kappa\}$ of cardinality $\lambda$ such that every $u \subseteq
\lambda$ of cardinality $< \kappa$ is the union $< \sigma$ members of $\cP$.
\end{enumerate}
\mn
Let $M_* \in (K_{\kappa})_\lambda$ and we shall prove that it is not
universal; \wilog \, the universe of $M_*$ is $\lambda$.
Now for each $u \in \cP$ and $\alpha < \lambda$ let
\[ v(\alpha,u,M_*) = \{\varp < \kappa: \text{ for some }\beta \in u
 \text{ we have } (\alpha,\beta) \in R^{M_*}_\varp\}, \]
  so $v(\alpha,u,M_*)
\subseteq \kappa$ has cardinality $< \kappa$.
Let 
\[ \cP_1 = \{w \in [v(\alpha,u,M_*)]^{\leq \cf(\kappa)}: \ u \in \cP,  \alpha \in \lambda \}, \] so
\mn
\begin{enumerate}
\item[$(*)_2$]   $\cP_1 \subseteq [\kappa]^{\leq \cf(\kappa)}$.
\end{enumerate}
\mn
Now
\mn
\begin{enumerate}
\item[$(*)_3$]   $|\cP_1| \le |\cP| + 2^{<\kappa} \le \lambda
< 2^\kappa = \kappa^{\cf(\kappa)}$.
\end{enumerate}
\mn
Hence
\mn
\begin{enumerate}
\item[$(*)_4$]   we can find $v \subseteq \kappa$ of cardinality $\cf(\kappa)$
which is not in $\cP_1$, moreover, $u \in \cP_1
\Rightarrow |u \cap v| < \cf(\kappa)$,
\end{enumerate}
\mn
which is justified by the following argument: Let $\langle v_\gamma:\gamma < 2^\kappa\rangle$ be a sequence
of members of $[\kappa]^{\cf(\kappa)}$ with any two having intersection of
cardinality $<\cf(\kappa)$, hence for every $u \in \cP_1,\{\gamma <
2^\kappa:|u \cap v_\gamma| = \cf(\kappa)\}$ has cardinality $\le
2^{\cf(\kappa)} < \kappa$, so all but $\le \lambda$ sets of $v_\gamma$'s
are as required.

Now consider the following $N$:
\mn
\begin{enumerate}
\item[$(*)_5$]
\begin{enumerate}
\item[(a)]  $N = (A \cup B,\ldots,R^N_\varp,\ldots)_{\varp <
  \kappa}$ belongs to $(K_{\kappa})_{\sigma^{\cf(\kappa)}}$, where  $|A| = \sigma$, $|B| = \sigma^{\cf(\kappa)}$, $A \cap B = \emptyset$, 
\sn
\item[(b)]  $R^N_\varp \ne \emptyset$ \Iff \, $\varp \in v$,
\sn
\item[(c)]  letting $\langle \varp_i:i < \cf(\kappa)\rangle$ list $v$ (from $(*)_4$),
  for every sequence $\langle \alpha_i:i < \cf(\kappa)\rangle$ in $A$ with no repetitions
  there is $\beta \in B$ such that $(\alpha_i,\beta) \in R^N_{\varp_i}$ for
$i < \cf(\kappa)$.
\end{enumerate}
\end{enumerate}
\mn
Now if $g$ embeds $N$ into $M_*$ then for some $\{ u_\varp: \varp < \partial < \sigma\} \subseteq  \cP$,
 we have $\Rang(g \upharpoonright A) = \cup\{u_\varp:\varp <
\partial\}$.  Now as $|A| = \sigma = \cf(\sigma)$, there is $\varp <
\partial$ such that $|u_\varp \cap \Rang(g \upharpoonright A)| \ge \sigma \ge \cf(\kappa)$ so we can choose pairwise distinct $\alpha_i \in A$ ($i < \cf(\kappa)$) such
that $\{g(\alpha_i):i < \cf(\kappa)\}\subseteq u_\varp$.
Let $\beta \in B$ as in $(*)_5(c)$.  So
$g(\beta)$ is well defined and we get an easy contradiction by $(*)_4$.

So $N$ cannot be embedded into $M_*$, hence we are done.
\end{PROOF}
\bigskip

\section {Proving known forcings fit the frame} \label{2}
\bigskip

\subsection {Near a Large Singular}\
\bigskip

Here we do not collapse cardinals, just change cofinalities.

\begin{claim}
\label{b2}
There is a nice $(\lambda,\kappa)$-system $\bfr$ such that $\bbR_{\bfr} =
\bbP$ \when \,:
\mn
\begin{enumerate}
\item[(A)]
\begin{enumerate}
\item[(a)]  $\kappa < \lambda < 2^\kappa$ are cardinals,
\sn
\item[(b)]  $D$ is a normal ultrafilter on $\kappa$,
\sn
\item[(c)]  if $\cA \subseteq D$ has cardinality $\le \lambda$, \then
  \, for some $B \in D$ we have $(\forall A \in \cA)(B \subseteq
  A \mod [\kappa]^{< \kappa})$, (e.g. $D$ is generated by a
$\subseteq^*_\kappa$-decreasing sequence of length of a regular
  cardinal $> \lambda$),
\sn
\item[(d)]  $\bbP$ is Prikry forcing for $D$ (so change the cofinality
  of $\kappa$ to $\aleph_0$ and add no bounded subset of $\kappa$ and
  satisfies the $\kappa^+$-c.c).
\end{enumerate}
\end{enumerate}
\end{claim}

\begin{PROOF}{\ref{b2}}
Recalling the definition of Prikry forcing for $D$:
\mn
\begin{enumerate}
\item[$(*)_1$]
\begin{enumerate}
\item[(a)]  $p \in \bbP$ iff $p = (w,A) = (w_p,A_p)$, where $w_p \in
  [\kappa]^{< \aleph_0}$ and $A_p \in D$ and $[0,\max w_p] \cap A = \emptyset$,
\sn
\item[(b)]  $p \le_{\bbP} q$ \Iff \, $w_p \subseteq w_q
  \subseteq w_p \cup A_p$ and $A_p \supseteq A_q$.
\end{enumerate}
\end{enumerate}
\mn
We define the system $\bfr$ by letting:
\mn
\begin{enumerate}
\item[$(*)_2$]
\begin{enumerate}
\item[(a)]  $\kappa_{\bfr} = \kappa$,
\sn
\item[(b)]  $\lambda_{\bfr} = \lambda$,
\sn
\item[(c)]  $\bbR_{\bfr} = \bbP$,
\sn
\item[(d)]  $\name X_{\bfr} =$ the generic $= \cup\{w_p:p \in \bfG_{\bbP}\}$,
\sn
\item[(e)]  $\le_{\pr} =\le_{\bfr,\pr}$ is defined by $p \le_{\pr} q$
  iff $w_p = w_q \wedge A_p \supseteq A_q$ (and $p,q \in \bbR_{\bfr}$),
\sn
\item[(f)]  for $p \in \bbR_{\bfr} = \bbP$ let $\cS_p = \cS_{\bfr,p} :=
  \{\bar q:\bar q = \langle q_\varp:\varp < \kappa\rangle$ and for
  some $B \in D$ we have $B \subseteq A_p$ and
  $\{A_{q_{\varp}}:\varp < \kappa\}$ list $\{A:A \subseteq A_p$
  and $A=B \mod [\kappa]^{< \kappa}\}$.
\end{enumerate}
\end{enumerate}
\mn
We still have to prove that $\bfr$ is as required, namely, that $\bfr$ satisfies conditions listed in Definition $\ref{h2}$ 1).

Now in Definition \ref{h2} 1), clauses (a)-(f) hold trivially.  For
clause (g) recall that by  4) from Definition $\ref{sydf}$  it is enough to check for only $p$, $\name\tau$, where $``\tau \in \{0,1\}"$ is forced by $p$. Recall the following well-known fact:
\mn
\begin{enumerate}
\item[$(*)_3$]  if $p \in \bbP$, $p \Vdash_{\bbP}
  ``\name\tau \in \{0,1\}"$ \then \, for some $A' \seq A_p$ from $D$ we have:
\sn
 if $\alpha \in \kappa$ and $u \subseteq A_p \cap
  \alpha$ is finite then $(w_p \cup u,A' \backslash \alpha)$ forces a value for
  $\name\tau$.

\end{enumerate}
\mn
[For the sake of completeness we prove $(*)_3$: by the Prikry-lemma, for each $s \in [A_p ]^{<\aleph_0}$ there exists $A_s \seq A_p \se ((\max s) +1)$, $A_s \in D$, such that $(w \cup s, A_s)$ decides the value of $\name \tau$. Now let $A'$ be the diagonal intersection of $A_s$'s ($s \in [A_p]^{<\aleph_0}$), pedantically $\Delta_{\alpha < \kappa} (\bigcap_{s \in [\alpha+1]^{<\aleph_0}} A_s)$, it is straightforward to check that $A'$ works.]

So given $p \in \bbP,\gamma$ and $\name\tau$ as in clause (g) from Definition $\ref{h2}$, let $A' \subseteq A_p$
be as in $(*)_3$ and let $\bar q = \langle q_\varp:\varp <
\kappa\rangle$ be defined by: $q_\varp \in \bbP,w_{q_\varp} = w_p$ and
$\{A_{q_\varp}:\varp < \kappa\}$ list $\{A \subseteq A_p:A \equiv A'
\mod [\kappa]^{< \kappa}\}$.

We still have to choose the $\gamma_\varp,F_\varp$. For each $\varp$
choose $\zeta_\varp \in A_{q_\varp}$ such that $A_{q_\varp} \backslash \zeta_\varp =
A' \backslash \zeta_\varp$. Clause $(*)_3$ ensures that there is a function $f : [A_p \cap \zeta_\varp]^{<\aleph_0} \to \{0,1\}$ in $\bfV$ such that $q_\varp \Vdash \name \tau = f( \name X \cap \zeta_\varp)$.

Lastly, for clause (h), assume $p \in \bbR_{\bfr} = \bbP$ and $\bar q =
\langle \bar q_\alpha:\alpha < \lambda\rangle$ satisfies $\bar q_\alpha \in
\cS_p$.  So for each $\alpha < \lambda$ there exists $B_\alpha \seq A_p$ such that $\{A_{q_{\alpha,\varp}}:\varp < \kappa\}$ lists $\{ A \in D: A \seq A_p, \ A \equiv B_\alpha \mod [\kappa]^{<\kappa}\}$, hence by clause (A)(c) of the assumption of the claim,
there is $B \in D$, a subset of $A_p$ such that
$B \subseteq B_\alpha \mod
[\kappa]^{< \kappa}$ for each $\alpha \in \lambda$ and
let $q_*= (w_p,B)$ so clearly $p \le_{\pr} q_*$.  Also for each $\alpha < \lambda$, for
some $\zeta < \kappa$ we have
$B \backslash \zeta \subseteq B_\alpha$ hence because $\bar{q}_\alpha \in \cS_p$ for some $\varp < \kappa$ we have
$A_{q_{\alpha,\varp}} = (B_\alpha \backslash \zeta) \cup (A_p \cap
\zeta) \supseteq B$ hence $q_{\alpha,\varp} \le q_*$.

We still have to prove that $\bfr$ is nice but as $\bbP$ satisfies
the $\kappa^+$-c.c., and the Prikry lemma this is obvious.
\end{PROOF}

\begin{claim}
\label{b23}
There is a $(\lambda,\kappa)-1$-system $\bbR_\bfr$ with $\bfV^{\bbR_\bfr} \models \cf(\kappa) = \theta$ \when \,:
\mn
\begin{enumerate}[label = (\Alph*), ref = \Alph*]
\stepcounter{enumi}
\item \label{mgd}
\begin{enumerate}[label = (\alph*), ref = \alph*]
\item   $\theta = \cf(\theta) < \theta_* < \kappa < \lambda < 2^\kappa$,
\sn
\item \label{minc}  $\bar D = \langle D_i:i < \theta\rangle$ is a sequence
  of normal ultrafilters on $\kappa$, increasing in Mitchell order,
i.e. $i < j \Rightarrow D_i \in \Mos \Col(^\kappa \bfV/D_j)$,
\sn
\item \label{ldir} each $D_i$ ($i \leq \theta$) is $<\lambda^+$-directed mod $[\kappa]^{<\kappa}$, i.e. satisfies the condition in Claim $\ref{b2}$(A)(c),
\sn
\item the forcing  $\bbR_\bfr$ changes  the cofinality of $\kappa$ to $\theta$, preserves each cardinal and
the function $\mu \mapsto 2^\mu$, satisfies the $\kappa^+$-c.c. Moreover, we can prescribe that in $\bfV^\bbP$ there is no new subset of $\theta_*$.
\end{enumerate}
\end{enumerate}
\end{claim}

\begin{PROOF}{\ref{b23}}
			Using \cite[Proposition 2.1]{Krue}, condition \ref{minc} implies the following.
			\bscl If $\bar D = \langle D_i:i < \theta\rangle$ is an increasing (w.r.t.\ the Mitchell order) sequence
			of normal ultrafilters on $\kappa$, $\theta \leq \kappa$, then there exists a coherent sequence $\langle \bar{U}_\varp: \ \varp < \kappa+1 \rangle$, $\bar{U}_\varp = \langle U_{\varp}(\alpha): \ \alpha < o^U(\varp) \rangle$ for some function  $o^U: \kappa+1 \to \kappa$ such that $\bar{D} = \bar{U}_\kappa$, which means:
			\sn
			\begin{enumerate}[label = $(\intercal)_\alph*$, ref = $(\intercal)_\alph*$]
				\item  for each $\varp \leq \kappa$, $\alpha < o^U(\varp)$  $U_{\varp}(\alpha)$ is an $\varp$-complete normal ultrafilter on $\varp$,
				\item \label{cos} moreover, for each $\varp \leq \kappa$, $\alpha < o^U(\varp)$ letting $\bfj_{\varp,\alpha}: \bfV \to \Mos\Col(^\varp \bfV / U_{\varp,\alpha})$ be the associated elementary embedding, we have
				\[ (\bfj_{\varp, \alpha} (\bar{U} \um \varp))_\varp = \langle U_{\varp}(\beta): \ \beta < \alpha \rangle, \]
				\item $\langle U_{\kappa}(\alpha): \ \alpha < o^U(\kappa) \rangle = \langle D_\alpha: \ \alpha < \theta\rangle$.
			\end{enumerate}
			\escl
			
	Now we define the forcing $\bbP_{\overline{U}}$ to be the Magidor forcing associated to the sequence $\bar{D} = \bar{U}_\kappa = \langle U_{\kappa}(\alpha): \alpha \leq \theta \rangle$ , (see also \cite{Mg78}, or \cite{Gi10}), here we use the definition from \cite[Definition 5.22]{Gi10}
	\bdd \label{mdf}
	Define $\bbP_{\overline{U}}$ to be the following (auxiliary) poset.
	 \newcounter{sMa}
	\setcounter{sMa}{0}
	\begin{enumerate}[label = $(*_{\arabic*})$, ref = $*_{\arabic*}$]
			\setcounter{enumi}{\value{sMa}}
			
		\item
		 	 Let $p = \left\langle d_0, d_1, \dots, d_n, d_{n+1} = \langle \kappa, A_\kappa \right\rangle \rangle    \in \bbP_{\overline{U}}$, iff
		\begin{enumerate}
			\item $A_\kappa \in \bigcap \bar{U}_\kappa = \bigcap_{\alpha < \theta} U_{\kappa,\alpha}$,
			\item each $d_j$ ($j \leq n$) is of the form 
			\begin{itemize}
				\item either $\langle \varp, A_\varp \rangle$ for some $\varp < \kappa$, where $o^{U}(\varp) > 0$, moreover,
				$$A_\varp \in \bigcap \bar{U}_{\varp} = \bigcap_{\gamma < o^{U}(\varp)} U_{\varp,\gamma},$$ 
				(this case we define $\kappa(d_j) = \varp$),
				\item or $d_j = \varp$, when $o^{U}(\varp) = 0$ (and we let $\kappa(d_j) = \varp$).
			\end{itemize}
			\item $\kappa(d_0) < \kappa(d_1) < \dots < \kappa(d_n) < \kappa(d_{n+1}) = \kappa$, 
			\item moreover, for each $j \leq n$ if $d_{j+1}$ is a pair, then  $\kappa(d_j) < \min A_{\kappa(d_{j+1})}$.
		\end{enumerate}
		\stepcounter{sMa} \stepcounter{sMa}
		
		\item
	 
	 We define 
	 $$p = \langle d_0, d_1, \dots, d_n, d_{n+1} = \langle \kappa, A_\kappa \rangle\rangle \leq q = \langle e_0, e_1, \dots, e_n, e_{m+1} = \langle \kappa, B_\kappa \rangle \rangle,$$ if
	\begin{enumerate}
		\item $m \geq n$, and
		\item there exists a sequence $0 \leq i_0 < i_1 < \dots < i_n < j_{n+1} =m$ such that for each $j \leq n+1$ we have
	 	\begin{itemize}
	 		\item $\kappa(d_j) = \kappa(e_{i_j})$, and
	 		\item $B_{\kappa(d_j)} \seq A_{\kappa(d_j)}$,
	 	\end{itemize}
	\item	moreover, for each $k \leq m$ not of the form $i_j$ ($j \leq n+1$), if $i_l = \min \{ i_j: \ j \leq n+1, \ i_j > k\}$, then 
	 $$B_{\kappa(e_k)} \cup \{ \kappa(e_k) \} \seq A_{\kappa(d_l)}. $$
	 \end{enumerate}
 \item \label{Magfor} Now there are pairwise disjoint sets $Y_\alpha$ ($\alpha <\theta$) by $\delta \in Y_\alpha$ iff $o^U(\delta) = \alpha$, and 
 $$\{ p \in \bbP_{\overline{U}}: \ p \geq \langle \langle \kappa, \bigcup_{\alpha < \theta} Y_\alpha \rangle \rangle \}$$ is the Magidor forcing  changing the cofinality of $\kappa$ to $\min\{ \omega,\cf(\theta) \}$.
 \end{enumerate}
	
 	\edd

 	\bdd We define $p \leq_* q$ to be true iff $p \leq q$ and  $\lh(p) = \lh(q)$.
 	\edd

 	We define the system $\bfr$ by letting:
 	\mn
 
 \begin{enumerate}[label = $(*_{\arabic*})$, ref = $*_{\arabic*}$]
 	\setcounter{enumi}{\value{sMa}}
 	\item \label{rM}
 		\begin{enumerate}[label = (\alph*), ref = \alph*]
 			\item  $\kappa_{\bfr} = \kappa$,
 			\sn
 			\item  $\lambda_{\bfr} = \lambda$,
 			\sn
 			\item  $\bbR_{\bfr} = \{ p \in \bbP_{\overline{U}}: \ p \geq \langle \langle \kappa, \bigcup_{\alpha < \theta} Y_\alpha \rangle \rangle \}$,
 			\sn
 			\item  let $\name X_{\bfr}$ be  the generic sequence, i.e. 
 			$$\name X_\bfr = \cup\{ \{ \kappa(d_j): j < \lh(p) \} : \ p = \langle d_0, d_1, \dots, d_{\lh(p)-1} \rangle  \in \bfG_{\bbP}\} \se \{ \kappa \},$$
 			\item  $\le_{\pr} =\le_{\bfr,\pr}$ is defined by $p \le_{\pr} q$
 			iff $p \leq_* q$,
 			\sn
 			\item \label{sF}  for $p = \langle d_0, d_1, \dots, d_n, d_{n+1} = \langle \kappa, A_{p,\kappa} \rangle \rangle \in \bbR_{\bfr} = \bbP$ let 
 			$$\cS_p = \cS_{\bfr,p} := \left\{ \begin{array}{l}
 			\bar q:\bar q = \langle q_\varp:\varp < \kappa\rangle \text{, where} \\
 			 (\bullet_1) q_\varp = \langle d_0, d_1, \dots, d_{n}, \langle \kappa, A_{q_\varp,\kappa} \rangle\rangle, \text{ and} \\
 			\text{for some } B \in \bigcap \bar{U}_{\kappa}  \text{ we have } \\
 			 (\bullet_2) \ B \subseteq A_{p,\kappa},  \text{ and} \\
 			 (\bullet_3) \  \{A_{q_{\varp}, \kappa}:\varp < \kappa\} \text{  lists }\{A_*:A_* \subseteq A_{p,\kappa} \ \wedge \ A_*=B \mod [\kappa]^{< \kappa}\} \end{array} \right\}. $$
 		\end{enumerate}
	\stepcounter{sMa} 
	\end{enumerate}
	It is known that $\name X$ is a club of $\kappa$ of order type $\theta$, moreover, if condition $\langle \langle \beta \rangle, \langle \kappa, A \rangle \rangle$ is in the generic filter (for some $\beta < \kappa$, $o{^U(\beta)} = 0$, then the forcing adds  no new subset of $\beta$. This means that by \ref{cos} $\{ \beta < \kappa: \ o^U(\beta) = 0\} \in U_{\kappa,0}$ is of cardinality $\kappa$, so there is no problem assuming that $\langle \langle \beta \rangle, \langle \kappa, A \rangle \rangle \in \bfG$ for some $\beta \geq \theta_*$.
	In order to finish the proof of Claim $\ref{b23}$ it suffices to verify that the forcing defined in \ref{Magfor} is a $(\lambda,\kappa)-1$-system.
	\bscl \label{bbPU} 
		If $\langle \bar{U}_\varp: \ \varp < \kappa+1 \rangle$ is a coherent sequence, where the ultrafilters $\{ U_\kappa(\alpha): \ o^{\bar{U}}(\kappa) \}$ are $< \lambda^+$-directed mod $[\kappa]^{<\kappa}$, then the forcing $\bbP_{\bar{U}}$ from Definition $\ref{mdf}$ is a $(\lambda,\kappa)-1$-system.
	\escl
	\begin{PROOF}{Subclaim \ref{bbPU}}
	
	Now we have only to check the requirements of Definition $\ref{h2}$ 1).
	Recall the following properties of the Magidor forcing, see \cite[Sec. 5.1 and 5.2]{Gi10}.
	\bfact (Prikry Lemma)
		For each $p \in \bbP_{\overline{U}}$ and each formula $\sigma(\name x_0, \dots, \name x_m)$ there exists $q \geq_* p$, $q \parallel \sigma(\name x_0, \dots, \name x_m)$ (i.e. either $q \Vdash \sigma(\name x_0, \dots, \name x_m)$, or $ q \Vdash \neg \sigma(\name x_0, \dots, \name x_m)$).
	\efact
	
	\bfact \label{reszgen} Suppose that $\bfG\subseteq \bbP_{\overline{U}}$ is generic over $\bfV$, $p = \langle d_0, d_1, \dots, d_n, d_{n+1} = \langle \kappa, A_{p,\kappa} \rangle \rangle \in \bfG$, $d_i = \langle \kappa(d_i), A_{\kappa(d_i)}$, then the filter 
	$\bfG \um (\kappa(d_i)+1) := \{q \um (\kappa(d_i) +1): \ q \in \bfG \}$ is $\bfV$-generic over the Prikry forcing $\bbP_{\overline{U} \um (\kappa(d_i)+1)}$ associated to the coherent sequence $\langle \overline{U}_\delta = \langle U_\delta(\gamma): \ \gamma <o^U(\delta) \rangle : \ \delta \leq \kappa(d_i) \rangle$.
	\efact
	These two fact clearly imply the following.
	\bl \label{torfor} 
		For each formula $\sigma(\name x_0, \dots, \name x_m)$ and $\delta$ with  for some $p = \langle d_0, d_1, \dots, d_n, d_{n+1} \rangle \in \bfG$ and $m \leq n$ $\delta = \kappa(d_m)$ (i.e.\ $\delta \in X_\mathfrak{r}$) there exists $q \in \bbP_{\overline{U} \um (\delta+1)} \geq p \um (\delta+1)$, such that $q = q^* \um (\delta+1)$ for some $q^* \in \bfG$, and for some $A^* \in \bigcap \overline{U}_\kappa$
		\beeq \label{qqcsi} q \tieconcat \langle \kappa, A^* \rangle \parallel_{\bbP_{\overline{U}}} \sigma(\name x_0, \dots, \name x_m).\eeq

	\el
	\begin{PROOF}{}
		By Fact $\ref{reszgen}$ it is enough to show that there are densely many such $q$'s in 
		$\bbP_{\overline{U} \um (\delta+1)}$. But by the Prikry Lemma for $\bbP_{\overline{U} \um (\delta+1)}$ there are in fact $\geq_*$-densely many such $q's$.
		\end{PROOF}
	
	Similarly to the case of Prikry forcing, this has the following consequence.
 	\bcl \label{kisebbk}
 		For each $p = \langle d_0, d_1, \dots, d_n, \langle \kappa, A_{p,\kappa} \rangle \in \bbP$, $\name \tau$ (with $p \Vdash \name \tau \in \{0,1 \}$) there exists a set $A' \in \bigcap \bar{U}_{\kappa}$, $A' \subseteq A_{p,\kappa}$, such that whenever $q = \langle e_0,e_1, \dots, e_m, \langle \kappa, A_{q,\kappa} \rangle \rangle \geq p' = \langle d_0, d_1, \dots, d_n, \langle \kappa, A' \rangle \rangle$, $\alpha \in A_{p,\kappa}$ are given with $\kappa(e_m) \leq \alpha$, and $q$ forces  a value for $\name \tau$, then so does $$q' = \left\langle e_0,e_1, \dots, e_m, \langle \kappa, A' \cap (\alpha, \kappa) \rangle \right\rangle.$$
 	\ecl
 	\begin{PROOF}{} For each $\alpha \in A_{p,\kappa}$ define $B_\alpha \subseteq A_{p,\kappa}$ so that
 	whenever $q = \langle e_0,e_1, \dots, e_{m} = \langle \kappa, A_{q,\kappa} \rangle \rangle  \geq p$ (with $\kappa(e_0),\kappa(e_1), \dots, \kappa(e_m) \leq \alpha$) decides the value of $\name \tau$, then so does $q' = \langle e_0,e_1, \dots, e_{m+1} = \langle \kappa,  B_\alpha  \rangle \rangle$. This can be done easily: first for each possible $e_0,e_1, \dots, e_{m}$ choose a set $B_{e_0,e_1, \dots, e_{n}} \seq (\alpha, \kappa)$ with $\langle e_0,e_1, \dots, e_{n}, \langle \kappa, B_{e_0,e_1, \dots, e_{n}} \rangle \rangle$ deciding the value of $\name \tau$ if such a $B_{e_0,e_1, \dots, e_{n}}$ exists, otherwise just let $B_{e_0,e_1, \dots, e_{n}} = A_{p,\kappa} \cap (\alpha, \kappa)$. Second, let $B_\alpha = \bigcap_{e_0,e_1, \dots, e_{n}} B_{e_0,e_1, \dots, e_{n}}$.
 	Now it is easy to check that the diagonal intersection  $A' = \Delta_{\alpha \in A_{p,\kappa}} B_\alpha \in \bigcap \bar {U}_\kappa$ works (note that the intersection of normal measures is a normal filter).
 	\end{PROOF}
 	
 	\bcl \label{kul} For every $p  \in \bbP$ and $\name\tau$, if $p \Vdash \name \tau \in \{0,1 \}$, then 
 	 we can choose $\bar{q} = \langle q_\varp: \ \varp < \kappa \rangle \in \cP_p$, $\langle \gamma_\varp: \ \varp< \kappa \rangle$,
 	  $\langle c_\varp: \ \varp < \kappa \rangle$ where each $c_\varp$ is a code for a $\gamma_\varp$-Borel function from $ \cP(\cH^\bfV(\gamma_\varp))$ to $\{0,1 \}$ so that
 	  \[ q_\varp \Vdash \ \name\tau = f_{c_\varp}(\name X \cap \gamma_\varp). \]
 	 \ecl
 	\begin{PROOF}{Claim \ref{kul}}
 	 First if $p = \langle d_0,d_1, \dots, d_n, \langle \kappa, A_{p,\kappa} \rangle\rangle \in \bbP$, $\name \tau$  are as in clause (g) from Definition $\ref{h2}$, let $A' = A'(p,\name\tau) \subseteq A_{p,\kappa}$
 	be given by Claim $\ref{kisebbk}$ and 
 	 \begin{enumerate}[label = $(*_{\arabic*})$, ref = $*_{\arabic*}$]
 		\setcounter{enumi}{\value{sMa}}
 		\item \label{qv} let $\bar q = \langle q_\varp:\varp <
 		\kappa\rangle \in \cS_p$ be defined by: $q_\varp \in \bbP, q_\varp =\langle d_0,d_1, \dots, d_n, \langle \kappa, A_{q_\varp,\kappa} \rangle\rangle$  where
 		$\{A_{q_\varp, \kappa}:\varp < \kappa\}$ lists $\{A_* \subseteq A_{p,\kappa}: A_* \equiv A'
 		\mod [\kappa]^{< \kappa}\}$.
 		
 	\end{enumerate}
 	\stepcounter{sMa} 
 
 	We still have to choose the $\gamma_\varp,c_\varp$. For each $\varp$
 	choose $\zeta_\varp \in A_{q_\varp,\kappa}$ such that $A_{q_\varp,\kappa} \cap (\zeta_\varp, \kappa) =
 	A' \cap (\zeta_\varp, \kappa)$. Now this with Claim $\ref{kisebbk}$ imply that $q_\varp$ forces that $\name \tau$ only depends on $\bfG \um (\zeta_\varp+1)$, in the following way.
 	\bscl \label{kuls}
 		If $q_\varp \in \bfG$, then for some $q^* \in \bfG$, $\delta \leq \zeta_\varp$
 		$$ q^* \um (\delta+1) \tieconcat \langle \kappa, A_{q_\varp,\kappa} \setminus [0,\zeta_\varp] \rangle \parallel ``\name \tau = 1 ".$$
 	\escl
 	\begin{PROOF}{Subclaim \ref{kuls}}
 	Indeed, if $q_\varp \in \bfG$, then by genericity there is some $\delta \leq \zeta_\varp$, and $q'= \langle e_0,e_1, \dots, e_m \rangle \geq q_\varp$, $q' \in \bfG$, and $\kappa(e_k) = \delta$, $A_{q', \kappa(e_{k+1})} \cap [0, \zeta_\varp] = \emptyset$ (i.e. $q'$ forces $\max(\name X \cap [0,\zeta_\varp]) = \delta)$), w.l.o.g.\ $q' \geq q_\varp$. Now  by Lemma $\ref{torfor}$ there is some $q^* \in \bfG$, $A^* \in \bigcap \overline{U}_\kappa$ with \beeq \label{foegy} q^* \um (\delta+1) \tieconcat \langle \kappa, A^* \rangle \parallel ``\name \tau = 1", \eeq 
 	w.l.o.g.\ $q^* \geq q' \geq q_\varp$. But then by the construction of $A' = A(p,\name \tau)$
 	we have 
 	\beeq \label{foeggy} q^* \um (\delta+1) \tieconcat \langle \kappa, A' \setminus [0,\delta] \rangle \parallel ``\name \tau = 1", \eeq 
 	therefore as $A_{q_\varp,\kappa} \cap (\zeta_\varp, \kappa) = A' \cap (\zeta_\varp, \kappa)$ (and $\delta \leq \zeta_\varp$), 
 	\beeq \label{foegggy} q^* \um (\delta+1) \tieconcat \langle \kappa, A' \setminus [0,\delta] \rangle \leq q^* \um (\delta+1) \tieconcat \langle \kappa, A_{\zeta_\varp} \setminus [0,\zeta_\varp] \rangle \parallel ``\name \tau = 1", \eeq 
 	as desired.
 \end{PROOF}
 	It is not difficult to check that this  implies that for every $q^{**} = \langle e_0,e_1, \dots, e_m \rangle \in \bigcup_{\delta \leq \zeta_\varp} \bbP_{\overline{U} \um (\delta+1)}$ with  
 	\[ q^{**} \tieconcat  \langle \kappa, A_{\zeta_\varp} \setminus [0,\zeta_\varp] \rangle \Vdash \name \tau = j_{q^{**}},  \] 
 	we have 
 	\[  \begin{array}{rl} q_\varp  \Vdash & \left(\{\kappa(e_i):  i \leq m \} \seq \name X \cap [0,\zeta_\varp] \seq \{\kappa(e_i): \ i \leq m \} \cup (\cup\{A_{q^{**},\kappa(e_i)}: \ i \leq m \})\right) \\ &  \rightarrow 
 		\name \tau = j_{q'}.\end{array} \]
 	
 	Now one can define a code $c_{\varp}$ for a partial $|2^{\zeta_\varp}|$-Borel function from $\cP([0,\zeta_\varp])$ to $\{0,1\}$ requiring that 
 	\[ q_\varp \Vdash  \name \tau = f_{c_\varp}(\name X \cap [0,\zeta_\varp])\]  
 	(in fact it is even a $\zeta_\varp$-Borel function).

	\end{PROOF} 	
 	\end{PROOF}

 	Finally it is left to verify clause (h) from Definition $\ref{h2}$. Fix $p \in \bbP$ and $\bar{q}_\alpha = \langle q_{\alpha, \varp}: \ \varp < \kappa \rangle \in \cS_p$ ($\alpha <\lambda$). Now recall $\eqref{rM}$ $\eqref{sF}$, and let $A_\alpha' \in \bigcap \bar{U}_\kappa = \bigcap_{\beta \leq \theta} U_{\kappa, \beta}$ the set corresponding to the sequence $\bar{q}_\alpha$, i.e.
 	(if $d_0,d_1,\dots, d_n, d_{n+1}= \langle \kappa, A_{p,\kappa} \rangle$ denote the components of $p$)
 	\beeq \label{qk} \begin{array}{l} \bar{q}_\alpha = \langle q_{\alpha,\varp} :\varp < \kappa\rangle   \text{ where } q_{\alpha, \varp} = \langle d_0,d_1, \dots, d_n, \langle \kappa, A_{q_{\alpha,\varp}, \kappa} \rangle \rangle \text{ and} \\
 		\{A_{q_{\alpha,\varp}, \kappa}: \ \varp < \kappa\} \text{  lists }\{A_*:A_* \subseteq A_{p,\kappa} \text{ and } A_*=A_\alpha' \mod [\kappa]^{< \kappa} \}. \end{array} \eeq
 	
 	 Then for each fixed $\beta \leq \theta$ as $A'_\alpha \in U_{\kappa,\beta}$ ($\forall \alpha < \lambda$), using \eqref{mgd} $\eqref{ldir}$  there is $B_\beta \in U_{\kappa,\beta}$ such that $B_\beta \seq A_{p,\kappa}$, and
 	\begin{enumerate}[label = $(*_{\arabic*})$, ref = $*_{\arabic*}$]
 		\setcounter{enumi}{\value{sMa}}
 		\item \label{kiv} for each $\alpha < \lambda$ $|B_\beta \se A'_\alpha|< \kappa$, 		
 	\end{enumerate}
 	\stepcounter{sMa} 
 	Let
 	\begin{enumerate}[label = $(*_{\arabic*})$, ref = $*_{\arabic*}$]
 		\setcounter{enumi}{\value{sMa}}
 		\item $B_*= \bigcup_{\beta \leq \theta} B_\beta \in \bigcap \bar{U}_\kappa. $		
 	\end{enumerate}
 	\stepcounter{sMa} 
 	Therefore $\eqref{kiv}$ implies (recalling $\theta < \kappa$)
 	\begin{enumerate}[label = $(*_{\arabic*})$, ref = $*_{\arabic*}$]
 		\setcounter{enumi}{\value{sMa}}
 		\item for each $\alpha < \lambda$: $|B_* \se A'_\alpha|< \kappa$, therefore for some $\zeta_\alpha< \kappa$ we $B_* \cap (\zeta_\alpha, \kappa) \seq A_\alpha'$.
 	\end{enumerate}
 	\stepcounter{sMa} 
 	Defining $q_*$ as $ \langle d_0,d_1, \dots,d_n, \langle \kappa, B_*\rangle \rangle$ clearly $p \leq q_*$ as $B_* \seq A_{p,\kappa}$. Moreover, for any fixed $\alpha < \lambda$ by $\eqref{qk}$ there exists some $\varp < \kappa$ with the property that
 	\begin{enumerate}[label = $(*_{\arabic*})$, ref = $*_{\arabic*}$]
 		\setcounter{enumi}{\value{sMa}}
 		\item  $(A_{q_{\alpha,\varp}})_\kappa \cap (\zeta_\alpha, \kappa) = A'_\alpha  \cap (\zeta_\alpha, \kappa) \supseteq B_* \cap (\zeta_\alpha, \kappa)$, and
 		\item  $(A_{q_{\alpha,\varp}})_\kappa \cap [0,\zeta_\alpha] = B_* \cap [0,\zeta_\alpha]$,
 	\end{enumerate}
 	\stepcounter{sMa} \stepcounter{sMa} 
 	so $B_* \seq (A_{q_{\alpha,\varp}})_\kappa$, thus concluding $q_{\alpha,\varp} \leq q_*$.

 \end{PROOF}

Next we will give another example of a $(\lambda, \kappa)$-system, the Radin forcing, provided the measure sequence has some similar $\lambda^+$-directedness property.
\bdd \label{MS}
In order to state the following claim we need to prepare with the notions below.
\begin{enumerate}[label = (\roman*), ref = (\roman*)] 
	\item Let $\kappa$ be a cardinal $\bfj: \bfV \to \bfM$ be an elementary embedding (into a transitive inner model $\bfM$) with $\crr(\bfj) = \kappa$. We call the sequence $\overline{F} = \langle F(\alpha): \ \alpha < \dom(\overline{F}) \rangle$   a $\bfj$-sequence of ultrafilters, if 
	\begin{enumerate}
		\item $F(0) = \kappa$,
		\item  every $F(\alpha) \seq \cP(\bfV_\kappa)$,
		\item and for each $0<\alpha < \dom(\overline{F})$, $\forall X \seq \bfV_\kappa$: $[X \in F(\alpha)$ iff $(\overline{F} \um \alpha) \in \bfj(X)]$.
	\end{enumerate}
\item for each ultrafilter sequence $\overline{F}$ that are $\bfj$-sequence for some suitable $\bfj$ let $\kappa(\overline{F})$ denote the critical point of the witnessing $\bfj$, thus the $F_\alpha$'s are concentrated on $\bfV_{\kappa(\overline{F})}$, while for an ordinal $\alpha$ we mean $\alpha$ under $\kappa(\alpha)$. 
\end{enumerate}
Therefore for each $\alpha < \dom(\overline{F})$ $F(\alpha)$ is a $\kappa$-complete normal ultrafilter on $\bfV_\kappa$, where under normality we mean that for each sequence $\langle X_\beta: \ \beta < \kappa \rangle$ in $F(\alpha)$ the diagonal intersection
\[ \Delta_{\beta < \kappa} X_\beta = \{ \overline{f}: \  \forall \gamma <\kappa(\overline{f}): \overline{f} \in X_\gamma \} \ \in F(\alpha).\]
\begin{enumerate}[label = (\roman*), ref = (\roman*)]
	\stepcounter{enumi} \stepcounter{enumi}

\item 
Let $A^{(n)}$ ($n \in \omega$) be the following sequence of classes
\[ A^{(0)} = \{\overline{F}: \ \overline{F} \text{ is a } \bfj \text{-sequence of ultrafilters for some } \bfj:\bfV\to \bfM \}, \]
and
\[ A^{(n+1)} = \{\overline{F} \in A^{(n)}: \forall 0<\alpha \in \dom(\overline{F}) \ V_{\kappa(\overline{F})} \cap A^{(n)} \in F(\alpha) \}. \]

Finally let
\[ \bfA = \bigcap_{n \in \omega} A^{(n)}. \]
\item For any set $X \seq A^{(0)}$ and a set $I$ of ordinals let 
\[ X \um I = \{ \overline{F} \in X: \ \kappa(\overline{F}) \in I \}.\]
\end{enumerate}
\edd

\begin{claim}
\label{bRad}
There is a $(\lambda, \kappa)$-system such that $\bbR_{\bfr} =
\bbP$ \when \,:
\mn
\begin{enumerate}[label = $(\Alph*)$, ref = $(\Alph*)$]
	\stepcounter{enumi} \stepcounter{enumi}
	\item \label{rad}
	\begin{enumerate}[label = $(\alph*)$, ref = $(\alph*)$]
		\item   $\theta_* < \kappa < \lambda < 2^\kappa$,
		\sn
		\item  $\overline{F}_*$ is an ultrafilter sequence consisting of $\kappa$-complete ultrafilters on $\bfV_\kappa$, $\overline{F}_* \in \bfA$.
		\sn
		\item \label{radc} there exists $f: \kappa \to \kappa$ such that 
		$$\{\overline{F}: \ \dom(\overline{F}) < f(\kappa(\overline{F})) \} \in \bigcap \overline{F}_* = \bigcap_{0<\alpha < \dom(\overline{F}_*)} F_*(\alpha),$$
		(i.e.\ when for a witnessing $\bfj$ for $\overline{F}_*$ $\bfj(f)(\kappa) \geq \dom(\overline{F}_*)$, for example if $\dom(\overline{F}_*) \leq (2^{2^\kappa})^\bfM$),
		\sn	
		\item \label{ldirra} $\bigcap \overline{F}_* = \bigcap_{0<\alpha < \dom(\overline{F}_*)} F_*(\alpha)$ is $<$$\lambda^+$-directed in the following sense.
		For every sequence $\langle X_\alpha: \alpha < \lambda \rangle$ in $\bigcap \overline{F}_*$ there exists $X_* \in \bigcap \overline{F}_*$ such that
		\[ \forall \alpha < \lambda \ \exists \beta < \kappa: \ X_* \um (\beta, \kappa) \seq X_\alpha.\]
		
		\sn
		\item  $\bbP = \bbP_{\overline{F}_*}$ is the Radin forcing for $\overline{F}_*$ (see Definition $\ref{drd}$ below), 
		so  preserves
		the function $\mu \mapsto 2^\mu$, moreover, we can prescribe that in $\bfV^\bbP$ there is no new subset of $\theta_*$,
		and $\bbP$ satisfies the $\kappa^+$-c.c.
		\end{enumerate}
\end{enumerate}	 
\end{claim}
\begin{PROOF}{Claim \ref{bRad}}
	We will use the definition of the Radin forcing from \cite[Definition 5.2]{Gi10}.
	Observe that the definition only depends on $\bigcap \overline{F}_*$.

	\bdd \label{drd} For an ultrafilter sequence $\overline{F}_* \in \bfA$ we define the Radin forcing $\bbP$ 
	to be the collection of finite sequences of the form $p = \langle d_0, d_1, \dots, d_n, d_{n+1} = \langle \overline{F}_*, A_{p,\kappa} \rangle \rangle$, where
	\newcounter{sRa}
	\setcounter{sRa}{0}
	\begin{enumerate}[label = $(*_{\arabic*})$, ref = $*_{\arabic*}$]
		\setcounter{enumi}{\value{sRa}}
		
		\item
		\begin{enumerate}
			\item $A_{p,\kappa} \in \bigcap \overline{F}_* = \bigcap_{0<\alpha < \dom(\overline{F}_*)} F_*(\alpha)$, $A_{p,\kappa} \in \bfA$,
			\item each $d_j$ ($j \leq n$) is either of the form 
			\begin{itemize}
				\item $\langle \overline{F}_{d_j}, A_{d_j} \rangle$ where $\overline{F}_{d_j} \in \bfA$, $A_{d_j} \seq \bfA$, moreover,
				$$A_{d_j} \in \bigcap \overline{F}_{d_j} = \bigcap_{0< \gamma < \dom(\overline{F}_{d_j})} F_{d_j
				}(\gamma).$$ 
				If $\varp = \kappa(\overline{F}_{d_j}$) we may refer to $\langle \overline{F}_{d_j}, A_{d_j} \rangle$ as $\langle \overline{F}_{p,\varp}, A_{p,\varp} \rangle$, and we also define $\kappa(d_j) = \kappa(\overline{F}_{d_j})$.
				\item or $d_j = \varp$ for some $\varp< \kappa$ (when we let $\kappa(d_j) = \varp$).
			\end{itemize}
			\item $\kappa(d_0) < \kappa(d_1) < \dots < \kappa(d_n) < \kappa(d_{n+1}) = \kappa$,
			\item moreover, for each $j \leq n$ if $d_{j+1}$ is a triplet, then $A_{p,\kappa(d_{j+1})} \cap V_{\kappa(d_j)} = \emptyset$.
		\end{enumerate}
		\stepcounter{sRa} \stepcounter{sRa}
		
		\item For the sequences $$p = \langle d_0, d_1, \dots, d_n, d_{n+1} = \langle  \overline{F}_*, A_{p,\kappa} \rangle\rangle,$$ $$q = \langle e_0, e_1, \dots, e_n, e_{m+1} = \langle \overline{F}_*, A_{q,\kappa} \rangle \rangle$$    we let $p \leq q$, if
			\begin{enumerate}
			\item  $m \geq n$, and
			\item  there exists a sequence $0 \leq i_0 < i_1 < \dots < i_n < j_{n+1} =m$ such that for each $j \leq n+1$ we have
			\begin{itemize}
				\item $\kappa(d_j) = \kappa(e_{i_j})$,
				\item and $$\text{either }\overline{F}_{p, \kappa(d_j)}= \overline{F}_{q, \kappa(e_{i_j})} \text{ and } A_{q,\kappa(e_{i_j})} \seq A_{p, \kappa(d_j)},$$
				 or $d_j= e_{i_j} = \kappa(d_j) = \kappa(e_{i_j}),$
			\end{itemize}
			\item 
			moreover, for each $k \leq m$ not of the form $i_j$ ($j \leq n+1$), if $i_l = \min \{ i_j: \ j \leq n+1, \ i_j > k\}$, then 
			$$A_{q,\kappa(e_k)} \cup \{ \overline{F}_{q,\kappa(e_k)}\} \seq A_{p,\kappa(d_l)}.$$
		\end{enumerate}
	\end{enumerate}
	
	\edd

	\bdd We define $p \leq_* q$ to be true iff $p \leq q$ and $\lh(p) = \lh(q)$.
	\edd
	
		We define the system $\bfr$ by letting:
	\mn
	
	\begin{enumerate}[label = $(*_{\arabic*})$, ref = $*_{\arabic*}$]
		\setcounter{enumi}{\value{sRa}}
		\item \label{rR}
		\begin{enumerate}[label = (\alph*), ref = \alph*]
			\item  $\kappa_{\bfr} = \kappa$,
			\sn
			\item  $\lambda_{\bfr} = \lambda$,
			\sn
			\item  $\bbR_{\bfr} = \bbP$,
			\sn
			\item  let $\name X_{\bfr}$ be  the generic sequence, i.e. 
			$$\name X_\bfr = \cup\{ \{ \kappa(d_j): j < \lh(p) \} : \ p = \langle d_0, d_1, \dots, d_{\lh(p)-1} \rangle  \in \bfG_{\bbP}\} \se \{ \kappa \},$$
			\item  $\le_{\pr} =\le_{\bfr,\pr}$ is defined by $p \le_{\pr} q$
			iff $p \leq_* q$,
			\sn
			\item \label{sFs}  for $p = \langle d_0, d_1, \dots, d_n, d_{n+1} = \langle \overline{F}_*, A_{p,\kappa} \rangle \rangle \in \bbR_{\bfr} = \bbP$ let 
			$$\cS_p = \cS_{\bfr,p} := \left\{ \begin{array}{l}
				\bar q:\bar q = \langle q_\varp:\varp < \kappa\rangle \text{, where} \\
				(\bullet_1) q_\varp = \langle d_0, d_1, \dots, d_{n}, \langle \overline{F}_*, A_{q_\varp,\kappa} \rangle\rangle, \text{ and} \\
				\text{for some } B \in \bigcap \overline{F}_*  \text{ we have } \\
				(\bullet_2) \ B \subseteq A_{p,\kappa},  \text{ and} \\
				(\bullet_3) \  \{A_{q_{\varp}, \kappa}:\varp < \kappa\} \text{  lists }\{A_*:A_* \subseteq A_{p,\kappa} \ \wedge \ A_*=B \mod [\kappa]^{< \kappa}\} \end{array} \right\}. $$
		\end{enumerate}
		\stepcounter{sRa} 
	\end{enumerate}

	Now we check the requirements of Definition $\ref{h2}$.
	
	It is known that  if a condition $\langle \langle \beta \rangle, \langle \overline{F}_*, A_\kappa \rangle \rangle$ is in the generic filter (for some $\beta < \kappa$) then the forcing adds  no new subset of $\beta$. This implies that as $\bigcap \overline{F}_* \seq F_*(0)$, which is concentrated on the ordinals, i.e.\ on $\kappa$ itself, w.\ l.\ o.\ g.\ we can assume that $\langle \beta, \langle \overline{F}_*, A \rangle \rangle \in \bfG$ for some $\beta \geq \theta_*$.
	
	Now we have only to check the requirements of Definition $\ref{h2}$.
	Recall the following properties of the Radin forcing, see \cite[Sec. 5.1]{Gi10}.
	\bfact (Prikry Lemma)
	For each $p \in \bbP$ and each formula $\sigma(\name x_0, \dots, \name x_m)$ there exists $q \geq_* p$, $q \parallel \sigma(\name x_0, \dots, \name x_m)$ (i.e. either $q \Vdash \sigma(\name x_0, \dots, \name x_m)$, or $ q \Vdash \neg \sigma(\name x_0, \dots, \name x_m)$).
	\efact

	The following claims, which complete the proof of Claim \ref{bRad} have the same proofs as in the case of Magidor forcing. In Claim $\ref{kisebbkR}$ condition \ref{rad}/\ref{radc} is essential for the argument.
	\bcl \label{kisebbkR}
	For each $p = \langle d_0, d_1, \dots, d_{n+1} = \langle \overline{F}_*, A_{p,\kappa} \rangle  \rangle\in \bbP$, $\name \tau$ (with $p \Vdash \name \tau \in \{0,1 \}$) there exists a set $A' \in \bigcap \overline{F}_*$, $A' \subseteq A_{p,\kappa}$, such that whenever $q = \langle e_0,e_1, \dots, e_m, \langle \overline{F}_*, A_{q,\kappa} \rangle \rangle \geq p' = \langle d_0, d_1, \dots, d_n, \langle \overline{F}_*, A'\rangle\rangle $, $\alpha \geq \kappa(e_m)$ are given and $q$ forces  a value for $\name \tau$, then so does $$q' = \left\langle e_0,e_1, \dots, e_m, \langle \overline{F}_*,  A' \um (\alpha, \kappa) \rangle \right\rangle.$$
	\ecl

	\bcl For every $p  \in \bbP$ and $\name\tau$, if $p \Vdash \tau \in \{0,1 \}$, then 
	there exists $\bar{q} = \langle q_\varp: \ \varp < \kappa \rangle \in \cP_p$, $\langle \gamma_\varp: \ \varp< \kappa \rangle \in \ ^\kappa \kappa$,
	$\langle c_\varp: \ \varp < \kappa \rangle$ such that each $c_\varp$ is a code for a $\gamma_\varp$-Borel function from $ \cP(\gamma_\varp)$ to $\{0,1 \}$, and
	\[ q_\varp \Vdash \ \name\tau = f_{c_\varp}(\name X). \]
	\ecl

\end{PROOF}

\bigskip

\section {The preparatory forcing} \label{3}
\subsection{The general frame} \label{3a}
\bigskip

This subsection is dedicated for the preparation, in Claim $\ref{c10}$ we provide a general frame to force a $(\lambda,\kappa)-1$ system.

\begin{definition}
\label{c7}
For $D$ a $<\kappa$-directed system (i.e. generating a $\kappa$-complete filter $D^*$) on $\cup D \subseteq V_\kappa$ (so $D^* \subseteq \cP(\cup D)$) let $\bbQ = \bbQ_D$ be
the following forcing notion (with the notations from Definition $\ref{MS}$, also applying to filters concentrated on $\kappa$):
\mn
\begin{enumerate}
\item[(A)]   $p \in \bbQ$ \Iff \,
\sn
\begin{enumerate}
\item[(a)]  $p = (w,A) = (w_p,A_p)$, and for some $\sigma_p < \kappa$ we have
\sn
\item[(b)]  $w_p \subseteq V_\kappa$,  $w_p = w_p \upharpoonright [0, \sigma_p)$ ,
\sn
\item[(c)]  $A_p \subseteq \cup D$, $A_p\in D^*$ and $A_p = A_p \upharpoonright [\sigma_p, \kappa)$.
\end{enumerate}
\sn
\item[(B)]  $\bbQ \models p \le q$ \Iff \,
\sn
\begin{enumerate}
\item[(a)]  $p,q \in \bbQ$,
\sn
\item[(b)]  $w_p \subseteq w_q \subseteq w_p \cup A_p$,
\sn
\item[(c)]  $A_p \supseteq A_q$,
\end{enumerate}
\sn
\item[(C)]  $\name w = \cup\{w_p:p \in \name{\bfG}\}$.
\end{enumerate}
\end{definition}

\begin{claim}
\label{c10}
If (A) and (B) hold, then (C) where:
\mn
\begin{enumerate}[label = $(\Alph*)$, ref = $(\Alph*)$]
\item  \label{c10A} $\bfv = (\bfV_0,\kappa,\bfh,\bfp,\bfG_\kappa,\bfV_1)$
  satisfies:
\sn
\begin{enumerate}
\item[(a)]  $\bfV_0$ is a universe of set theory,
\sn
\item[(b)]  in $\bfV_0$ $\kappa$ is supercompact and $\bfh:\kappa
  \rightarrow \cH(\kappa)$ is a Laver diamond,
\sn
\item[(c)]   $\bfp$ is the Easton support iteration $\langle
  \bbP_{\bfp,\alpha},\name{\bbQ}_{\bfp,\beta}:\alpha \le \kappa,\beta <
  \kappa\rangle = \langle \bbP^0_\alpha,\name{\bbQ}^0_\beta:\alpha \le
  \kappa,\beta < \kappa\rangle$ built essentially as in Laver \cite{L78} from $\bfh$
  and let $\bbP_{\bfp} = \bbP_{\bfp,\kappa}$ (hence for $\alpha < \kappa$ also $\bbP^0_\alpha \in V_\kappa^{\bfV_0}$), see Definition $\ref{P0def}$,
\sn
\item[(d)]   $\bfG_\kappa = \bfG_{\bfp,\kappa} \subseteq \bbP_{\bfp}$
  is generic over $\bfV_0$ and $\bfV = \bfV_1 = \bfV_0[\bfG_\kappa]$.
\end{enumerate}
\sn
\item \label{B}
\begin{enumerate}[label=(\alph*), ref = \alph*]
\item \label{(a)} $\kappa < \lambda < \chi = \chi^\lambda$,

\sn

\item \label{(c)}  $\bbP^1_\chi = \langle \bbP^1_\alpha,\name \bbQ^1_\beta:\alpha \le \chi,\beta <
  \chi\rangle \in \bfV_1$ is an iteration with $< \kappa$ support such that $\bbP^1_\chi$ is
  $\lambda^+$-c.c.\ and $<$$\kappa$-directed closed, preserving cardinals,
  \item for each $\alpha < \chi$ $$\bfV_1^{\bbP^1_\alpha} \models \ |\bbQ^1_\alpha| \leq \chi.$$
\sn

\item \label{d}  for the set $S^* \seq \chi$ there is a system $\langle \name D_\delta: \ \delta \in S^* \rangle \in \bfV_1$, $\name D_\delta$ is a $\bbP^1_\delta$-name of a subset of $\cP^{\bfV_1^{\bbP^1_\delta}}(V_\kappa)$, and if $\name{D}_\delta$  generates a $\kappa$-complete filter satisfying $(\forall \alpha < \kappa) \ |\cup D_\alpha) \um \alpha|<\kappa$, then the forcing $\name{\bbQ}^1_\delta$, $\delta \in S^*$ is of the form
$Q_{\name D_\delta}$, the forcing from Definition $\ref{c7}$. Moreover, $\langle \name D_\delta: \ \delta \in S^* \rangle $ satisfies the following
\begin{enumerate}[label = ($\#$), ref = ($\#$)]

 \item \label{cf0}  $ \bfV_1^{\bbP^1_\chi} \models \  \forall X \subseteq V_\kappa \ \forall D \in [\cP(X)]^{\leq \lambda}:$
 \begin{itemize}
 	\item $ (\forall \alpha < \kappa \ |X \um \alpha|<\kappa)\  \wedge $
 	\item $ (D \text{ generates a proper  }<\kappa\text{-complete filter})$
 \end{itemize} 
  	$\ \ \ \longrightarrow (D =  D_\delta  \text{ for some } \delta \in S^*)$.
 \end{enumerate}
\end{enumerate}
\sn
\item \label{C} in $\bfV_1^{\bbP^1_\chi}$ we have $2^\kappa$ is $\chi$, and the following.
 \begin{enumerate}[label = (\alph*), ref = \alph* ]
 	\item \label{Prikry}  There is a $\kappa$-complete normal ultrafilter $U$, which is $<\lambda^+$-directed mod $[\kappa]^{<\kappa}$.
 	\item \label{Mag} (Setting for Magidor forcing:) There is a sequence $\bar U = \langle U_i:i < \kappa\rangle$
 	of normal ultrafilters on $\kappa$, strictly increasing in the Mitchell order,
 	i.e. $i < j \Rightarrow U_i \in \Mos \Col(^\kappa(\bfV^{\bbP^1_\chi})/U_j)$, such that each  $U_i$ is $<\lambda^+$-directed mod $[\kappa]^{<\kappa}$.
 	\item \label{Rad} (Setting for Radin forcing:) For any $\Upsilon \geq \kappa$ and $\eta$ there is a $\kappa$-complete fine normal ultrafilter $W$ on $[\Upsilon]^{<\kappa}$ such that for the elementary embedding $\bfj_W$ of $\bfV_1^{\bbP^1_\chi}$ with critical point $\kappa$ we have (letting $\overline{U}$ denote the  measure sequence associated to $\bfj_W$):
 	 \begin{enumerate}[label = $(\star)$, ref = $\star$ ]
 	 	\item for every $\sigma \leq \min(\dom(\overline{U}, \eta))$ if the filter $\bigcap (\overline{U} \um \sigma) = \bigcap_{\gamma < \sigma} U_\gamma$ concentrates on a set $X \subseteq V_\kappa$ with $(\forall \alpha < \kappa)$ $|X \um \alpha|<\kappa$,  then  $\bigcap (\overline{U} \um \sigma)  \text{ is } <\lambda^+$-directed in the following sense: 
 	 	 	Whenever $\langle A_i: \ i < \lambda \rangle$ ($\forall i < \lambda$ $A_i \in \bigcap (\overline{U} \um \sigma)$) is given, there exists $A_* \in \bigcap (\overline{U} \um \sigma)$ such that
 	\beeq \forall i \in \lambda \ \exists \delta_i < \kappa: \ A_*\upharpoonright [\delta_i, \kappa) \subseteq A_i. \eeq
 	
 	\end{enumerate}
 	In particular $\kappa$ is supercompact.
 \end{enumerate}
 \end{enumerate}
\end{claim}

\begin{remark}
\label{b8}
This continues D{\v z}amonja-Shelah \cite{Sh:659}.
\end{remark}

\begin{PROOF}{Lemma \ref{c10}}
	First we construct the iteration $\bbP^0$ using the Laver function $\bfh: \kappa \to \cH(\kappa) \in \bfV_0$. The construction $\bbP^0 = \langle \bbP^0_\alpha, \name \bbQ^0_\beta: \ \alpha \leq \kappa, \beta < \kappa \rangle$ goes by induction, we follow \cite{L78}, only with a slight technical modification which we will need in the proof of $\eqref{C}$$\eqref{Mag}$.
	
	Let $\bfh$ be as in \cite{L78} (i.e.
	\newcounter{bcounter}
	\begin{enumerate}[label = $(\bullet)_{\arabic*}$, ref=$(\bullet)_{\arabic*}$]
		\item \label{bfhd} for each $\lambda \geq\kappa$, $x \in \cH(\lambda^+)$ there exists a $\kappa$-complete fine normal ultrafilter $U$ on $[\lambda]^{<\kappa}$ such that for the associated elementary embedding $\bfj_U$
		\[\bfj_U(\bfh)(\kappa) = x).\]
		\stepcounter{bcounter}
	\end{enumerate}
		\newcounter{btc}
	\setcounter{btc}{0}

	\bdd \label{P0def}
		We define $\bbP^0 = \langle \bbP^0_\alpha, \name \bbQ^0_\beta: \ \alpha \leq \kappa, \beta < \kappa \rangle$ and $\langle \mu_\alpha: \ \alpha < \kappa \rangle$ by induction. If $\langle \bbP^0_\alpha, \name \bbQ^0_\beta: \ \alpha <\gamma, \beta < \gamma \rangle$ are already defined, then 
		\begin{enumerate}[label = $(\bullet)_{\Roman*}$, ref=$(\bullet)_{\Roman*}$]
			\item if $\gamma$ is strongly inaccessible then $\bbP^0_\gamma$ is the direct limit (i.e.\ we use bounded support),
			\item otherwise let $\bbP^0_\gamma$ be the inverse limit of $\bbP^0_\beta$'s ($\beta < \gamma$) (i.e.\ for a function $p$ with $\dom(f) = \gamma$ $p \in \bbP^0_\gamma$ iff $(\forall \beta < \gamma)$ $ p \um \beta \in \bbP^0_\beta$).  .
		\end{enumerate}
		Second, 
		\begin{enumerate}[label = $(\bullet)_{\alph*}$, ref=$(\bullet)_{\alph*}$]
		\item if $\sup \{ \mu_\alpha: \alpha < \gamma\} \leq \gamma$, $\gamma$ is inaccessible, and  $\bfh(\gamma)= \langle \name Q_*, \mu_*, \name U \rangle$, where $\name Q_*$ is a $\bbP^0_\gamma$-name for $<$$\gamma$-directed closed notion of forcing, $\mu_*$ is an ordinal, $\name U$ is a (possibly trivial) $\bbP^0_\gamma$-name, then let 
		\[ \name \bbQ^0_\gamma = \name \bbQ_*, \ \ \mu_\gamma = \mu_*. \]
		\item In the remaining case let $\name \bbQ^0_\gamma$ be the trivial forcing, $\mu_\gamma = \gamma$.
		\end{enumerate}
	\edd

	Recall $\bfG^0_\kappa
	\subseteq \bbP^0_\kappa$ is generic over $\bfV_0$ such that
	$\bfV_0[\bfG^0_\kappa] = \bfV_1$, and let $\bfG^1_\chi \subseteq \bbP^1_\chi$ be generic
	over $\bfV_1$, let $\bfV_2 = \bfV_1[\bfG^1_\chi] = \bfV_0[\bfG^0_\kappa \ast \bfG^1_\chi]$. 
	Note that as $|\bbP^0_\kappa| = \kappa$, $\kappa < \lambda$, $\eqref{B}$ implies that 
	\begin{enumerate}[label = $(\bowtie)_{\arabic*}$, ref = $(\bowtie)_{\arabic*}$]
		
		\setcounter{enumi}{\value{btc}}
		\item 	\label{chilam}  $\bfV_1 = \bfV_0[\bfG^0_\kappa]  \models  \ \chi^\lambda = \chi$, thus $\cf(\chi)>\lambda$ holds too.
		\stepcounter{btc}
	\end{enumerate}
	Also note that as $\bbP^1_\chi$ is $<\kappa$-closed
	\begin{enumerate}[label = $(\bowtie)_{\arabic*}$, ref = $(\bowtie)_{\arabic*}$]
		
		\setcounter{enumi}{\value{btc}}
		\item 	\label{Vk}  $V_\kappa^{\bfV_2} = V_\kappa^{\bfV_1}$, and $\bfV_2 \models  ``\kappa$ is still inaccessible."
		\stepcounter{btc}
	\end{enumerate}

	First observe that because of our cardinal arithmetic assumptions $\chi^\kappa \leq \chi^\lambda = \chi$ ($\eqref{B} \eqref{(a)}$), and as $|\bbP^0_\kappa| = \kappa$, we have $(\chi^\lambda)^{\bfV_1} = \chi^{\lambda \cdot \kappa} = \chi$, so by an easy induction (and by the $\lambda^+$-cc) $|\bbP_\chi^1|^{\bfV_1} = \chi$, so 
	\begin{enumerate}[label = $(\bowtie)_{\arabic*}$, ref = $(\bowtie)_{\arabic*}$]
		
		\setcounter{enumi}{\value{btc}}
		\item 	 $|\bbP^0_\kappa \ast \name \bbP^1_\chi| = \chi$ up to equivalence.
		\stepcounter{btc}
	\end{enumerate}
	Recalling $\chi^\lambda = \chi$ again, clearly
	\begin{enumerate}[label = $(\bowtie)_{\arabic*}$, ref = $(\bowtie)_{\arabic*}$]
		
		\setcounter{enumi}{\value{btc}}
		\item 	 \label{bt0} $\bfV_0[\bfG^0_\kappa \ast \bfG^1_\chi] \models$ $2^\chi = (2^\chi)^{\bfV_0}$,
		\stepcounter{btc} \stepcounter{btc}
		\item \label{bt00}
		$\bfV_0[\bfG^0_\kappa \ast \bfG^1_\chi] \models$ $2^\kappa = \chi$.
	\end{enumerate}
		

\bdd \label{1def}  We have to introduce the following objects.

\begin{enumerate}[label = $(\bullet)_{\arabic*}$, ref=$(\bullet)_{\arabic*}$]
	\setcounter{enumi}{\value{bcounter}}
\item \label{jM} Let $\bfj:\bfV_0 \rightarrow \bfM$ be an elementary embedding with
critical point $\kappa$ such that $(\bfj(\bfh))(\kappa) = \langle \name{\bbP}^1_\chi, \chi^+, \check{\emptyset} \rangle$ ($\check{\emptyset} = \emptyset$ is the canonical name for the empty set)
and $\bfj(\kappa)> \chi$, $^{\chi} \bfM \subseteq \bfM$, 
\item Let $\langle \bbP^0_\alpha,\name{\bbQ}^0_\beta:\alpha \le
\bfj(\kappa),\beta < \bfj(\kappa)\rangle = \bfj(\langle
\bbP^0_\alpha,\bbQ^0_\beta:\alpha \le \kappa,\beta < \kappa\rangle)$
so $\name{\bbQ}^0_\kappa = \name{\bbP}^1_\chi$, and
\item let $\name \bbP'_{\bfj(\chi)} = \bfj(\name \bbP^1_\chi)$, i.e.\ 
\[ \text{(a }\bbP^0_{\bfj(\kappa)}\text{-name for a }<\bfj(\kappa)\text{-directed closed notion of  forcing)}^\bfM. \]
  (Recall that $\name \bbP^1_\chi$ is a $\bbP^0_\kappa$-name for the iteration $\langle
\name{\bbP}^1_\alpha,\name{\bbQ}^1_\beta:\alpha \le \chi,\beta <
\chi\rangle) \in \bfV_0^{\bbP^0_\kappa}$.)
\stepcounter{bcounter}
\stepcounter{bcounter}
\stepcounter{bcounter}
\end{enumerate}
\edd
Similarly to \ref{Vk}
	\begin{enumerate}[label = $(\bowtie)_{\arabic*}$, ref = $(\bowtie)_{\arabic*}$]
	
	\setcounter{enumi}{\value{btc}}
	\item 	\label{Mk}  $V_\kappa^{\bfM[\bfG^0_{\kappa+1}]} = V_\kappa^{\bfM[\bfG^0_\kappa]} = V_\kappa^{\bfV_2}$, ($\kappa$ is inaccessible)$^{\bfM[\bfG^0_{\kappa+1}]}$.
	\stepcounter{btc}
\end{enumerate}

From now on we will identify $\bbP^0_\kappa \ast \name \bbP^1_\chi$ with the $(\kappa+1)$-step iteration $\bbP^0_{\kappa+1}$, and also

\begin{enumerate}[label = $(\bowtie)_{\arabic*}$, ref = $(\bowtie)_{\arabic*}$]
	
	\setcounter{enumi}{\value{btc}}
	\item 	  $\bfG^0_{\kappa +1} = \bfG^0_\kappa * \bfG^1_\chi$ is a generic subset of
	$\bbP^0_{\kappa +1} = \bbP^0_\kappa \ast \name \bbP^1_\chi$ (over $\bfV_0$).
	\stepcounter{btc}
\end{enumerate}



\br Having completed the requirements of Claim $\ref{c10}$ we remark that given a scheme for an iteration fitting all our assumptions other than $\eqref{B}$$\eqref{d}$, it is easy to adapt it to have \ref{cf0} using $\chi^\lambda = \chi$ \ref{chilam}.
\er

\noindent Now we can prove the statements in \ref{c10}(C). \\
\underline{Case 1}:  First we verify \ref{c10}(C)(a). \\ We would like a suitable $\kappa$-complete ultrafilter in $\bfV_0[\bfG^0_\kappa \ast \bfG^1_\chi]$, for which   we will use the basic trick: using the elementary embedding $\bfj: \bfV_0 \to \bfM$, extending $\bfV_0$ with $\bfG^0_\kappa \ast \bfG^1_\chi$, $\bfM$ with $\bfG^0_{\kappa+1} (= \bfG^0_\kappa \ast \bfG^1_\chi)$, and finding a single condition in $\bbP^0_{\bfj({\kappa})}  \ast \bbP'_{\bfj(\chi)} / \bfG^0_{\kappa+1}$ compatible with $\{\bfj(p \um \{\kappa\}) = \bfj(p)\um \{ \bfj(\kappa) \}: \ p \in \bfG^0_\kappa \ast \bfG^1_\chi \}$ giving us sufficient  information
  (just as if there existed some extension $ \tilde{\bfj}: \bfV_0[\bfG^0_\kappa \ast \bfG^1_\chi] \to \bfM[\mathbf{H}_{\bfj(\kappa)}^0 \ast \mathbf{H}'_{\bfj(\chi)}]$).

We will need the following facts.
\bfact \label{mf}
 The filter $\bfG^0_{\kappa +1}$ is generic over $\bfM$ as well, and the forcing notions
$\bbP^0_{j(\kappa)}/ \bfG^0_{\kappa +1}$  and $(\bbP^0_{j(\kappa)} *
\name{\bbP}'_\gamma)/\bfG^0_{\kappa +1}$ ($\gamma \le \bfj(\chi)$) are well defined and
$<\chi^+$-directed closed in $\bfM[\bfG^0_{\kappa +1}]$.
\efact
\begin{PROOF}{} 
	First recall that a pair $(p,\name q) \in  (\bbP^0_{j(\kappa)} *
	\name{\bbP}'_\bfj(\chi))/\bfG^0_{\kappa +1}$ iff $p = p_0 \um (\kappa, \bfj(\kappa))$ for some $p_0 \in \bbP^0_{\bfj(\kappa)}$, and $( \Vdash_{\bbP^0_{\bfj(\kappa)}} \name q \in \name \bbP'_{\bfj(\chi)})^\bfM$.
	We only have to refer to the construction of the iteration Definition $\ref{P0def}$ i.e.\ recall that 
	\begin{enumerate}[label = $(\roman*)$, ref = $(\roman*)$]
		\item  we have $\Vdash_{\bbP^0_{\kappa}}  `` \name \bbP^1_\chi$ is $<\kappa$-support iteration of $<\kappa$-directed closed forcing notions", and
		\item  for each $\alpha \leq \beta < \kappa$ we have that $\Vdash_{\bbP^0_{\beta}} `` \name \bbQ^0_{\beta}$
		 is $<\beta$-directed closed", and is the trivial forcing if $\beta < \sup \{ \mu_\varrho: \ \varrho < \beta \}$ (in particular, if $\beta < \sup \{ \mu_\varrho: \ \varrho < \alpha \}$),
		 \item for each  $\alpha < \beta  < \kappa$, where $\beta$ is limit and $\cf(\beta) < \mu_\alpha$ the iteration $\bbP^0_\beta$ is the inverse limit of $\bbP^0_\delta$'s ($\delta < \beta$).
	\end{enumerate}
		So using \cite[Thm. 5.5]{B3}  for each $\alpha < \beta < \kappa$ the quotient
	$(\bbP^0_{\kappa} *
	\name{\bbP}^1_\chi)/\bfG^0_{\alpha}$ (of the $\kappa+1$-long iteration $\bbP^0_{\kappa} *
	\name{\bbP}^1_\chi = \bbP^0_{\kappa+1}$)   is $<\beta$-directed closed in $\bfV_0[\bfG^0_{\alpha}]$ if $\beta \leq \sup \{ \mu_\varrho: \ \varrho < \alpha \}$, and $\bbP^0_\alpha$ has the $\beta$-cc. Thus by elementarity (letting $\alpha = \kappa+1$, $\beta = \chi^+ = \mu_\alpha$, recalling $(\chi^+)^\bfM = \chi^+$ by $^{\chi} \bfM \subseteq \bfM$):
	\[ \bfM[\bfG^0_{\kappa +1}] \models \ \text{"}(\bbP^0_{j(\kappa)} *
	\name{\bbP}'_\bfj(\chi))/\bfG^0_{\kappa +1} \text{ is } <\chi^+ \text{-directed closed."} \]
\end{PROOF}

\bfact \label{Cf} $\bfV_1 \models \ \text{``}(\bbP^0_{j(\kappa)} * \name{\bbP}'_\gamma)/\bfG^0_{\kappa +1} \text{ is } <\chi^+\text{-directed closed."}$
\efact 
Fact $\ref{Cf}$ follows from the fact below.

  \bfact	 \label{bt1}
  	 $\bfV[\bfG^0_\kappa \ast \bfG^1_\chi] \models \ ``^\chi\bfM[\bfG^0_{\kappa +1}] \subseteq \bfM[\bfG^0_{\kappa +1}]"$.
  \efact 	 
  \begin{PROOF}{Fact \ref{bt1}}
	
	For \ref{bt1}  pick a name $\name f$ for a function $\name f : \chi \to \bfM[\bfG_{\kappa+1}^0]$, and observe that w.l.o.g.\ we can assume that $\name f: \chi \to \ord$, i.e. for each $\alpha < \chi$ $f(\alpha)$ is an ordinal, in particular $\ran(f) \subseteq \bfM$.
	Now for each $\alpha$ there exists a maximal antichain $A_\alpha = \{a^\alpha_i: \ i < |A_\alpha| \} \seq \bbP^0_{\kappa+1}$, and $\{  x^\alpha_i: \ i< |A_\alpha \} \seq \bfM$, s.t. $a^\alpha_i \Vdash \name f(\alpha) = x^\alpha_i$. Now as $\bbP^0_{\kappa+1} = \bbP^0_\kappa \ast \name \bbP^1_\chi$ is of power $\chi$, we have $|A_\alpha| \leq \chi$, therefore as $\bfM$ is closed under sequences of length $\chi$ (\ref{jM}, Definition $\ref{1def}$) there is indeed a name $\name g \in \bfM$, such that $\Vdash \name f = \name g$.
\end{PROOF}

\bdd \label{fdf}
In $\bfV_0[\bfG^0_{\kappa +1}]$, $\zeta \in S^*$ let 
\begin{enumerate}
	\item $\varp_\zeta \in \bfV[\bfG^0_\kappa \ast \name \bfG^1_{\zeta+1}]$ denote the generic subset of $V^{\bfV_1}_\kappa$ (or just $\kappa$) given by $\bbQ^1_\zeta$, i.e.
	\[ \Vdash_{\bbP^0_{\kappa} * \name{\bbP}^1_{\zeta+1}} \ \name\varp_\zeta = \cup \{ \name \varp: \ \exists \name{A}: (\name \varp,\name A) \in \bfG_{\name \bbQ^1_\zeta}\}\]
	(after identifying $\bbP^0_{\kappa} * \name{\bbP}^1_{\zeta+1} =\bbP^0_{\kappa} * (\name{\bbP}^1_{\zeta} \ast \name \bbQ^1_\zeta)$ with $(\bbP^0_{\kappa} * \name{\bbP}^1_{\zeta}) \ast \name \bbQ^1_\zeta$).
	\item Define $\cN_\zeta$ to be a set of $\bbP^0_\kappa \ast \name \bbP^1_\zeta$-names of subsets of $V_\kappa$ containing exactly one name from each equivalence class, i.e.\ no $\name A \neq \name B \in \cN_\zeta$ $\Vdash_{\bbP^0_\kappa \ast \name \bbP^1_\zeta} \name A = \name B$, but each set in the extension is represented.
\end{enumerate}

Observe that by \ref{Mk} we can assume that 
\begin{enumerate}[label = $(\bowtie)_{\arabic*}$, ref = $(\bowtie)_{\arabic*}$]
	
	\setcounter{enumi}{\value{btc}}
	\item 	  $\cN_\zeta \subseteq \bfM$,
	\stepcounter{btc}
\end{enumerate}
and as  $|V^{\bfV_2}_\kappa|  = \kappa$, and by the $\lambda^+$-cc \ref{B}\ref{(c)}   \begin{enumerate}[label = $(\bowtie)_{\arabic*}$, ref = $(\bowtie)_{\arabic*}$]
	
	\setcounter{enumi}{\value{btc}}
	\item \label{nze}	  $|\cN_\zeta| \leq |\bbP^0_\kappa \ast \bbP^1_\zeta|^\lambda = \chi$,
	\stepcounter{btc}
\end{enumerate}
so by $^{\chi} \bfM \subseteq \bfM$ we can assume that
\begin{enumerate}[label = $(\bowtie)_{\arabic*}$, ref = $(\bowtie)_{\arabic*}$]
	
	\setcounter{enumi}{\value{btc}}
	\item 	  $\cN_\zeta \in \bfM$, and  $\bfj \um \cN_\zeta \in \bfM$.
	\stepcounter{btc}
\end{enumerate}
	
\begin{enumerate}
	\setcounter{enumi}{2}	
	\item Using the notation
	\[  \cA_{\name{\bbQ}^1_\zeta} =  \{ \name A \in \cN_\zeta: \ (\name\varp,\name A) \in \bfG_{\name{\bbQ}^1_\zeta} \text{ for some } \name\varp \}, \] 
	note that  $\cA_{\name{\bbQ}^1_\zeta} \in \bfM[\bfG^0_\kappa \ast \bfG^1_{\zeta+1}]$
		(so $\name \cA_{\name \bbQ^1_\zeta}$ is a $\bbP^0_\kappa \ast \name \bbP^1_{\zeta+1}$-name for a set of $\bbP^0_\kappa \ast \name \bbP^1_\zeta$-names).
	Now similarly $\bfj`` \cA_{\name{\bbQ}^1_\zeta} = \{ \bfj(\name A): \ \name A \in \cA_{\name \bbQ^1_\zeta}\} \in \bfM[\bfG^0_\kappa \ast \bfG^1_{\zeta+1}]$ is a set of $\bbP^0_{\bfj(\kappa)} \ast \name \bbP'_{\bfj(\zeta)}$-names, and we can define the $\bbP^0_{\bfj(\kappa)} \ast \name \bbP'_{\bfj(\zeta)}$-name  $\name A'_{\bfj(\zeta)} \in \bfM$ for a subset of $V_{\bfj(\kappa)}$ so that
	\[ \bfM \models \ \  ``\Vdash_{\bbP^0_{\bfj(\kappa)} \ast \name\bbP'_{\bfj(\zeta)}} \  \name A'_{\bfj(\zeta)} = \cap \{\bfj(\name A): \ \name A \in \cA_{\name{\bbQ}^1_\zeta} \}''. \]
	
\end{enumerate}


\edd

\begin{claim} \label{cind} There is a sequence $\langle q_\zeta: \ \zeta \leq \chi \rangle \in \bfV[\bfG^0_\kappa \ast \bfG^1_\chi]$ such that:
\mn
\begin{enumerate}
	\item[$(*)_{1.1}$]
	\begin{enumerate}[label = (\alph*), ref = \alph*]
		\item \label{C(a)} $q_\zeta \in (\bbP^0_{\bfj(\kappa)} *
		\name \bbP'_{j(\chi)})/\bfG^0_{\kappa +1}$, and if $\varp < \zeta$ then $q_\varp \le q_\zeta$,
		\sn
		\item $q_\zeta \in (\bbP^0_{\bfj(\kappa)} *
		\name \bbP'_{j(\zeta)}) /\bfG^0_{\kappa +1}$, i.e.\ $q_\zeta \um \bfj(\kappa) \Vdash_{\bbP^0_{\bfj(\kappa)}} q_\zeta(\bfj(\kappa)) \in \name \bbP'_{\bfj(\zeta)}$),
		\item \label{C(b)} whenever $p \in \bfG^0_{\kappa+1} \cap (\bbP^0_\kappa \ast \name{\bbP}^1_\zeta)$ then $\bfj(p \um \{\kappa\}) =  \bfj(p) \um \{\bfj(\kappa)\} \le q_\zeta$ (in the quotient forcing $(\bbP^0_{\bfj(\kappa)} *
		\name \bbP'_{j(\chi)})/\bfG^0_{\kappa +1}$),
		\sn
		\item \label{C(c)} whenever $\name A$ is a $\bbP^0_\kappa * \name{\bbP}^1_\zeta$-name
		of a subset of $\kappa$ (so $\bfj(\name A)$ is a $\bbP^0_{\bfj(\kappa)} * \name{\bbP}'_{\bfj(\zeta)}$-name for a subset of $\bfj(\kappa)$) then 
		$$q_\zeta \parallel_{(\bbP^0_{\bfj(\kappa)} \ast \name \bbP'_{\bfj(\zeta)})/ \bfG^0_{\kappa+1}} \ \kappa \in \bfj(\name A).$$
		\sn
		
		\item \label{C(d)} if $\zeta \in S^*$ for $D_\zeta = \name D_\zeta[\bfG^1_{\zeta}]$ (from \ref{cf0} of \ref{d})  we have: If $D_\zeta$ 
		generates a $\kappa$-complete filter on $V_\kappa$ (in $\bfV_1[\bfG^1_{\zeta}]= \bfV_0[\bfG^0_{\kappa}*\bfG^1_{\zeta}]$)
		 \then  \ 
		 we have (that $q_\zeta$ forces)
		\beeq \label{C(e)k} \left(q_{\zeta +1}(\bfj(\kappa)\right)(\bfj(\zeta)) \geq  \left(\name \varp_\zeta \cup \left( \name A'_{\bfj(\zeta)} \um \{ \kappa \} \right), \name A'_{\bfj(\zeta)} \um (\kappa+1, \bfj(\kappa))  \right). \eeq
	\end{enumerate}
	(In the proof of only \ref{c10} $D_\zeta$'s for that $D_\zeta \subseteq cP(\kappa)$ are relevant, and it is enough to ensure that if for each $A \in D_\zeta$ we have $\kappa \in A$ (forced by $q_\zeta$), then $q_{\zeta+1} \Vdash  \ ``\kappa \in \bfj(\varp_\zeta)''$.)
\end{enumerate}
\end{claim}
\begin{PROOF}{Claim \ref{cind}}
	Working in $\bfV_2 = \bfV_0[\bfG_{\kappa+1}]$ we can define the $q_\eta$'s ($\eta < \chi$, $q_\eta \in (\bbP^0_{\bfj(\kappa)} *
	\bbP'_{\bfj(\eta)})/\bfG^0_{\kappa +1}$) by induction on $\eta$.
	Assume that $q_\xi$'s ($\xi < \eta$) are chosen and $\eqref{C(a)}- \eqref{C(d)}$ hold. First we choose $q_\xi'$ satisfying $\eqref{C(a)}$, $\eqref{C(b)}$, $\eqref{C(d)}$ which we will then further strengthen.
	
	Let $q'_0 \in (\bbP^0_{\bfj(\kappa)} *
	\name \bbP'_{\bfj(0)})/\bfG^0_{\kappa +1}=\bbP^0_{\bfj(\kappa)}/\bfG^0_{\kappa +1}$ be the empty condition. For $\eta$ limit we choose $q'_\eta \in (\bbP^0_{\bfj(\kappa)} *
	\name  \bbP'_{\bfj``\eta})/\bfG^0_{\kappa +1}$ to be an upper bound of the increasing sequence $\langle q_\xi: \ \xi < \eta \rangle$ satisfying $\eqref{C(b)}$. Now  it is easy to see that $\eqref{C(b)}$ holds for $q'_\eta$, even if $\bbP^0_\kappa \ast \bbP^1_\eta$ is bigger than the direct limit of $\bbP^0_\kappa \ast \name \bbP^1_\xi$'s ($\xi < \eta$); also recall Fact $\ref{Cf}$.	
	 
	 If $\eta = \xi+1 $ is a successor and
	 \begin{itemize} \item if $\xi \notin S^*$, \end{itemize} \noindent then using simply the $<$$(2^{\chi})^+$-directed closedness of $\bbP^0_{\bfj(\kappa)} *
	 \bbP'_{\bfj(\chi)})/\bfG^0_{\kappa +1}$ (by Fact $\ref{Cf}$)  define $q'_\eta \in (\bbP^0_{\bfj(\kappa)} *
	\bbP'_{\bfj(\xi+1)})/\bfG^0_{\kappa +1}$ to be an upper bound of $q_\xi \in \bbP^0_{\kappa} *
	\bbP^1_{\xi} $ and the set $\{\bfj(p): \ p \in (\bbP^0_{\kappa} *
	\bbP^1_{\xi+1}) \cap \bfG^0_{\kappa+1}  \}$.
	
	 Otherwise, 
	 \begin{itemize} \item if $\xi \in S^*$, \end{itemize} (where $\eta = \xi+1$) then recall that by the definition of $\name{\bbQ}^1_{\xi+1}$ each $p \in (\bbP^0_{\kappa} *
\name	\bbP^1_{(\xi+1)})$ the coordinate $p(\xi+1)$ is a $(\bbP^0_{\kappa} *
\name	\bbP^1_{\xi})$-name for a pair $(\varp, A)$ with $\varp = \varp \um (0,\gamma)$ for some $\gamma < \kappa$ and $A \seq V_\kappa^{\bfV_0[\bfG^0_{\kappa} \ast \bfG^1_\xi]}$, $A = A \um [\gamma,\kappa)$. As the fact that $\name D_\xi$ generates a $<$$\kappa$-closed filter on $V_\kappa$ implies that $\bfj(\name D_\xi)$ generates a $<$$\bfj(\kappa)$-closed filter in $V_{\bfj(\kappa)}$, we can choose $q'_{\xi+1}$ so that  $q'_{\xi+1}(\bfj(\xi))$ satisfies $\eqref{C(e)k}$ (with $\zeta = \xi$), hence $\eqref{C(d)}$ (with $\zeta = \xi+1 = \eta$) as well.
	
	Finally, for $\eqref{C(c)}$, first note that we can assume $\name A \in \cN_\eta$, so there are at most $\chi$-many of them. Now  choosing an increasing sequence of conditions $\langle q''_\gamma : \gamma < \chi \rangle$ in 
		$(\bbP^0_{\bfj(\kappa)} *
	\name	\bbP'_{\bfj(\eta)})/\bfG^0_{\kappa +1}$  with $q''_0 = q'_\eta$, we can decide for each name $\name{X}$ the statement $\kappa \in \bfj(\name{X})$. So using the $<\chi^+$-directed closedness of 	$(\bbP^0_{\bfj(\kappa)} *
		\name \bbP'_{\bfj(\eta)})/\bfG^0_{\kappa +1}$ in $\bfV_0[\bfG^0_{\kappa +1}]$ (Fact $\ref{Cf}$), we can choose $q_\eta$ to  be an upper bound of the sequence $\langle q''_\gamma : \gamma < \chi \rangle$, yielding $\eqref{C(c)}$ as desired.

	Finally, $q_\chi$ is defined to be an upper bound of the $q_\eta$'s ($\eta < \chi$).

\end{PROOF}

\bfact
By the definition of $\bbP^0_\kappa \ast \name \bbP^1_\chi$, and the way we $q_\chi$ we constructed we have: 
\begin{enumerate}[label = $(\bowtie)_{\arabic*}$, ref = $(\bowtie)_{\arabic*}$]
	
	\setcounter{enumi}{\value{btc}}
	\item  For each  $\delta \in S^*$ if $D_\delta$ generates a $\kappa$-complete ultrafilter on $V_\kappa$, then 
	\[ \Vdash_{\bbP^0_{\kappa+1}} \ \forall \name A \in \name D_\delta \ \exists \alpha < \kappa \text{ s.t. } (\name \varp_\delta \um (\alpha, \kappa) \seq \name A), \]
	\stepcounter{btc}
	\item moreover, (in $\bfM[\bfG^0_{\kappa+1}]$)
	\[ q_\chi \Vdash_{(\bbP^0_{\bfj(\kappa)}\ast \name \bbP'_{\bfj(\chi)})/ \bfG^0_{\kappa+1}} \forall d \ \left(\kappa(d) = \kappa \ \wedge \ d \in \bigcap_{\name A \in \name D_\delta} \bfj(\name A)\right) \rightarrow \left(d \in \bfj(\name \varp_\delta)\right), \]
	
	in particular, this  defines the normal ultrafilter
		\stepcounter{btc}
\end{enumerate}
\begin{enumerate}[label = $(\bullet)_{\arabic*}$, ref=$(\bullet)_{\arabic*}$]
	\setcounter{enumi}{\value{bcounter}}
	\item $ D^\bullet = \{ \name
	A[\bfG^0_{\kappa+1}]: \ \Vdash_{\bbP^0_{\kappa+1}} \name A \seq \kappa, \ q_\chi \Vdash ``\kappa
	\in \bfj(\name A)"\}, $
	\stepcounter{bcounter}
\end{enumerate}
\begin{enumerate}[label = $(\bowtie)_{\arabic*}$, ref = $(\bowtie)_{\arabic*}$]
	
	\setcounter{enumi}{\value{btc}}
	\item  if $\delta \in S^*$ is such that $D_\delta \subseteq D^\bullet$, then $\varp_\delta$ is  the pseudointersection of $D_\delta$, moreover, $\varp_\delta \in D^\bullet$.
	\stepcounter{btc}
\end{enumerate}

\efact
\noindent This together with \ref{cf0} complete the proof of $\eqref{C}$$\eqref{Prikry}$.

\noindent

\underline{Case 2}:  For \ref{c10}(C)(b) we proceed as follows.
\mn
In $\bfV_1^{\bbP^1_\chi}$ we have to find a sequence $\bar U = \langle U_\alpha: \ \alpha < \kappa \rangle$ of normal measures on $\kappa$ increasing in the Mitchell order, such that each $U_\alpha$ satisfies our closedness properties, namely, whenever $\langle X_\nu : \nu < \lambda \rangle$ is a sequence in $U_\alpha$,  there exists $X' \in U_\alpha$, $|X' \setminus X_\nu| < \kappa$ for each $\nu < \lambda$.  Let $U_0$ be the normal ultrafilter using Case 1, i.e.\ $\eqref{C} \eqref{Prikry}$.

Working in $\bfV_1[\bfG^1_\chi] = \bfV_0[\bfG^0_\kappa \ast \bfG^1_\chi]$ we will construct the sequence by induction, so fix $\alpha < \kappa$, and assume that $U_\beta$'s are already defined for $\beta < \alpha$. So 
\begin{enumerate}[label = $(\bullet)_{\arabic*}$, ref=$(\bullet)_{\arabic*}$]
	\setcounter{enumi}{\value{bcounter}}
	\item let $\overline{\name U}$ be a $\bbP^0_\kappa \ast \name\bbP^1_\chi = \bbP^0_{\kappa+1}$-name for $\langle U_\beta: \ \beta < \alpha \rangle \in \bfV_0[\bfG^0_\kappa \ast \bfG^1_\chi]$,
	where $1_{\bbP^0_{\kappa+1}}$ forces that $\overline{\name U} = \langle \name U_\beta: \ \beta < \alpha \rangle$ is an increasing sequence of $\kappa$-complete normal ultrafilters w.r.t. the Mitchell-order of length $\alpha$, each $\name U_\beta$ is $<$$\lambda^+$-directed modulo $[\kappa]^{<\kappa}$.
	\stepcounter{bcounter}
\end{enumerate}  and fix an elementary embedding $\bfj_*: \bfV_0 \to M_*$ with critical point $\kappa$, $^{\chi} M_* \subseteq M_*$ with 
\beeq \label{j*} \bfj_*(\bfh)(\kappa) = \left\langle  \name{\bbP}^1_\chi, \chi^+, \name U \right\rangle \eeq
(recall the definition of $\bfh$ \ref{bfhd}, this is possible).

Defining $\name \bbP'_* = \bfj_*(\name \bbP^1)$, and letting $(\bbP^0_*)_{\bfj_*(\kappa)} = \bfj_*(\bbP^0_\kappa)$ observe that by the definition of $\bbP^0_\kappa$ (Definition $\ref{P0def}$)
\[ \bfj_*(\bbP^0_\kappa \ast \name \bbP^1_\chi) = (\bbP^0_*)_{j_*(\kappa)} \ast \name (\bbP'_*)_{\bfj_*(\chi)},\] 
and
\[ (\bbP^0_*)_{\kappa+1} = \bbP^0_\kappa \ast \name \bbP^1_{\chi}. \]
Now our fixed $\bfG^0_{\kappa+1} \subseteq \bbP^0_{\kappa+1}$ is generic over $\bfV_0$ and also over $\bfM_*$.

   With a slight abuse of notation (in the proof of \underline{Case 2} from now on, in order to avoid notational awkwardness) we will refer to $(\bbP^0_*)_{\bfj_*(\kappa)}$ as $\bbP^0_{\bfj_*(\kappa)}$, and to $(\name \bbP'_*)_{\bfj_*(\chi)}$ as $\bbP'_{\bfj_*(\chi)}$; moreover, observe that all the preceding facts and claims hold in this setting, we only used that $\bfj(\bfh(\kappa)) = \langle \name \bbP^1_\chi, \chi^+, \name x \rangle$ for some name $\name x$, which obviously holds for $\bfj_*$ as well.
 In this new setting we appeal to Claim $\ref{cind}$, obtaining the condition $q^*_\chi \in \bbP^0_{\bfj_*(\kappa)+1} / \bfG^0_{\kappa+1}$, and the $\kappa$-complete normal ultrafilter
 \beeq \label{Dufdf}  D^\bullet_* = \{ \name A[\bfG^0_{\kappa+1}]: \ \bfM_*[\bfG^0_{\kappa+1}] \models \ `` q^*_\chi \Vdash_{\bbP^0_{\bfj_*(\kappa)} \ast \bbP'_{\bfj_*(\chi)} / \bfG^0_{\kappa+1}} \kappa \in \bfj_*(\name A) \} \text{"} \eeq
(which is $\kappa$-complete normal ultrafilter over $\bfV_0[\bfG^0_{\kappa+1}]$, belonging to $\bfV_0[\bfG^0_{\kappa+1}]$) and $<\lambda^+$-directed w.r.t.\ $\supseteq^*$.
We only need to prove the following claim, implying that the filter $D_*^\bullet$ dominates $\{ U_\beta: \beta < \alpha\}$ in the Mitchell order.
\bcl 
For each $\beta < \alpha$ there exists a sequence $\langle W_\gamma: \ \gamma < \kappa \rangle \in \bfV_0[\bfG^0_{\kappa+1}]$, where 
\begin{itemize}
	\item for $D_*^\bullet$-many $\gamma < \kappa$ the set $W_\gamma$ is an ultrafilter over $\gamma$,
	\item for each $X \in \cP(\kappa) \cap \bfV_0[\bfG^0_{\kappa+1}]$
	\[ X \in U_\beta \ \iff \ \{ \gamma < \kappa: \ (X \cap \gamma) \in W_\gamma \} \in D_*^\bullet. \]
\end{itemize}
\ecl
\begin{PROOF}{}
Using $\eqref{j*}$
\[ \left\{ \begin{array}{rl} \gamma < \kappa: & \ \bfh(\gamma)= \langle \name x_\gamma, \alpha_\gamma,\name y_\alpha \rangle, \text{ where } \name y_\alpha \text{ is a } \bbP^0_{\gamma+1}\text{-name} \\ 
	& \text{for a sequence of subsets of } \cP(\gamma) \text{ (of length } \alpha),  \\
	&  \name x_\alpha = \name \bbQ^0_\alpha
	\end{array} \right\} \in D_*^\bullet,\]
so we can fix $Y \in D_*^\bullet \cap \bfV_0$, and the sequence $\langle \name W_\gamma: \ \gamma < \kappa \rangle$ such that 
\begin{enumerate}[ref=$\blacktriangle_\arabic*$, label = ($\blacktriangle_\arabic*$)]
	\item for each $\gamma \in Y$ $\name W_\gamma$ is a $\bbP^0_{\gamma+1}$-name for a subset of $\cP(\gamma)$,
	\item \label{ub} $\bfj_*(\langle \name W_\gamma: \ \gamma < \kappa \rangle)(\kappa) = \name U_\beta$.
\end{enumerate}
We will prove that $W_\gamma = \name W_\gamma [\bfG^0_{\kappa+1}]$ ($\gamma < \kappa$) works.

For a fixed $\bbP^0_\kappa \ast \bbP^1_\chi$-name $\name X \in \bfV_0$ (for a subset of $\kappa$)
define the $\bbP^0_\kappa \ast \bbP^1_\chi$-name $\name Z_X \in \bfV_0$ as follows.
\beeq \label{Z} 1_{\bbP^0_\kappa \ast \bbP^1_\chi} \Vdash  \name Z_X = \{ \gamma < \kappa: \ \name X \upharpoonright \gamma \in \name W_\gamma\}, \eeq
We only have to verify that
\beeq \label{miss} \name X [\bfG^0_{\kappa+1}] \in \name U_\beta[\bfG_{\kappa+1}^0] \ \text{ iff } \  \name Z_X[\bfG^0_{\kappa+1}] \in D_*^\bullet. \eeq
But the latter is defined (by \eqref{Dufdf}) as 
\[ \begin{array}{l}
\name Z_X[\bfG^0_{\kappa+1}] \in D_*^\bullet, \\
 \ \ \Updownarrow \\ (\text{in } \bfM_*[\bfG_{\kappa+1}^0]) \ q^*_\chi \Vdash_{\bbP^0_{\bfj_*(\kappa)}\ast \bbP'_{\bfj_*(\chi)} / \bfG^0_{\kappa+1}} \kappa \in \bfj_*(\name Z_X),\\
 \ \ 
\end{array} 
\]

Therefore, as $\bfj_*(\name W)(\kappa) = \name U_\beta$ by $\eqref{Z}$, for $\eqref{miss}$ it suffices to show
\beeq  \label{m2}\name X [\bfG^0_{\kappa+1}] \in \name U_\beta[\bfG_{\kappa+1}^0] \ \text{ iff } \ q^*_\chi \Vdash \bfj_*(\name X) \upharpoonright \kappa \in \name U_\beta. \eeq
But then by the elementarity of $\bfj_*$ (and $\crr(\bfj_*) = \kappa$)
\[ \forall \alpha < \kappa, \forall p \in \bbP^0_{\kappa} \ast \name \bbP^1_\chi: \ \ p \Vdash_{\bbP^0_{\kappa} \ast \name \bbP^1_\chi} \check{\alpha} \in \name X \ \iff \ \bfj_*(p) \Vdash_{\bbP^0_{\bfj_*(\kappa)} \ast \name \bbP'_{\bfj_*(\bfj_*(\chi))}} \check{\alpha} \in \bfj_*(\name X), \]
 and recalling $p \in \bfG_{\kappa+1}^0$ implies $q^*_\chi \geq \bfj_*(p)$) we get that 
\begin{enumerate}
	\item[$(*)_1$] $q^*_\chi$ forces $\bfj_*(\name X) \upharpoonright \kappa$ to be equal to $\name X [\bfG^0_{\kappa+1}]$.
\end{enumerate}
This  yields $\eqref{m2}$, completing the proof of \underline{Case 2}.
 
\end{PROOF}

\underline{Case 3}:  For \ref{c10}(C)(c). 
First we redefine the elementary embedding $\bfj$ from Definition $\ref{1def}$ (and as well $\bbP^0_{\bfj(\kappa)}$, $\name \bbP'_{\bfj(\chi)}$):
\setcounter{bcounter}{1}
\bdd {\ }
\begin{enumerate}[label = $(\bullet)_{\arabic*}$, ref=$(\bullet)_{\arabic*}$]
	\setcounter{enumi}{\value{bcounter}}
	\item Let $\rho = |2^{(\Upsilon\cdot \chi)^\kappa}+\eta|$, and
	\item \label{jM'} define $\bfj:\bfV_0 \rightarrow \bfM$ be an elementary embedding with
	critical point $\kappa$ such that $(\bfj(\bfh))(\kappa) = \langle \name{\bbP}^1_\chi, \rho^+, \check{\emptyset} \rangle$ ($\check{\emptyset} = \emptyset$ is the canonical name for the empty set)
	and $\bfj(\kappa)> \rho$, $^{\rho}\bfM \subseteq \bfM$, 
	\item Let $\langle \bbP^0_\alpha,\name{\bbQ}^0_\beta:\alpha \le
	\bfj(\kappa),\beta < \bfj(\kappa)\rangle = \bfj(\langle
	\bbP^0_\alpha,\bbQ^0_\beta:\alpha \le \kappa,\beta < \kappa\rangle)$
	so $\name{\bbQ}^0_\kappa = \name{\bbP}^1_\chi$, and let $\name \bbP'_{\bfj(\chi)} = \bfj(\name \bbP^1_\chi)$.
\end{enumerate}
\edd
Similarly as in Facts $\ref{mf}$, $\ref{Cf}$, $\ref{bt1}$ we can get the following.

\bfact \label{mfcase3}
The filter $\bfG^0_{\kappa +1}$ is generic over $\bfM$ as well, and the forcing notions
$\bbP^0_{j(\kappa)}/ \bfG^0_{\kappa +1}$  and $(\bbP^0_{j(\kappa)} *
\name{\bbP}'_\gamma)/\bfG^0_{\kappa +1}$ ($\gamma \le \bfj(\chi)$) are well defined and
$<|2^\Upsilon+\eta|^+$-directed closed in $\bfM[\bfG^0_{\kappa +1}]$.
\efact

\bfact \label{Cfcase3} $\bfV_1 \models \ \text{``}(\bbP^0_{j(\kappa)} * \name{\bbP}'_\gamma)/\bfG^0_{\kappa +1} \text{ is } <|2^\Upsilon+\eta|^+\text{-directed closed."}$
\efact 
Fact $\ref{Cf}$ follows from the fact below.

\bfact	 \label{bt1case3}
$\bfV[\bfG^0_\kappa \ast \bfG^1_\chi] \models \ ``\ ^{2^\Upsilon+\eta}\bfM[\bfG^0_{\kappa +1}] \subseteq \bfM[\bfG^0_{\kappa +1}]"$.
\efact 	
 
 Using this new $\bfj$ we will extract the ultrafilter $W \subseteq \cP([\Upsilon]^{<\kappa})$ (in the sense of $\bfV_0[\bfG_{\kappa+1}^0]$), and the sequence of ultrafilters $\overline{U}$ as well from the information provided by $\bfG_{\kappa+1}^0 = \bfG_\kappa^0 \ast \bfG_\chi^1$, and $q_\chi \in (\bbP^0_{\bfj(\kappa)} *
\bbP'_{j(\chi)})/\bfG^0_{\kappa +1}$ (given by Claim $\ref{cind}$), and then we will prove that it is indeed a measure sequence corresponding to the elementary embedding $\bfj_W$.
Obviously
\begin{enumerate}[label = $(\circledcirc_1)$, ref = $\circledcirc_1$ ]
	\item \label{cc0}
	 $\bfj(\kappa)> \chi$, $^\chi  M \subseteq M$.
\end{enumerate}
Observe that Claim $\ref{cind}$ is true in this setting as
 and let the master condition $q_\chi \in (\bbP^0_{\bfj(\kappa)} *
\bbP'_{j(\chi)})/\bfG^0_{\kappa +1}$ be given by that.
First we claim that by possibly extending $q_\chi$, we can assume that
\begin{enumerate}[label = $(\circledcirc_2)$, ref = $\circledcirc_2$ ]
	\item \label{cc}
	For each $A \in \cP([\Upsilon]^{<\kappa}) \cap \bfV_2$ the condition $q_\chi \in (\bbP^0_{\bfj(\kappa)} *
	\bbP'_{\bfj(\chi)})/\bfG^0_{\kappa +1}$ decides about ``$\bfj``\Upsilon \in \bfj(A)$" (in $\bfM[\bfG_{\kappa+1}^0]$). 
\end{enumerate}
For this first we count the possible $A$'s. Recall that $\bbP^1_\chi$ is $<\kappa$-closed ($\eqref{B}$/ $\eqref{(c)}$ 
$$[\chi]^{<\kappa} \cap \bfV_2 = [\chi]^{<\kappa} \cap \bfV_1 = [\chi]^{<\kappa} \cap \bfV_0[\bfG_\kappa^0],$$
and as $|\bbP^0_\kappa|= \kappa$, 
\beeq 
	|[\Upsilon]^{<\kappa} \cap \bfV_2| \leq (\Upsilon \cdot \chi)^\kappa.
\eeq
Second, as $|\bbP_\kappa^0 \ast \bbP^1_\chi| = \chi$, we have
\beeq \label{chikap}  \bfV_2 = \bfV_0[\bfG^0_\kappa \ast \bfG^1_\chi] \models  \cP([\chi]^{<\kappa})| \leq (2^{(\chi \cdot \Upsilon)^\kappa})^{\bfV_0} \leq \rho. \eeq
Now using Fact $\ref{Cf}$  we can extend $q_\chi$ to another condition $q_*$ in (at most) $\rho$-many steps (in $(\bbP^0_{\bfj(\kappa)} *
\bbP'_{j(\zeta)})/\bfG^0_{\kappa +1}$) so that 
\begin{enumerate}[label = $(\circledcirc_3)$, ref = $\circledcirc_3$] 
	\item \label{eldont}for each name $\name A$ for a subset of $[\chi]^{<\kappa}$ \[ \bfM[\bfG^0_{\kappa+1}] \models  \ q_* \parallel \ \bfj``\Upsilon \in \name A, \]
\end{enumerate}

and so (by possibly replacing $q_\chi$ by $q_*$) $\eqref{cc}$ holds, indeed.
Now we can  define the $\kappa$-complete, fine, normal ultrafilter
\beeq \label{Wdf} W = \{ \name A[\bfG^0_\kappa \ast \bfG^1_\kappa] \in [\Upsilon]^{<\kappa}: \ q_\chi \Vdash \bfj``\chi \in \bfj(A) \} \in \bfV_2, \eeq
Now let $\bfj_{W}: \bfV_2 \to \bfM_W = \Mos(^{[\chi]^{<\kappa}} \bfV_2 /W)$ be the corresponding elementary embedding, and let $\overline{U} = \langle U_\alpha: \ \alpha < \dom(\overline{U}) \rangle$ be the ultrafilter sequence of maximal length associated to $\bfj_W$, that is, the following holds in $\bfV_2$.
\begin{enumerate}[label = $(\boxminus_\arabic*)$, ref = $\boxminus_\arabic*$]
	\item $U_0 = \kappa$, and for each $\alpha \in \dom(\overline{U})$, $\alpha > 0$ the set $U_\alpha \subseteq \cP(V_\kappa)$ is a $\kappa$-complete normal ultrafilter satisfying
	\[ \forall A \subseteq V_\kappa: \ A \in U_\alpha \ \iff \ U \upharpoonright \alpha \in \bfj_W(A) \]
	(therefore for each $\alpha < \dom(\overline{U})$ $U \upharpoonright \alpha \in\bfM_W$),
	\item $\overline{U} \notin \bfM_W$.
\end{enumerate}
The following two claims complete the proof of $\ref{c10}\eqref{C}\eqref{Rad}$ as by our assumptions $\dom(\overline{U})= \eta \leq \rho$.
\bcl \label{zc1} For every ultrafilter sequence 
$\overline{F} \in  \bfM_W$ with $\kappa(\overline{F}) = \kappa$ there exists an ultrafilter sequence $\overline{F}' \in \bfM[\bfG_{\kappa+1}^0]$ with $\kappa(\overline{F}') = \kappa$ 
such that  for each name $\name A$ for a subset of $V^{\bfV_0[\bfG_\kappa^0 \ast \bfG^1_\chi]}_\kappa$ we have
\[ \overline{F} \in \bfj_W(\name A[\bfG_\kappa^0 \ast \bfG_\chi^1]) \ \iff \ \bfM[\bfG^0_{\kappa+1}]\models \ q_\chi \Vdash \overline{F}' \in \bfj(\name A).  \]
\ecl

\bcl \label{zc2}
For every set of ($\bbP_{\bfj(\kappa)} \ast \bbP'_{\bfj(\chi)}$-names for) ultrafilter sequences 
$\{ \overline{\name F_i}: i<\sigma\} \subseteq  \bfM$ with $\kappa(\name F_i) = \kappa$ ($i < \sigma$) if the filter 
\[ F_*= \bigcap_{i < \sigma} \{A \subseteq V_\kappa^{\bfV_2}: \ q_\chi \Vdash \overline{\name F_i} \in \bfj(A) \} \]
satisfies $(\forall \alpha < \kappa): |\cup F_* \um \alpha| < \kappa$, then $F_*$
 is $<\lambda^+$-directed in the sense that
for any system $\langle X_\alpha: \alpha < \lambda \rangle$ in $F_*$ there is a set $X'   \in F_*$ s.t.\ for each $\alpha < \lambda$ there exists  $\delta < \kappa$ with $X' \upharpoonright [\delta, \kappa) \subseteq X_\alpha$.
\ecl
\begin{PROOF}{Claim \ref{zc1}} Instead of factoring through our elementary embeddings (after forcing) we provide a direct calculation.
Fix the ultrafilter sequence $F \in \bfM_W$, and pick a function 
$f \in \bfV_2$, $\dom(f) = [\Upsilon]^{<\kappa}$, $\bfj_W(f)(\bfj_W``\Upsilon) = F$, where we can assume that 
\beeq \label{ufs} \forall x \in \dom(f) \ f(x) \text{ is an u.f. sequence with } \kappa(f(x)) = \otp(\kappa \cap x) .\eeq
Now we can fix a $\bbP^0_\kappa \ast \bbP^1_\chi$-name $\name f \in \bfV_0$ of $f$, such that $1_{\bbP^0_\kappa \ast \bbP^1_\chi}$ forces $\eqref{ufs}$.
Now as in $V_0$ $\name f$ is a $\bbP^0_\kappa \ast \bbP^1_\chi$-name for a function with $\dom(f) = [\Upsilon]^{<\kappa}$, by elementarity $\bfj(\name f)$ is a $\bbP^0_{\bfj(\kappa)} \ast \bbP'_{\bfj(\chi)}$-name for a function with domain $[\bfj(\Upsilon)]^{<\bfj(\kappa)}$, thus  there is name $\name F' \in \bfM$ such that 
\beeq \label{F'} \bfM \models \ \ \Vdash_{\bbP^1_{\bfj(\kappa)} \ast \bbP'_{\bfj(\chi)}} \bfj(\name f)(\bfj``\Upsilon) =\name F'. \eeq

It is only left to check that for each $X \subseteq V_\kappa^{\bfV_2}$ the conditions
"$F \in \bfj_W(X)$" and "$\name F' \in \bfj(X)$" are equivalent. More precisely, we prove the following.
\begin{enumerate}[label = $(\circ)$, ref = $\circ$]
	\item \label{ekv} For every fixed $\bbP^0_\kappa \ast \bbP^1_\chi$-name $\name X$ for a subset of $V^{\bfV_2}_\kappa$ 
	\[ F \in \bfj_W(\name X[\bfG_\kappa^0 \ast \bfG^1_\chi]) \iff q_\chi \Vdash \name F' \in \bfj(\name X). \]
\end{enumerate}
As $F = \bfj_W(f)(\bfj_W``\Upsilon)$ we can reformulate the lhs.\ as the statement \[\bfV[\bfG_\kappa^0 \ast \bfG_\chi^1] \models \ \{y \in [\Upsilon]^{<\kappa}: \ f(y) \in X  \} \in W, \] 
i.e. for some $p \in \bfV_0[\bfG_\kappa^0 \ast \bfG_\chi^1]$
\[ p \Vdash_{\bbP^0_\kappa \ast \bbP^1_\chi} \{y \in [\Upsilon]^{<\kappa}: \ \name f(y) \in \name X  \} \in \name W.  \]
Now for the the $\bbP^0_\kappa \ast \bbP^1_\chi$-name $\name C:=\{y \in [\Upsilon]^{<\kappa}: \ \name f(y) \in \name X  \}$ we have (by $\eqref{cc}$ and $\eqref{Wdf}$)
\[ \name C[\bfG_\kappa^0\ast \bfG_\chi^1] \in W \ \iff \ q_\chi \Vdash \bfj``\Upsilon \in \bfj(\name C). \]
(Recall that $q_\chi$ decides this by $\eqref{eldont}$ as $\name C$ is a name for a subset of $[\Upsilon]^{<\kappa}$.) This latter is equivalent to 
\[ q_\chi \Vdash \bfj(\name f)(\bfj``\Upsilon) \in \bfj(\name X),\]
so recalling $\Vdash \bfj(\name f)(\bfj``\Upsilon) = F'$ by $\eqref{F'}$
this is clearly equivalent to $q_\chi \Vdash \name F' \in \bfj(\name X)$, therefore  $\eqref{ekv}$ holds, as desired.
\end{PROOF}

\begin{PROOF}{Claim \ref{zc2}}
Fix $\langle \name F'_i: i < \sigma \rangle$ given by Claim $\ref{zc1}$

We only have to recall how we constructed $q_\chi$, which ensures the existence of the desired pseudointersection. Fix a sequence $\langle X_\alpha: \alpha < \lambda \rangle$ in the filter $F_*$. Now let $D' = \{ X_\alpha: \ \alpha < \lambda\}$, which is equal to $D_\zeta$ for some $\zeta < \chi$ by \ref{cf0} from our assumptions \ref{B}/\ref{d}. Finally, recalling Definition $\ref{fdf}$ and $\eqref{C(e)k}$ from Claim $\ref{cind}$ we get that for the generic sequence $\varepsilon_\zeta$ (which is a pseudointersection of the $D' = D_\zeta$) \[ q_{\zeta+1} \Vdash \bfj(\varepsilon_\zeta) \um (\kappa+1) = \varepsilon_\zeta \cup (A'_{\bfj(\zeta)} \um [\kappa, \kappa+1)),\] which means that by the definition of $A'_{\bfj(\zeta)}$ (Definition $\ref{fdf}$)
	\[ \forall i < \sigma \ (\forall X \in D_\zeta \ q_\chi \Vdash \name F_i \in \bfj(X)) \Rightarrow (q_\chi \Vdash \overline{\name F_i} \in A'_{\bfj(\zeta)} )  \]

\end{PROOF}

\end{PROOF}

\subsection {The preliminary forcing for obtaining $(\kappa,\lambda)-1$ systems together with a universal in $(K_{\kappa})_\lambda$} \label{3B}\

This subsection deals with the application of Claim $\ref{c10}$, we show that it is possible to force a universal object in $(K_\kappa)_\lambda$ with a notion of forcing satisfying requirements from Claim $\ref{c10}$.
\bigskip
\begin{conclusion}
	\label{h28}
	Assume $\kappa$ is supercompact
	 $\kappa < \lambda < \chi =
	\chi^\lambda$, $\lambda$ is regular, $(\forall \theta)(\theta \in \card \wedge \kappa \le \theta
	< \lambda \Rightarrow 2^\theta = \theta^+)$ 
	and $\sigma = \cf(\sigma)
	< \kappa$.
	
	\noindent
	Then in some forcing extension $\bfV^{\bbP}$ preserving cardinals and preserving
	cofinalities $> \kappa$ and in $\bfV^{\bbP},2^\kappa = \chi,\kappa$
	strong limit singular of cofinality $\sigma$ and there is a universal graph in cardinality $\lambda$.
\end{conclusion}

\begin{PROOF}{\ref{h28}}
	We shall use \ref{h2}, but we have to justify it.
	That is, we need a forcing fitting in the scheme in Claim $\ref{c10}$ with $\bfV_0 = \bfV$, specifying the $(<
	\kappa)$-directed-complete iteration   $\bbP^1_\chi = \langle \bbP^1_\alpha, \bbQ^1_\beta: \ \alpha \leq \chi, \beta < \chi \rangle \in \bfV_1 =  \bfV^{\bbP^0_\kappa}$ in which we are free to choose $\name \bbQ_\beta$'s on $\beta'$s outside $S^* \subseteq \chi$. (And then conclusion \ref{C}/ \ref{Prikry} or \ref{Mag} with Claim \ref{h8} together with Claim $\ref{b2}$ or $\ref{b23}$ give the desired consistency result.) Our task is to construct (in $\bfV_1$) a suitable iteration $\bbP^1_\chi$, and checking that $\bbP^1_\chi$ is
	\begin{enumerate}[label = $(\intercal)_{\arabic*}$, ref =  $(\intercal)_{\arabic*}$]
		\item \label{e1}  $<$$\kappa$-directed closed,
		\item \label{e2}  is of  cardinality
		$\chi$ (up to equivalence),
		\item \label{e3}  has the $\lambda^+$-c.c.,
		\item  \label{e4} not collapsing cardinals, and
		\item  \label{e5}	$\bfV_1 \models \ \ \Vdash_{\bbP^1_\chi}$ ``there is a universal graph in
		$(K_{\kappa})_\lambda$",
		\item \label{eu} and we can choose $S^* \in [\chi \setminus \{0,1\}]^{\chi}$, $S^* \in \bfV_1$, $|\chi \setminus S^*| = \chi$,  and the $\bbP^1_\delta$-names $\name D_\delta$ ($\delta \in S^*$) satisfying \eqref{B}\eqref{d} from Claim $\ref{c10}$.
	\end{enumerate}
	 
  We will do the same as  in \cite{Sh:175a}, we define (in $\bfV_1$)
  \begin{enumerate}
	\item $\bbQ^1_0$ to be the forcing of $\chi$-many stationary sets of $\lambda$, any two intersecting in a set of size smaller than $\kappa$,
  	\item $\name \bbQ^1_\beta$ for $\beta \in \chi \se (S^* \cup \{0\})$ the main iteration from \cite{Sh:175a} just with $\kappa$-many colors: forcing a generic random graph, and the embeddings into it with $<\kappa$-support partial functions.
  \end{enumerate}
	We need to check that the iteration is indeed $\lambda^+$-cc, which will be ensured by showing that (in $\bfV_1$) $\bbQ^1_0$ is $\lambda^+$-cc, and in $(\bfV_1)^{\bbQ^1_0}$ the iteration of $\bbQ^1_\alpha$'s ($0< \alpha < \chi$), i.e. $\bbP^1_\chi / \bfG^1_1$ has the $\kappa^+$-cc .\\
	First for future reference we we have to remark that
	 \newcounter{pcounter} \setcounter{pcounter}{0}
	\begin{enumerate}[label = $(*)_{\arabic*}$, ref= $(*)_{\arabic*}$]
		\item \label{apr} in $\bfV_1 = \bfV_0^{\bbP^0_\kappa}$  $\kappa$ is still inaccessible as according to Claim $\ref{c10}$ (\ref{c10A}) $\bbP^0_\kappa$ is an Easton-support iteration with $\bbP^0_\alpha \in V_\kappa^{\bfV_0}$ for each $\alpha < \kappa$. As $|\bbP_\kappa^0| = \kappa$ our cardinal arithmetic assumptions above $\kappa$ are also preserved.
		\stepcounter{pcounter}
	\end{enumerate}

		Working in $\bfV_1$ we define the first step $\bbQ^1_0$ to be $Q(\lambda,\chi,\kappa)$ as in \cite[Sec. 6.]{Ba}, see below \ref{aQd} in Definition $\ref{aD}$.
		\bl \label{Baulem} In $\bfV_1$ there exists a forcing poset $\bbQ^1_0$ that is $<$$\kappa$-directed closed, of power $\chi$, having $\lambda^+$-cc, preserving cardinals from $(\kappa, \lambda]$, and 
		\[ \bfV_1^{\bbQ^1_0} \models \ \exists   \{S_\alpha : \ \alpha < \chi \} \subseteq \cP(\lambda), \text { a system of stationary sets s.t. }
		 \ \forall \alpha < \beta < \chi: \ |S_\alpha \cap S_\beta|<\kappa. \]
		\el
	\begin{PROOF}{Lemma \ref{Baulem}}
		\bdd \label{aD} First we define the following auxiliary posets.
		\newcounter{Qba}
		\setcounter{Qba}{0}
		\begin{enumerate}[label = $(\intercal)_\arabic*$, ref =$(\intercal)_\arabic*$]
			\setcounter{enumi}{\value{Qba}}
			
			\item For a regular cardinal $\mu$ we let $Q'(\lambda,\chi,\mu)$ be the set of functions $f$ satisfying
			\begin{enumerate}[label = $(\roman*)$, ref =$(\roman*)$]  
				\item $\dom(f) \in [\chi]^{<\mu}$,
				\item for each $\alpha \in \dom(f)$ $f(\alpha) \in [\lambda]^{<\mu}$,
			\end{enumerate}
			
			with $f \leq g$, iff
			\begin{enumerate}[label = $(\roman*)$, ref =$(\roman*)$]
				\setcounter{enumii}{2}  
				\item $\dom(f) \subseteq \dom(g)$,
				\item $\forall \alpha \in \dom(f)$: $f(\alpha) \subseteq g(\alpha)$,
				\item for each $\alpha \neq \beta \in \dom(f)$ $f(\alpha) \cap f(\beta) = g(\alpha) \cap g(\beta)$.
			\end{enumerate}
			
			\item \label{aQd} Let $Q(\lambda,\chi,\kappa) \subseteq \prod_{\mu \in \reg \cap [\kappa,\lambda]} Q'(\lambda,\chi, \mu)$ be the collection of the following functions $f$
			\begin{enumerate}[label = $(\roman*)$, ref =$(\roman*)$]  
				\item $\forall \mu <\nu \in \reg \cap [\kappa,\lambda]$, $\forall \alpha \in \dom(f_\mu)$: $ f_\mu(\alpha) \subseteq f_{\nu} (\alpha)$
			\end{enumerate}
			with the pointwise ordering inherited from the full product \\ $\prod_{\mu \in \reg \cap [\kappa,\lambda]} Q'(\lambda,\chi, \mu)$.
			
			\stepcounter{Qba} \stepcounter{Qba}	
		\end{enumerate}
		\edd
		\bdd \label{Q10} We let $\bbQ^1_0 = Q(\lambda,\chi,\kappa) \in \bfV_1$.  	
		\edd
		For later reference we note the following. Recall that $\chi^\lambda = \chi$ by our assumptions.
		\begin{observation} \label{|Q_0|}
			For each $\mu \in \reg \cap [\kappa, \lambda]$ $|Q'(\lambda,\chi,\mu)| \leq \chi^{<\mu}\cdot \lambda^{<\mu} = \chi$. Therefore $|\bbQ^1_0| = \chi$.
		\end{observation}
		
		By \cite[Lemma 6.3]{Ba}, recalling $(\sigma \in \card \cap [\kappa, \lambda)) \to (2^\sigma = \sigma^+)$, so $\lambda^{<\lambda} = \lambda$ we have the following.
		\bcl \label{lambdapcc} $Q(\lambda,\chi,\kappa)$ is $\lambda^+$-cc,  $<$$\kappa$-directed closed, preserving cofinalities and cardinals.
		\ecl
		Clearly
		\newcounter{Ssz}
		\setcounter{Ssz}{0}
		\begin{enumerate}[label = $(\ddagger)_\arabic*$, ref =$(\ddagger)_\arabic*$]
			\setcounter{enumi}{\value{Ssz}}
			\item \label{lubi} every directed subset of power less than $\kappa$ in	$\bbQ^1_0 = Q(\lambda,\chi,\kappa)$ has a least upper bound.
			\stepcounter{Ssz}
		\end{enumerate}
		
		Now obviously in $\bfV_1^{\bbQ^1_0}$

		\begin{enumerate}[label = $(\ddagger)_\arabic*$, ref =$(\ddagger)_\arabic*$]
			\setcounter{enumi}{\value{Ssz}}
			\item \label{S'df} the generic subsets $S_\alpha$ ($\alpha < \chi$) defined by	$\Vdash_{\bbQ^1_0} \name S_\alpha = \cup \{ f_\kappa(\alpha): \ f \in \name\bfG\}$ form a $\kappa$-almost disjoint system, i.e. if $\alpha <\beta$, then $\Vdash |\name S_\alpha \cap \name S_\beta| < \kappa$,
			\stepcounter{Ssz}
		\end{enumerate}
		we only need to verify that 
		
		\begin{enumerate}[label = $(\ddagger)_\arabic*$, ref =$(\ddagger)_\arabic*$]
			\setcounter{enumi}{\value{Ssz}}
			\item \label{ssst} for each $\alpha < \chi$ the subset
			\[ S_\alpha  \ \text{ is stationary subset of } \lambda,\]
			\stepcounter{Ssz}
		\end{enumerate}
		which is a standard argument, but for the sake of completeness we elaborate.
		(In fact, recalling \cite[Lemmas 6.3-6.5.]{Ba} with the aid of the following it is easy to argue that  $(S_\alpha \cap E^\lambda_{\geq \kappa})$ i.e.\ restricting $S_\alpha$ to points of cofinality at least $\kappa$ is stationary.)
		\bcl \label{clab}
		The notion of forcing $Q(\lambda, \chi, \kappa)$ is equivalent to the two-step iteration
		$Q(\lambda, \chi,\kappa^+) \ast \name Q'(\lambda,\chi,\kappa, \name F)$ where
		\[ \begin{array}{rl} \bfV_1^{Q(\lambda, \chi,\kappa^+)} \models & \bullet \  F_\alpha \ (\alpha \in \chi) \text { is the generic sequence in } [\lambda]^\lambda, \\
			& \bullet \ Q'(\lambda,\chi,\kappa,  F) \seq Q'(\lambda,\chi,\kappa) \text{ defined by}\\ & \ \ [f \in \name Q'(\lambda,\chi,\kappa, \name F)  \iff  \forall \alpha \in \dom(f) \ f(\alpha) \subseteq \name F_\alpha]. 
		\end{array}
		\]
		Moreover, $Q(\lambda, \chi,\kappa^+)$ is $<$$\kappa^+$-closed, (in $\bfV_1^{Q(\lambda, \chi,\kappa^+)}$), and $Q'(\lambda,\chi,\kappa, F)$ has the $\kappa^+-cc$).
		\ecl
		%
		Looking at the definition of the forcing $Q(\lambda, \chi,\kappa)$ if we are given a condition $p$, and a $Q(\lambda, \chi,\kappa)$-name $\name C_*$ for a club set in $\lambda$ first recall that $Q(\lambda, \chi,\kappa)$ is $<$$\kappa$-closed (Claim $\ref{clab}$), in particular $<\omega_1$-closed, as $\kappa$ is inaccessible. We can define an increasing sequence $p^j $ ($j < \omega$) in $Q(\lambda, \chi,\kappa)$ with $p^0 = p$,  and an increasing sequence of ordinals   $\varrho_j$ ($j < \kappa$) satisfying $p^j \Vdash \varrho_j \in \name C_*$, and if $j < k$, then  $\sup \cup \{ p^j_\lambda(\beta): \  \beta \in \dom(p^j_\lambda) \} < \varrho_k$. This is possible, as $|\dom(p_j)| < \lambda$, as well as $|p^j_\lambda(\beta)| < \lambda$, and $\lambda$ is regular. Then clearly any upper bound of the $p^j$'s forces $\varrho_\omega := \sup\{ \varrho_j: \ j < \omega \} \in \name C_*$,  but as the least upper bound $p^\omega$ does not say anything about the statements $\varrho_\omega \in \name S_\beta$ ($\beta < \chi$) we can extend it to a condition $(p^\omega)'$ with  $\varrho_\omega \in (p^\omega)'_{\mu}(\alpha)$ for each $\mu \in \reg \cap [\kappa, \lambda]$ (thus $(p^\omega)'  \Vdash \varrho_\omega \in \name S_\alpha \cap \name C_*$). 
		This completes the proof of Lemma $\ref{Baulem}$.
	\end{PROOF}
 	 	
 	 As $\bbQ^1_0$ as already defined in Definition $\ref{Q10}$ we can define the iteration $\langle \bbP^1_\alpha, \ \name \bbQ^1_\beta: \ \alpha \leq \chi, \beta < \chi \rangle$ for which we have to define $S^*$.
 	 \bdd We let 
 	 $0,1 \notin S^* \subseteq \chi$ be such that $|S^*| = \chi$, $|\chi \setminus S^*| = \chi$.
 	\edd
 	 \bdd \label{itdf}
 	 	We let $\langle \bbP^1_\alpha, \ \name \bbQ^1_\beta: \ \alpha \leq \chi, \beta < \chi \rangle$ be the following $<\kappa$-support iteration.
 	 	The definition of the $\bbP^1_\beta$-name $\name \bbQ^1_\beta$ goes by induction on $\beta$ as follows, distinguishing three cases. But first
 	 	\begin{enumerate}[label = $\circledast$, ref = $\circledast$] 
 	 		\item \label{nev}   we have to remark that in steps with $\beta \in S^*$ we will only assume that $\name D_\beta$ is a $\bbP^1_\beta$-name for a system of subsets if $V_\kappa^{\bfV_1}$, where
 	 		\[ \Vdash_{\bbP^1_{\beta}} \ \name D_\beta \in [\cP(V_\kappa^{\bfV_1})]^{\leq \lambda}, \]
 	 		first we will deduce some properties of $\bbP^1_\chi$ based on only this weak assumption up until the end of the proof of Lemmas $\ref{cre}$ and $\ref{cru}$ and then we will verify that the   $\name D_\beta$'s ($\beta \in S^*$) can be suitably chosen (while defining the iteration $\bbP^1_\chi$) so that  the iteration fulfills all our remaining demands from \ref{e1}-\ref{eu}.
 	 		Similarly, for steps in $\chi \setminus S^* \setminus \{0,1\}$ up until the end of the proof of Lemmas $\ref{cre}$ and $\ref{cru}$ we only assume that $\Vdash_{\bbP^1_\beta} \name M_\beta \in (K_{\kappa})_\lambda$, i.e.\ is a $\bbP^1_\beta$-name for a $\kappa$-colored graph on $\lambda$.
 	 	\end{enumerate}
  	\begin{itemize}
 	 		\item For every $M = \langle |M|, R^M_\alpha: \ \alpha < \kappa \rangle \in (K_{\kappa})_\lambda$ we will use the notation $c_M: [\lambda]^2 \to \kappa$ denoting the color of the edge between $i$ and $j$, i.e. 
 	 		\[ c_M(i,j) = \alpha \ \iff \ (i,j) \in R^M_\alpha.\]
 	 	\end{itemize}
 	 
 	 	\begin{enumerate}[label = Case ($1$):, ref = Case ($1$)]
 	 		\item \label{cs1} $\beta = 1$. 	 		
 	 	\end{enumerate}  		Let $\bbQ^1_1 \in \bfV_1^{\bbQ^1_0}$ be the forcing for obtaining a random $\kappa$-colored graph on $\lambda$ with conditions of power $<\kappa$, i.e. $q \in \bbQ^1_1$ iff
  		\begin{enumerate}[label = $(\roman*)$, ref = $(\roman*)$]
  			\item $q \subseteq \{ [i \ R_\gamma \ j], \ [i \ \neg R_\gamma \ j ]: \ i \neq j < \lambda, \ \gamma < \kappa\}$,
  			\item $\forall i \neq j < \lambda$ we have 
  			$$([i \ R_\gamma \ j], [i \ R_{\gamma'} \ j] \in q) \ \longrightarrow \ (\gamma = \gamma'),$$
  			\item $|q| < \kappa$,
  		\end{enumerate}
  	with the usual ordering. Then
  		\newcounter{indfcou} \setcounter{indfcou}{0}
  	\begin{enumerate}[label = $(\diamond)_{\arabic*}$, ref = $(\diamond)_{\arabic*}$]
    		\item the generic object $\name M_* =  \langle\lambda,  \name R^{M_*}_\alpha: \ \alpha < \kappa \rangle$ satisfies
    		\[ \Vdash_{\bbP^1_2} \ \langle \name R^{M_*}_\alpha: \ \alpha < \kappa \rangle \ \text{is a partition of } [\lambda]^2. \]
  		\stepcounter{indfcou}
  	\end{enumerate} 
  		\begin{enumerate}[label = Case ($2$):, ref = Case ($2$)]
 	 		\item \label{case2}  $\beta \in \chi \setminus S^* \setminus \{0,1\}$.
 	 	\end{enumerate}
  		In order to define $\bbQ^1_\beta \in \bfV_1^{\bbP^1_\beta}$ (formally a $\bbP^1_\beta$-name $\name \bbQ^1_\beta \in \bfV_1$) we first need to work in $\bfV_1^{\bbP^1_1}$ as preparation.
  		Let $\Upsilon$ be large enough regular cardinal, and define the continuous increasing chain $\overline{N}_\beta = \langle N_{\beta, \gamma}: \ \gamma < \lambda \rangle \in \bfV_1^{\bbQ^1_0}$ so that
  		\begin{itemize}
  			\item $\beta$,  $\bbP^1_\beta$, $\langle \overline{N}_\gamma: \ \gamma \in \beta \setminus S^* \setminus \{0,1\} \rangle$, $\bfG^1_0 $ $\in N_{\beta, 0}$,
  			\item $\kappa+1  \subseteq N_{\beta,0}$,
  			\item for each $\gamma < \lambda$:
  			\begin{enumerate}[label = $(\bullet)_{\alph*}$, ref = $(\bullet)_{\alph*}$]
  				\item  $N_{\beta, \gamma} \prec \left(\cH(\Upsilon)^{\bfV_1^{\bbP^1_1}}, \in \right)$,
  			\item  $|N_{\beta, \gamma}| < \lambda$,
  			\item $N_{\beta, \gamma} \cap \lambda$ is an initial segment of $\lambda$
  			\item $N_{\beta, \gamma} \cap \lambda < N_{\beta, \gamma+1} \cap \lambda$,
			\item \label{de} for $\varp< \lambda$ limit  $N_{\beta, \varp} = \bigcup_{\gamma < \varp} N_{\beta, \gamma}$,
			\end{enumerate}
  		\end{itemize}
  		and
	  	\begin{enumerate}[label = $(\diamond)_{\arabic*}$, ref = $(\diamond)_{\arabic*}$]
  			\setcounter{enumi}{\value{indfcou}}
  			\item let $\xi_\beta(\gamma) = N_{\beta,\gamma} \cap \lambda$ ($\gamma < \lambda$).
  			\stepcounter{indfcou}
  		\end{enumerate} 
  		So the set $\{ \xi_\beta(\gamma): \ \gamma < \lambda \}$ is a club subset of $\lambda$, and as $S_\beta$ is stationary (Lemma $\ref{Baulem}$) the set $C_\beta = \textrm{cl}(S_\beta \cap \{ \xi_\beta(\gamma): \ \gamma < \lambda \})$ (i.e. the smallest closed set containing $S_\beta \cap \{ \xi_\beta(\gamma): \ \gamma < \lambda \}$) is a club. Therefore the system $\langle N_{\beta,\gamma}: \ \gamma < \lambda \ \wedge \ \xi_\beta(\gamma) \in C_\beta \rangle$ clearly satisfies our requirements, hence (after reparametrization) we can assume that
  		\begin{enumerate}[label = $(\diamond)_{\arabic*}$, ref = $(\diamond)_{\arabic*}$]
  			\setcounter{enumi}{\value{indfcou}}
  			\item \label{tart} $\{\xi_\beta(\gamma+1) : \ \gamma \in \lambda\} \subseteq S_\beta$,
  			\stepcounter{indfcou}
  		\end{enumerate}
  		and we let
  		\begin{enumerate}[label = $(\diamond)_{\arabic*}$, ref = $(\diamond)_{\arabic*}$]
  			\setcounter{enumi}{\value{indfcou}}
  			\item $N_\beta^* = \{\xi_\beta(\gamma) : \ \delta \in \lambda\}$.
  			\stepcounter{indfcou}
  		\end{enumerate} 
  	For later reference we remark that the $\kappa$-almost disjointness of the $S_\alpha$'s and \ref{tart} together implies
  	\begin{enumerate}[label = $(\diamond)_{\arabic*}$, ref = $(\diamond)_{\arabic*}$]
  		\setcounter{enumi}{\value{indfcou}}
  		\item \label{adf} if $\beta \neq \gamma < \chi$ then $|\{\xi_\beta(\delta+1) : \ \delta \in \lambda\} \cap \{\xi_\gamma(\delta+1) : \ \delta \in \lambda\} |  < \kappa$.
  		\stepcounter{indfcou}
  	\end{enumerate} 
  		Now the forcing $\bbQ_\beta^1 \in \bfV_1^{\bbP^1_\beta}$ is defined so that it shall give an embedding $f_\beta$ of the $\kappa$-colored graph $M_\beta \in \bfV_1^{\bbP^1_\beta}$ into $M_*$, formally defined by
  			\begin{enumerate}[label = $(\diamond)_{\arabic*}$, ref = $(\diamond)_{\arabic*}$]
  			\setcounter{enumi}{\value{indfcou}}
  		 		\item $q \in \bbQ_\beta^1$, iff
  		 				\begin{enumerate}[label = $(\roman*)$, ref = $(\roman*)$]
  		 				\item \label{i} $q$ is a set of elementary conditions of the following form
  		 				\begin{itemize}
  		 					\item $[f_\beta(i) = j]$, where $j \in \{ \xi_\beta(\nu+1): \ \kappa i \leq \nu < \kappa(i+1)\}$ (so necessarily $i < j$),
  		 					\item $[ j \notin \ran(f_\beta)]$ for some $j < \lambda$,
  		 				\end{itemize}
  		 				\item \label{ellentmond} the collection $q$ corresponds to a partial injection, and free of any explicitly contradictory subset of terms, under which we mean that 
  		 				\begin{enumerate}[label = $(\alph*)$, ref = $(\alph*)$]
  		 					\item \label{ela}	there are no $i,j \in \lambda$ s.t.\ $[f_\beta(i) = j]$, $[ j \notin \dom(f_\beta)] \in q$,
  		 					\item \label{eb}  there are no $i, j_0 \neq j_1 \in \lambda$ s.t.\ $[f_\beta(i) = j_0]$, $[ f_\beta(i) = j_1] \in q$,
  		 					\item there are no $[f_\beta(i_0) = j_0], [f_\beta(i_1) = j_1] \in q$ s.t.\ $c_{M_\beta}(i_0,i_1) \neq c_{M_*}(j_0,j_1)$.
  		 				\end{enumerate}
  		 				Note that $f_\beta$'s are automatically injective by \ref{i}.
  		 				\item $|q| < \kappa$.
  		 			\end{enumerate}
  		 		
  			\stepcounter{indfcou}
  		\end{enumerate}

 	 	\begin{enumerate}[label = Case ($3$):, ref = Case ($3$)]
 	 		\item \label{case3} $\beta \in S^*$.
 	 	\end{enumerate}
  		As  $\name D_\beta$ is a $\bbP^1_\beta$-name for a system of subsets of $V_\kappa^{\bfV_1}$, if additionally for each $\alpha < \kappa$  $|(\cup \name D_\beta)  \um \alpha| < \kappa$ holds  (and if $\name D_\beta$ generates a proper $\kappa$-complete filter), then we define $\bbQ^1_\beta$ to be the Mathias forcing  $\bbQ_{D_\beta}$ from Definition $\ref{c7}$,  otherwise we can let $\bbQ^1_\beta$ to be the trivial forcing. Note that this requirement ensures that 
  		\begin{enumerate}[label = $(\diamond)_{\arabic*}$, ref = $(\diamond)_{\arabic*}$]
  			\setcounter{enumi}{\value{indfcou}}
  			\item \label{wkorl} if $(w,A) \in \bbQ^1_\beta$, then $|w| < \kappa$.
 		\end{enumerate}	
  		This completes Definition $\ref{itdf}$.
 	 \edd
 	 Now as $\bbP^1_\chi$ is a $<\kappa$-support iteration of $<\kappa$-directed closed posets, $\bbP^1_\chi$ itself is $<\kappa$-directed closed by \cite[Thm 2.7]{B3}, in particular not adding any new sequence of length $<\kappa$,  we have:
 \begin{observation}
 	For each $\beta \in \chi \setminus S^* \setminus \{0,1\}$ forcing with $\bbQ^1_\beta$ over $\bfV_1^{\bbP^1_\beta}$ adds an embedding $f_\beta : M_\beta \to M_*$. 
 \end{observation}
 We already saw that $\bbP^1_1 = \bbQ^1_0$ is $\lambda^+$-cc (Lemma $\ref{Baulem}$), now we prove that in $\bfV_1^{\bbP^1_1}$ the quotient forcing $\bbP^1_\chi / \bfG^1_1$ has the $\kappa^+$-cc (no matter how we choose the $\bbP^1_\beta$-name $\name D_\beta$, or $\name M_\beta$ only satisfying \ref{nev} for $2 \leq \beta <\chi$) after which not only the $\lambda^+$-ccness of $\bbP^1_\chi$ follows, but  some easy calculation will be sufficient for \ref{e2}-\ref{eu}.
 In order to prove the antichain property we will need some technical preparation, the same way as in \cite{Sh:175a}. Recalling that each $\bbP^1_\alpha$ is $<\kappa$-closed (and \ref{wkorl}) is straightforward to prove (by induction on $\alpha$) that
 
 \begin{enumerate}[label = $(*)_{\arabic*}$, ref= $(*)_{\arabic*}$]
 	\setcounter{enumi}{\value{pcounter}}
 	\item The set 
 	$$\begin{array}{rcl} D_\alpha^\bullet = \{ p \in \bbP^1_\alpha: \ \forall \gamma \in \dom(p) & \ \beta \in S^* \to   \exists w_\gamma \in \bfV_1 \text{ s.t.} & \Vdash_{\bbP^1_\gamma} p(\gamma) = ( \check{w}_\gamma, \name A_\gamma ) \\
 		& \text{otherwise:} \ \exists s_\gamma \in \bfV_1 \text{ s.t.} &   \Vdash_{\bbP^1_\gamma} p(\gamma) = \check{s}_\gamma  \} \end{array} $$
 	is a dense subset of $\bbP^1_\alpha$.
 	\stepcounter{pcounter}
 	\item Therefore, in the quotient forcing $\bbP^1_\alpha/ \bfG^1_1$ (as defined in \cite{B3}, or see below) the set
 	\[ D_\alpha^0 = \{ p \in \bbP^1_\alpha/ \bfG^1_1: \ \exists  q_0 \in \bfG^1_1: \ \langle q_0 \rangle \cup p \in D_\alpha^\bullet   \} \ \in \bfV_1^{\bbP^1_1} \]
 	is dense (where $\bbP^1_\alpha/ \bfG^1_1 = \{ p\um (\dom(p) \setminus \{0\}): \ p \in \bbP^1_\alpha\} \in \bfV_1^{\bbP^1_1}$, and $p \leq_{\bbP^1_\alpha/ \bfG^1_1} q$,  iff for some $r_0 \in \bfG_1^1 \subseteq \bbP^1_1$ $\langle r_0 \rangle \cup p \leq_{\bbP^1_\alpha} \langle r_0 \rangle \cup q$).
 	\stepcounter{pcounter}
 	\item With a slight abuse of notation (in order to avoid further notational awkwardness) we will identify each condition $p \in D^0_\alpha \subseteq \bbP^1_\alpha/ \bfG^1_1$ with the function on the same domain, but  for each $\gamma \in \dom(p)$
 	\begin{itemize} 
 		\item  if $\beta \in S^*$ then writing $p(\beta) = (w, \name A )$ (instead of 
 	 some $\bbP^1_\beta$-name satisfying $\langle q_0 \rangle \cup p\upharpoonright \gamma \Vdash_{\bbP^1_\beta} p(\beta) = ( \check{w}, \name A)$ for some $q_0 \in \bfG_1^1$), 
 	 \item or $p(\beta)= s$, where $s$ is a set of symbols as in Case (1), (2) in Definition $\ref{itdf}$ (instead of $\langle q_0 \rangle \cup p\upharpoonright \beta \Vdash_{\bbP^1_\beta} p(\beta) = \check{s}$ for some $q_0 \in \bfG^1_1$).
 	 \end{itemize}
  	\stepcounter{pcounter}
 \end{enumerate}
  Note as $\bbP^1_1$ is $<\kappa$-closed (recall that $D_\alpha^0 \subseteq \bfV_1$) that
  \begin{enumerate}[label = $(*)_{\arabic*}$, ref= $(*)_{\arabic*}$]
  	\setcounter{enumi}{\value{pcounter}}
  	\item (in $\bfV_1^{\bbP^1_1}$) for any $\alpha \leq \chi$, and increasing sequence $\overline{p} = \langle p_\zeta: \ \zeta < \varp < \kappa \rangle$  in $D^0_\alpha$  has a least upper bound in $\bbP^1_\alpha$, which we will denote by $\lim_{\zeta < \varp} p_\zeta$, and this limit is in $D^0_\alpha$. For the sake of completeness we include the formal definition of $ \lim_{\zeta < \varp} p_\zeta$. The limit of $\overline{p} = \langle p_\zeta: \ \zeta < \varp < \kappa \rangle$ is the function $p^*$, where
  	\begin{enumerate}[label = $(\alph*)$, ref =  $(\alph*)$]
  \item   $\dom(p^*) =  \bigcup_{\zeta<\varp} \dom(p_\zeta)$,
  	\item for $\beta \in S^* \cap \dom(p^*)$ $p^*(\beta) = (\bigcup_{\zeta < \varp} w_{p_\zeta(\beta)}, \name A_\beta)$, where $p_\zeta(\beta) = (w_{p_\zeta(\beta)}, A_{p_\zeta(\beta)})$, and $\name A_\beta$ is the $\bbP^1_\beta$-name defined so that $\Vdash_{\bbP^1_\beta} \name A_\beta = \bigcap_{\zeta < \varp} \name A_{p_\zeta(\beta)}$ holds,
	 \item for $\beta \in \chi \setminus S^* \setminus \{0,1\}$ set $p^*(\beta) = \bigcup_{\zeta < \varp} p_\zeta(\beta)$.
	\end{enumerate}
	\stepcounter{pcounter} 
 \end{enumerate}

\bdd \label{levag} In $\bfV_1^{\bbP^1_1}$ for $\alpha \leq \chi$, $\delta \leq \lambda$ for each condition $p \in D^0_\alpha$ we define $p^{[\delta]}$ to be the function with $\dom(p^{[\delta]}) = \dom(p)$,
\begin{enumerate}[label = $(\alph*)$, ref= $(\alph*)$]
\item if $1 \in \dom(p)$, then $p^{[\delta]}(1) = \{ [i \ R_\gamma \ j] \in p(1): \ i,j < \delta \}$,
\item for $1 < \beta \in \dom(p) \cap S^*$ we let $p^{[\delta]}(\beta) = p(\beta)$,
\item  \label{levd} otherwise (for $1< \beta \in \dom(p) \setminus S^*$) we let
$$\begin{array}{rc} p^{[\delta]}(\beta) = & \{ [f_\beta(i)= j] \in p(\beta): \ i,j <\max\{ \xi_\beta(\gamma): \ \gamma<\lambda, \ \xi_\beta(\gamma) \leq \delta   \} \} \\ 
	& \cup \\
	& \{ [j \notin \ran(f_\beta)] \in p(\beta): \ j  <\max\{ \xi_\beta(\gamma): \ \gamma<\lambda, \ \xi_\beta(\gamma) \leq \delta\} \}.
 \end{array}$$
\end{enumerate}
\edd 
Observe that, because of each $p$ and each $p(\beta)$ ($\beta \in \dom(p)$) has support of size $<\kappa$, and $\lambda > \kappa$ is regular,
\begin{enumerate}[label = $(*)_{\arabic*}$, ref= $(*)_{\arabic*}$]
	\setcounter{enumi}{\value{pcounter}}
	\item for each $\alpha \leq \chi$, $p \in D^0_\alpha \subseteq (\bbP^1_\alpha/ \bfG^1_1)$ we have   $p^{[\delta]} = p$ for every large enough $\delta$, and
	\item clearly $p^{[\delta]} \um \beta = (p \um \beta)^{[\delta]}$ (for $\beta < \alpha$).
	\item \label{obha}  for $p \leq q \in D^0_\alpha$ with $p^{[\delta]}, q^{[\delta]} \in D^0_\alpha$ we obviously have $p^{[\delta]} \leq q^{[\delta]}$.
		\stepcounter{pcounter} \stepcounter{pcounter} \stepcounter{pcounter} 	
\end{enumerate}
Note that for $p \in D^0_\alpha \subseteq \bbP^1_\alpha/ \bfG^1_1$ the reduced function $p^{[\delta]}$ is in $\bfV_1^{\bbP^1_1}$ (even in $\bfV_1$), but is not necessarily a condition in $\bbP^1_\alpha/ \bfG^1_1$.
It is straightforward to check by induction on $\alpha$ the following.
\begin{observation} \label{obsz}
For each $\alpha \leq \chi$, $p \in D^0_\alpha$ and $\delta < \lambda$ 
\begin{enumerate}[label = $\alph*)$, ref= $\alph*)$]
	\item 
	$p^{[\delta]}$ is an actual condition (i.e. belongs to $D^0_\alpha \subseteq \bbP^1_\alpha/\bfG^1_1$), iff for every $\beta \in \dom(p)$ 
	$$p^{[\delta]} \upharpoonright \beta \in \bbP^1_\beta,$$
	and (letting $\delta^-_\beta = \max (N^*_\beta \cap (\delta+1))$
	\beeq \label{eeg} \begin{array}{l} \forall [f_\beta(i_0) = j_0], [f_\beta(i_1) = j_1] \in p(\beta): \\
	 j_0,j_1 < \delta_\beta^- \ \longrightarrow \ p^{[\delta]} \upharpoonright \beta \Vdash_{\bbP^1_\beta / \bfG^1_1} c_{M_\beta}(i_0,i_1) = c_{M_*}(j_0,j_1). \end{array} \eeq
	\item In particular, for limit $\alpha$
	\[ p^{[\delta]} \in \bbP^1_\alpha / \bfG^1_1 \ \iff \ (\text{for cofinally many } \varp < \alpha:) \ p^{[\delta]} \um \varp \in \bbP^1_\varp, \]
	\item while for $\alpha = \beta+1$
	\[  p^{[\delta]} \in \bbP^1_\alpha/ \bfG^1_1 \ \iff \ p^{[\delta]} \um \beta \in \bbP^1_\beta/ \bfG^1_1 \text{ and }\eqref{eeg} \text{ holds.} \]
\end{enumerate}
\end{observation}
The following notion and lemma is of central importance.
\bdd
In $\bfV_1^{\bbP^1_1}$ for $\alpha \leq \chi$ define
\[ D^*_\alpha = \{ p \in D^0_\alpha: \forall \delta < \lambda \ p^{[\delta]} \in \bbP^1_\alpha / \bfG^1_1 \}.\]
\edd
Having Observation $\ref{obsz}$ in our mind it is easy to check the following.
\begin{enumerate}[label = $(*)_{\arabic*}$, ref= $(*)_{\arabic*}$]
	\setcounter{enumi}{\value{pcounter}}
	\item \label{*lim} Whenever $\langle p_\zeta: \ \zeta < \varp < \kappa \rangle$ is an increasing sequence in $D^*_\alpha$, then $\lim_{\zeta < \varp} p_\zeta \in D^*_\alpha$.
	\stepcounter{pcounter}
\end{enumerate}
This leads us to note how  the statements $p \in D^*_\alpha$ and $p\um \beta \in D^*_\beta$ ($\beta < \alpha$) relate to each other.
\begin{observation} \label{obsz2}
	For each $\alpha \leq \chi$, $p \in D^0_\alpha$
	\begin{enumerate}[label = $\alph*)$, ref= $\alph*)$]
		\item 
		$p\in D^*_\alpha$, iff for every $\beta \in \dom(p)$ 	and for every $\delta < \lambda$
		$$p \upharpoonright \beta \in D^*_\beta,$$
	 and (letting $\delta^-_\beta = \max (N^*_\beta \cap (\delta+1))$
		\beeq \label{eeg2} \begin{array}{l} \forall [f_\beta(i_0) = j_0], [f_\beta(i_1) = j_1] \in p(\beta): \\
			j_0,j_1 < \delta_\beta^- \ \longrightarrow \ p^{[\delta]} \upharpoonright \beta \Vdash_{\bbP^1_\beta / \bfG^1_1} c_{M_\beta}(i_0,i_1) = c_{M_*}(j_0,j_1). \end{array} \eeq
		\item \label{bsz2b} In particular, for limit $\alpha$
		\[ p \in D^*_\alpha \ \iff \ (\text{for cofinally many } \varp < \alpha): \ p \um \varp \in D^*_\varp, \]
		\item  \label{bsz2c} while for $\alpha = \beta+1$
		\[  p \in D^*_\alpha \ \iff \ [p \um \beta \in D^*_\beta] \text{ and }[ \text{for each } \delta < \lambda \ \eqref{eeg2} \text{ holds}.]  \]
	\end{enumerate}
\end{observation}

\bl \label{cru} For $\alpha \leq \chi$
\begin{enumerate}[label = $(\blacksquare)^1_\alpha$, ref = $(\blacksquare)^1_\alpha$ ]
	\item \label{cruu} $$\bfV_1^{\bbP^1_1} \models D^*_\alpha \text{ is dense in } \bbP^1_\alpha / \bfG^1_1.$$
\end{enumerate}
\el
\bl \label{cre} 
 For every $\alpha \leq \chi$
\begin{enumerate}[label = $(\blacksquare)^2_\alpha$, ref = $(\blacksquare)^2_\alpha$ ]
	\item \label{cree} $$\bfV_1^{\bbP^1_1} \models \ \bbP^1_\alpha / \bfG^1_1 \text{ has the }\kappa^+\text{-cc}.$$
\end{enumerate} 
\el
\begin{PROOF}{Lemmas \ref{cru} and \ref{cre}}
We proceed by induction, and prove Lemmas $\ref{cru}$ and $\ref{cre}$ simultaneously: More exactly we prove Lemma $\ref{cru}$ for $\alpha$ provided that both Lemmas holds for $\beta$'s less than $\alpha$, and we verify the $\kappa^+$-cc property for $\bbP^1_\alpha$ assuming that $D^*_\alpha$ is a dense subset of $\bbP^1_\alpha / \bfG^1_1$.   
For $\alpha \leq 2$ (when $\bbP^1_2 / \bfG^1_1$ is essentially the forcing $\bbQ^1_1$ of the random graph \ref{cs1} of Definition $\ref{itdf}$) the statement \ref{cruu} clearly holds.

Suppose that that we know that for each $\varp < \alpha$ \ref{cruu} and \ref{cree}  hold . Assume first that $\alpha$ is limit. If $\cf(\alpha) \geq \kappa$, then $\bbP^1_\alpha = \bigcup_{\varp < \alpha} \bbP^1_\varp$, $D^*_\alpha = \bigcup_{\varp < \alpha} D^*_\varp$, so the latter is dense, we are done.

Second, if $\alpha$ is limit, but $\cf(\alpha) < \kappa$, then let $\langle \eta_\theta: \ \theta < \cf(\alpha) \rangle$ be a continuous increasing sequence with limit $\alpha$, let $p_{-1} \in D^0_\alpha$ be arbitrary. We will choose the increasing sequence $\langle p_\theta: \ \theta < \cf(\alpha) \rangle$ in $D^0_\alpha$ with $p_0 \geq p_{-1}$, and $p_\theta \um \eta_\theta \in D^*_{\eta_\theta}$. This would suffice as for each $\theta < \cf(\kappa)$ the sequence $p_\varrho \um \eta_\theta$ ($\varrho \in \cf(\alpha)$) is eventually in $D^*_{\eta_\theta}$, so for $p^* = \lim_{\varrho < \cf(\alpha)} p_\varrho$ using \ref{*lim} we have $p^* \um \eta_\theta \in D^*_{\eta_\theta}$, leading to
$$ (\forall \theta < \cf(\alpha)) \ p^* \um \eta_\theta \in D^*_{\eta_\theta},$$
so by \ref{bsz2b} we are done. For the construction of the $p_\theta$'s,  as $D^0_\alpha$ and $D^*_{\eta_\theta}$'s are $<\kappa$-closed we only have to ensure that $p_\theta \in D^0_\alpha$ can be chosen so that not only $p_\theta \geq p_\varrho$ ($\varrho < \theta$), but 
$p_\theta \um \eta_\theta \in D^*_{\eta_\theta}$. Now applying the induction hypothesis, and extending $(\lim_{\varrho < \theta}p_\varrho) \um \eta_\theta \leq p^*_\theta \in D^*_{\eta_\theta}$ we can choose $p_\theta$ to be the least upper bound of $p^*_\theta$ and $(\lim_{\varrho < \theta}p_\varrho)$ (in fact for $\theta$ limit we did not even have to appeal to the induction hypothesis). 

Third, if $\alpha = \beta + 1$, let $p_{-1} \in D^0_\alpha$, and we will extend $p_{-1} \um \beta \leq p^* \in D^*_\beta$ (using $(\blacksquare)^1_\beta$) so that the right hand side of Observation $\ref{obsz2}$ \ref{bsz2c} holds for $p = p^* \cup \langle p_{-1}(\beta) \rangle$.

For this, let $\{j_\theta: \ \theta < \nu \}$ enumerate $\{j < \lambda: \ [f_\beta(i) = j] \in p_{-1}(\beta) \text{ for some } i <\lambda\}$ in increasing order, and we can fix the system $\{ i_\theta: \ \theta < \nu\}$ so that 
\newcounter{pnter} \setcounter{pnter}{0}
\begin{enumerate}[label = $(\odot)_{\arabic*}$, ref= $(\odot)_{\arabic*}$]  
	\item \label{1ee}  $\{ i_\theta: \ \theta < \nu \}$ is such that
	for each $\theta$ $[f_\beta(i_\theta) = j_\theta] \in p_{-1}(\beta)$.
	
	\stepcounter{pnter}
\end{enumerate}
Note that by Definition $\ref{itdf}$/\ref{case2}/ \ref{i}
\begin{enumerate}[label = $(\odot)_{\arabic*}$, ref= $(\odot)_{\arabic*}$]
	\setcounter{enumi}{\value{pnter}}
	\item \label{kisebbt}  for each $\theta$: $i_\theta < j_\theta$,
	\stepcounter{pnter}
\end{enumerate}
 and also we can choose $\gamma_\theta$ for each $\theta < \nu$ such that
$\xi_\beta(\gamma_\theta) = j_\theta$, thus
\begin{enumerate}[label = $(\odot)_{\arabic*}$, ref= $(\odot)_{\arabic*}$]
	\setcounter{enumi}{\value{pnter}}
	\item we have \[ \{j < \lambda: \ \exists i < \lambda \ [f_\beta(i) = j] \in p_{-1}(\beta)\} = \{j_\theta: \ \theta < \nu \} = \{\xi_\beta(\gamma_\theta): \ \theta < \nu \}.\]
	\stepcounter{pnter}
\end{enumerate}

Now we construct the increasing sequence $\langle p_\theta : \ \theta < \nu \rangle$ in $D^*_\beta$ with the properties
\begin{enumerate}[label = (\greek*), ref = (\greek*)]
	\item \label{ee} $p_{-1} \um \beta \leq p_0$,
	\item for each $\theta < \nu$, for each $\varp_0 < \varp_1 < \theta$
	\[ p_\theta^{[\xi_\beta(\gamma_{\varp_1}+1)]} \Vdash_{\bbP^1_\beta / \bfG^1_1} c_{M_\beta}(i_{\varp_0}, i_{\varp_1}) = c_{M_*}(j_{\varp_0}, j_{\varp_1}).\]
	\stepcounter{pnter}\stepcounter{pnter}
\end{enumerate}
This clearly suffices, as we can let $p^* = \lim_{\theta < \nu} p_\theta \in D^*_\beta$, and then $p^* \cup \langle p_{-1}(\beta) \rangle$ belongs to $D^*_\alpha$:  as $j_{\varp_1} = \xi_\beta(\gamma_{\varp_1})$ so $\xi_\beta(\gamma_{\varp_1}+1)$ is the minimal $\delta < \lambda$ s.t. 
$(p^* \cup \langle p_{-1}(\beta) \rangle)^{[\xi_\beta(\gamma_{\varp_1}+1)]}(\beta)$ contains the symbol $[f_\beta(i_{\varp_1}) = j_{\varp_1}]$, therefore by Observation $\ref{obsz2}$ \ref{bsz2c} we are done, \ref{cruu} follows, indeed.

Appealing to the induction hypothesis let $p_0 \in D^*_\beta$, $p_0 \geq p_{-1}$.
Using the $<\kappa$-closedness of $D^*_\beta$ (\ref{*lim}) it is enough to deal with the successor case, that is, for each $\theta$ choose $p_{\theta+1}$ so that $p_{\theta+1}^{[\xi_\beta(\gamma_\theta+1)]}$ forces that the partial function $i_\varp \mapsto j_\varp$ ($\varp \leq \theta$) is an embedding of $\name M_\beta \um \{i_\varp: \ \varp \leq \theta\}$ into $\name M_* \um \{j_\varp: \ \varp \leq \theta\}$. 
Using again \ref{*lim} 
\begin{enumerate}[label = $(\odot)_{\arabic*}$, ref= $(\odot)_{\arabic*}$]
	\setcounter{enumi}{\value{pnter}}
	\item \label{kcl} it suffices to show that for each $\varp < \theta$ and $q \geq p_{-1} \um \beta$, $q \in D^*_\beta$ there exists $q' \in D^*_\beta$, $q' \geq q$
	\[ q'^{[\xi_\beta(\gamma_{\theta}+1)]} \Vdash_{\bbP^1_\beta / \bfG^1_1} c_{M_\beta}(i_{\varp}, i_{\theta}) = c_{M_*}(j_{\varp}, j_{\theta}).\]
	\stepcounter{pnter}
\end{enumerate}
We will see that this follows from the following (formally) more general lemma, stated here for later reference.
\bl \label{tech} For every $\beta \leq \chi$, $q \in D^*_\beta$, $\delta < \lambda$,  $i'$, $i'' \in \max(N^*_\beta \cap (\delta+1))$  there exists $q' \in D^*_\beta$,  $q' \geq q$ such that
\[ q'^{[\delta]} \text{ forces a value for  } c_{M_\beta}(i',i''). \]
Moreover, if  $(\forall \gamma \in \dom(q) \setminus S^*) \ ([f_\beta(i) = j)] \in  q(\gamma) \setminus q^{[\delta]}(\gamma)) \longrightarrow (j = \max(N^*_\gamma \cap (\delta+1)) \ \wedge j < \delta)$ (hence $\delta \notin N^*_\gamma$ and $q(1) = q^{[\delta]}(1)$), then there exists  $q'$ with
$$ (\forall \gamma \in \dom(q') \setminus S^*): \ \ q'(\gamma) \setminus q'^{[\delta]}(\gamma) = q(\gamma) \setminus q^{[\delta]}(\gamma) \ . $$    

\el
(Here we remark that the proof of the lemma uses the $\kappa^+$-cc property of $\bbP^1_\beta / \bfG^1_1$, but we will only use it for proving \ref{kcl}, that is  to complete the proof of $((\blacksquare)^1_\beta \wedge (\blacksquare)^1_\beta) \to $\ref{cruu}.)
\begin{PROOF}{Lemma \ref{tech}}
	
	So fix $q \in D^*_\beta$, let $\varrho$ be chosen so that $\xi_\beta(\varrho) = \max(N^*_\beta \cap (\delta+1))$, so $i',i'' < \xi_\beta(\varrho) \leq \delta$, and recall that for the model $N_{\beta, \varrho} \prec (\cH^{\bfV_1^{\bbP^1_1}}(\Upsilon), \in)$ we know that $i',i'',\name M_\beta, \bbP^1_\beta, \bfG^1_1 \in N_{\beta,\varrho}$ (and thus $\bbP^1_\beta/ \bfG^1_1 \in N_{\beta,\varrho}$). So we can find $A \in N_{\beta,\varrho}$ such that $A$ is a maximal antichain in $N^0_\beta \subseteq \bbP^1_\beta/ \bfG^1_1$, each $p \in A$ decides the value of $c_{M_\beta}(i', i'')$. But as $\bbP^1_\beta / \bfG^1_1$ has the $\kappa^+$-cc, and $\kappa+1 \subseteq N_{\beta,\varrho}$ we have that $A \subseteq N_{\beta,\varrho}$.
	
	 So 
	 \newcounter{lemtr} \setcounter{lemtr}{0}
	\begin{enumerate}[label = $(\boxplus)_{\arabic*}$, ref= $(\boxplus)_{\arabic*}$]
		\setcounter{enumi}{\value{lemtr}}
		\item let $q' \in D^*_\beta$ be a common upper bound of $q$ and some $q'' \in A$.
		\stepcounter{lemtr}
	\end{enumerate}
	
	We have to argue that not only $ q' \Vdash_{\bbP^1_\beta / \bfG^1_1} c_{M_\beta}(i', i'') = c_*$ (for some $c_* < \kappa$) 
	but 
	\beeq \label{vf} q'^{[\delta]} \Vdash_{\bbP^1_\beta / \bfG^1_1} c_{M_\beta}(i', i'') = c_*. \eeq For $\eqref{vf}$  it is enough to prove that $q''^{[\delta]} = q''$, because then 
	$q'^{[\delta]} \geq q''^{[\delta]} = q''$ (by \ref{obha}), yielding $\eqref{vf}$, as we wanted.
	But as $q'' \in N_{\beta, \varrho}$, and $\lambda \cap N_{\beta, \varrho} = \xi_\beta(\varrho) \leq \delta$ for each $\zeta \in \dom(q'')\setminus S^* \setminus \{0,1\}$ we have $\langle N_{\zeta, \iota}: \iota < \lambda \rangle \in N_{\beta, \varrho}$ (recall \ref{case2} from Definition $\ref{itdf}$), so $\xi_\beta(\varrho)$ is an accumulation point of the $\xi_\zeta(\iota)$'s. Hence we get that

\begin{enumerate}[label = $(\boxplus)_{\arabic*}$, ref= $(\boxplus)_{\arabic*}$]
	\setcounter{enumi}{\value{lemtr}}
		\item \label{kmeeee} for each $\zeta \in \dom(q'')\setminus S^* \setminus \{0,1\}$ 
		$\xi_\beta(\varrho) = \xi_\zeta(\iota)$ for some $\iota < \lambda$ (in fact, for $\iota = \xi_\beta(\varrho)$),
		\stepcounter{lemtr}
	\end{enumerate}
	so $q''^{[\xi_\beta(\varrho)]} = q''^{[\delta]} = q''$, we are done.
	
	Finally, for the moreover part 
	let $\delta_\gamma^- = \max(N_\gamma \cap (\delta+1))$), and define $i_\gamma^-$ to be the unique ordinal s.t.\  
	\beeq \label{idelta} [f_\gamma(i_\gamma^-) = \delta_\gamma^-] \in q(\gamma) \eeq (if there exists). Note that our conditions on $q$ imply that if $i_\gamma^-$ is defined, then $i_\gamma^- < \delta_\gamma^-$.
	Now by induction and by the first part
	 define $q'' \geq q$ such that for every $\gamma \in \dom(q'') \setminus S^*$ with $i_\gamma^-$ defined 
	 $$ ([f_\gamma(i) = j] \in q''^{[\delta]}(\gamma)) \ \to \  q''^{[\delta]} \um \gamma \text{ decides the value } c_{\name M_\gamma}(i,i_\gamma^-),$$
	 and
	 $$ ([f_\gamma(i) = j] \in q''^{[\delta]}(\gamma)) \ \to \  q''^{[\delta]} (1) \text{ decides the value } c_{\name M_*}(j,\delta_\gamma^-)$$
	 (in fact this latter follows from $j,\delta_\gamma^-< \delta$ and $\eqref{idelta}$).
	Now clearly $q''^{[\delta]} \geq q^{[\delta]}$, and we want to define the condition $q'$ to be the least upper bound of $q''^{[\delta]}$ and $q$, which is possible, as for every $\gamma$ with $i^-_\gamma$ defined we have that $ q''^{[\delta]} \um \gamma$ forces that $q''^{[\delta]}(\gamma) \cup \{[f_\gamma(i^-_\gamma) = \delta^-_\gamma]\}$ is indeed a partial embedding.
	
\end{PROOF}
 Turning back to the statement from \ref{kcl}, as $j_\varp < j_\theta = \xi_\beta(\gamma_\theta) < \xi_\beta(\gamma_\theta+1)$ we also have $i_\varp,i_\theta < \xi_\beta(\gamma_\theta)$ (thus obviously $i_\varp,i_\theta < \xi_\beta(\gamma_\theta+1)$). 
 Apply the lemma with $\delta = \xi_\beta(\gamma_\theta+1)$, $i' = i_\varp$, $i'' = i_\theta$,
\begin{enumerate}[label = $(\odot)_{\arabic*}$, ref= $(\odot)_{\arabic*}$]
	\setcounter{enumi}{\value{pnter}}
	\item let $q' \in D^*_\beta$ be given by the lemma, so $ q' \Vdash_{\bbP^1_\beta / \bfG^1_1} c_{M_\beta}(i_{\varp}, i_{\theta}) = c_{M_*}(j_{\varp}, j_{\theta})$
		\stepcounter{pnter}
\end{enumerate}
 (which is obvious, as
	\begin{enumerate}[label = $(\odot)_{\arabic*}$, ref= $(\odot)_{\arabic*}$]
		\setcounter{enumi}{\value{pnter}}
		\item \label{kme} $q' \geq p_{-1} \um \beta$, and $p_{-1}$ is a proper condition in $D^0_\alpha$ with $[f_\beta(i_\theta) = j_\theta]$, $[f_\beta(i_\varp) = j_\varp] \in p_{-1}(\beta)$, hence $q' \tieconcat \langle p_{-1}(\beta)\rangle$, too)
		\stepcounter{pnter}
	\end{enumerate}
we have to argue that
\beeq \label{vfo} q'^{[\xi_\beta(\gamma_{\theta}+1)]} \Vdash_{\bbP^1_\beta / \bfG^1_1} c_{M_\beta}(i_{\varp}, i_{\theta}) = c_{M_*}(j_{\varp}, j_{\theta}). \eeq 
But $q'^{[\xi_\beta(\gamma_{\theta}+1)]} \Vdash_{\bbP^1_\beta / \bfG^1_1} c_{M_\beta}(i_\varp, i_\theta) = c_*$, and if $[j_\varp \ R_{c_*} \ j_\theta] \notin q'^{[\xi_\beta(\gamma_{\theta}+1)]}(1)$, so not in $q'(1)$, then adding $[j_\varp \ R_{c_*+1} \ j_\theta]$ to the first coordinate of $q' \tieconcat \langle p_{-1}(\beta)\rangle \in D^0_\alpha$ would lead to a contradiction. This verifies that assuming the induction hypotheses for $\beta$ the set  $D^*_{\beta+1}$ is dense in $\bbP^1_\beta / \bfG^1_1$.

Now assuming that $D^*_\alpha$ is dense we are ready to prove that $\bbP^1_\alpha / \bfG^1_1$ has the $\kappa^+$-cc.
So let $\langle p_\gamma: \ \gamma < \kappa^+ \rangle$ be an antichain in $D^*_\alpha$.
By  extending each $p_\gamma$
\begin{enumerate}[label = $(\odot)_{\arabic*}$, ref= $(\odot)_{\arabic*}$]
	\setcounter{enumi}{\value{pnter}}
		\item \label{form} we can assume that for each $\gamma < \kappa^+$
			\begin{enumerate}[label = $(\roman*)$, ref= $(\roman*)$]
				 \item \label{fo0}for each $\beta \in \dom(p_\gamma)$, for each $i_0,i_1,j_0<j_1$ with
					$[f_\beta(i_0) = j_0], [f_\beta(i_1) = j_1] \in p_\gamma(\beta)$
					the condition $p^{[j_1]} \um \beta$ decides the value $c_{\name M_\beta}(i_0,i_1)$,
				\item \label{fo1} for each $\gamma < \kappa^+$ the condition $p_\gamma(1)$ is a complete graph on the set $L_\gamma$ with its edges colored, i.e.\
				$$\begin{array}{l} L_\gamma = \{ i < \lambda: \ \exists i' < \lambda \ \exists \varp < \kappa \ [i \ R_\varp \ i'] \in p_\gamma(1) \},  \\
				\text{so }	(\forall  i,j \in L_\gamma) \ (\exists \delta < \kappa): \  [i \ R_\delta \ j] \in p_\gamma(1).
				\end{array}  $$
			
				\item \label{fo2} for each $\gamma < \kappa^+$ and $\beta \neq \beta'  \in \dom(p_\gamma) \setminus S^* \setminus \{0,1\}$ we have $\{\xi_\beta(\rho+1): \ \rho<\lambda \} \cap \{\xi_{\beta'}(\rho+1): \ \rho<\lambda \} \subseteq L_\gamma$
			 (recall that $|\{\xi_\beta(\rho+1): \ \rho<\lambda \} \cap \{\xi_{\beta'}(\rho+1): \ \rho<\lambda \}| < \kappa$ by \ref{adf}),
						
				\item \label{ggr} for each $\gamma < \kappa^+$ and $\beta  \in \dom(p_\gamma) \setminus S^* \setminus \{0,1\}$ and for each $j < \lambda$ if $[j \notin \ran(f_\beta)] \in p_\gamma(\beta)$, or $[f_\beta(i) = j] \in p_\gamma(\beta)$ for some $i < \lambda$, then $j \in L_\gamma$,

				\item \label{fo3}  for each $\gamma < \kappa^+$, $\beta \in \dom(p_\gamma)$ and $j < \lambda$, if $j \in L_\gamma$, then  either $[j \notin \ran(f_\beta)] \in p_\gamma(\beta)$, or (for some $i$) $[f_\beta(i) = j] \in p_\gamma(\beta)$,
				
				\item \label{fo+1} the set $L_\gamma \subseteq \lambda$ is closed, of limit order type,

		\end{enumerate}
	\stepcounter{pnter}
\end{enumerate}
[This is possible, a simple induction using Lemma $\ref{tech}$, the fact $$[f_\beta(i)= j] \in p_\gamma(\beta) \ \to \ j \in N^*_\beta$$  and \ref{*lim} yields that there is $p_\gamma' \geq p_\gamma$ in $D^*_\alpha$, with $(p_\gamma' \um \beta)^{[j_1]}$ determining the value $c_{\name M_\beta}(i_0,i_1)$ whenever  $[f_\beta(i_0) = j_0], [f_\beta(i_1) = j_1] \in p_\gamma(\beta)$, $j_0<j_1$. Now repeating this $\omega$-many times we get a condition satisfying \ref{fo0}. Then we can obtain an even stronger condition satisfying \ref{fo1}-\ref{fo+1} by only adding symbols of the form $[j \notin \ran(f_\beta)]$ at coordinates $1< \beta \in \chi \setminus S^*$ and extending also $p_\gamma'(1)$.]
As $\kappa$ is inaccessible in $\bfV_1$ by \ref{apr}, and in $\bfV_1^{\bbP^1_1}$ as $\bbP^1_1$ is $<\kappa$-closed 
we can apply the delta system lemma, so w.l.o.g.\ $\langle \dom(p_\gamma): \ \gamma < \kappa^+ \rangle$ forms a delta system. 
 By applying the delta system lemma again we can assume that for each $\beta \in \chi \setminus S^*$ each of the following systems of sets forms a delta system:
 \begin{itemize}
 	\item$ \begin{array} {llr} L_\gamma  \ (\gamma < \kappa^+) \end{array}$,
 	\item $ I_\gamma(\beta) =    \left\{ \begin{array} {ll}  i: & \ [f_\beta(i) = j]\in p_\gamma(\beta) \vee  \exists j \in [\xi_\beta(\kappa i), \xi_\beta(\kappa (i+1)))  \\
 	& [j \notin \ran(f_\beta)] \in p_\gamma(\beta) \end{array} \right\}  (\gamma < \kappa^+).   $
 \end{itemize}
 
Therefore (recalling that each $i < \lambda$ has $\kappa$-many possible images) there are $\xi \neq \zeta < \kappa^+$, such that $p_\xi$ and $p_\zeta$ has no explicitly contradictory terms on the coordinates concerning the $\kappa$-colored graphs, and agreeing in the first part of the condition on the coordinates dedicated to Mathias forcing, under which we  mean the following (w.l.o.g.\ we can assume that $\xi =0$, $\zeta =1$):
\begin{enumerate}[label = $(\odot)_{\arabic*}$, ref= $(\odot)_{\arabic*}$]
	\setcounter{enumi}{\value{pnter}}
	\item for $\beta = 1$ for each $i,j \in L_0(1) \cap L_1(1)$ there exists some $\varp < \kappa$ s.t.\ $[i \ R_\varp \ j] \in p_0(1) \cap p_1(1)$,
	\item \label{embcomp} for $\beta \in \chi  \setminus S^* \setminus \{0,1 \}$ (if $\beta \in \dom(p_0) \cap \dom(p_1)$) the set $p_0(\beta) \cup p_1(\beta)$ determines a partial injection from a subset of $\lambda$ to a subset of $\lambda$, i.e. satisfies  \ref{ellentmond} \ref{ela}, \ref{eb} (from Definition $\ref{itdf}$ \ref{case2}),
	\item for $\beta \in S^* \cap \dom(p_0) \cap \dom(p_1)$  $p_0(\beta) = (w_\beta, \name A_{0,\beta})$, $p_1(\beta) = (w_\beta, \name A_{1,\beta})$ for some $w_\beta \in [V_\kappa^{\bfV_1}]^{<\kappa}$, and $\bbP^1_\beta$-names $\name A_{0.\beta}$, $\name A_{1,\beta}$.
	
	\stepcounter{pnter}
\end{enumerate}
Now $p_0$ and $p_1$ seem good candidates for a compatible pair in our supposed antichain, but we cannot take just the upper bound coordinatewise, as for coordinates $\beta >1$ outside $S^*$ it will not necessarily force that $p_0(\beta) \cup p_1(\beta)$ is an embedding of $\name M_\beta$ to $\name M_*$.
Although it is not immediate, the following claim shows that we can construct a common upper bound, completing the proof of \ref{cree} for $\alpha$.
\bcl \label{kite} There exists a condition $q \in D^*_\alpha$ extending both $p_0$ and $p_1$.
\ecl
\begin{PROOF}{Claim \ref{kite}}
	Let $\{ j_\varp: \ \varp < \varrho\}$ enumerate $L_0 \cup L_1 = \{j: \ [j \ R_\nu \ j'] \in p_0(1) \cup p_1(1) \ \text{ for some } j'<\lambda, \ \nu < \kappa  \}$ (in increasing order).
	\newcounter{bcou} \setcounter{bcou}{0}
	\begin{enumerate}[label = $(\bullet)_{\arabic*}$, ref= $(\bullet)_{\arabic*}$]
		\item \label{clo1} As $L_0$, $L_1$ are closed sets of ordinals without maximal element (\ref{fo+1}) obviously so is $\{ j_\varp: \ \varp < \varrho\}$, let $j_\varrho$ be its supremum.
		\stepcounter{bcou}
		\item \label{2..} By adding symbols of the form $[j \notin \ran(f_\beta)]$ to $p_0(\beta)$, $p_1(\beta)$ we can assume the following (not harming \ref{embcomp})
		\begin{enumerate}[label = $(\bullet)_{2\alph*}$, ref= $(\bullet)_{2\alph*}$]
		\item \label{a3} for $1< \beta \in \dom(p_0) \cup \dom(p_1)$ if $[f_\beta(i) = \ j_\theta] \in p_0(\beta) \cup p_1(\beta)$ holds for no  $i$ then  $[j_\theta  \notin \ran(f_\beta)] \in p_0(\beta) \cap p_1(\beta)$,
		\item 	 \label{b3} whenever $\beta \neq\beta' \in  \dom(p_0) \cup \dom(p_1)$,  $j^* \in \{ \xi_\beta(\rho+1): \ \rho < \lambda \} \cap \{ \xi_{\beta'}(\rho+1): \ \rho < \lambda \} \cap j_\varrho$ and there is no $i$ with $[f_\beta(i) = j^*] \in p_0(\beta) \cup p_0(\beta)$ then $[j^* \notin \ran(f_\beta)] \in p_0(\beta) \cap p_1(\beta)$,
		\item \label{c3} observe that (recalling \ref{form}) whenever $\beta \in  \dom(p_0) \cup \dom(p_1) \setminus S^*$, and $j$ is such that either $[j \notin \ran(f_\beta)] \in p_0(\beta) \cup p_1(\beta)$, or $[f_\beta(i) = j] \in p_0(\beta) \cup p_1(\beta)$ for some $i$, then $j < j_\varrho$,
		\item \label{d3} also observe that $[f_\beta(i) = j] \in p_0(\beta) \cup p_1(\beta)$ implies that $j = j_\varp$ for some $\varp < \varrho$.
		\end{enumerate}
		\stepcounter{bcou}
	\end{enumerate}

	\begin{enumerate}[label = $(\bullet)_{\arabic*}$, ref= $(\bullet)_{\arabic*}$]
		\setcounter{enumi}{\value{bcou}}
		\item \label{2.-}	We construct the increasing sequence $\langle q_\varp: \ \varp < \varrho \rangle$ in $D^*_\alpha$ satisfying 
		$$q_\varp^{[j_\varp]} \geq p_0^{[j_\varp]}, p_1^{[j_\varp]},$$
		\item\label{2.}  and also	for each $\varp < \varrho$ the strict inequality $q_\varp(\beta) \gneq q^{[j_\varp]}_\varp(\beta)$ is only possible if $\beta \in \dom(p_0) \cup \dom(p_1) \setminus \{1\}$ and  $(\delta^\beta_\varp)^- = \max(N^*_\beta \cap (j_\varp+1)) <j_\varp$ hold and  for each such $\beta$ 
		the difference 
		$$q_\varp(\beta) \setminus q^{[j_\varp]}_\varp(\beta) = \left\{ \begin{array}{llr} \{ [f_\beta(i) =  (\delta^\beta_\varp)^-] \}, & \text{if}& [f_\beta(i) =  (\delta^\beta_\varp)^-] \in p_0(\beta) \cup p_1(\beta), \\
		 	\{ \left[(\delta^\beta_\varp)^- \notin \ran(f_\beta)\right] \}, & \text{if} & \left[(\delta^\beta_\varp)^- \notin \ran(f_\beta)\right] \in p_0(\beta) \cup p_1(\beta), \end{array} \right.$$
		  otherwise, if $[(\delta^\beta_\varp)^- \notin \ran(f_\beta)] \notin p_0(\beta) \cup p_1(\beta)$ and for no $i$ we have $[f_\beta(i) =(\delta^\beta_\varp)^-] \in p_0(\beta) \cup p_1(\beta)$, then $q_\varp(\beta) = q_\varp^{[j_\varp]}(\beta)$.
		  (Since for the generic embedding $f_\beta$ $\ran(f_\beta) \subseteq N^*_\beta$ must hold, roughly speaking $q_\varp$ contains all the information from $p_0$ and $p_1$ below $j_\varp$.)
		
		\stepcounter{bcou} 
		\stepcounter{bcou}

	\end{enumerate}

	Now we claim that provided the sequence $\langle q_\varp: \ \varp < \varrho \rangle$ exists there is a common upper bound of $p_0$ and $p_1$. 
	\bcl \label{befcl} The least upper bound of $\langle q_\varp: \ \varp < \varrho \rangle$ (denoted by $q_\varrho \in D^*_\alpha$) can be extended to an upper bound of $p_0$ and $p_1$.
	\ecl
	\begin{PROOF}{Claim \ref{befcl}}
		As the sequence $\langle j_\varp: \ \varp < \varrho\rangle$ has no maximal element, and $q_\varrho \geq q_\varp \geq q^{[j_\varp]}_\varp \geq p_0^{[j_\varp]}, p_1^{[j_\varp]}$ by \ref{2.-}, \ref{2.} clearly $q_\varrho(1) \geq p_0(1), p_1(1)$, and similarly $q_\varrho(\beta) \geq q_0(\beta) \geq p_0(\beta), p_1(\beta)$ for $\beta \in S^*$.
		
		Now fix $1< \beta \in \dom(p_0) \cup \dom(p_1) \setminus S^*$, let $\iota_\beta^- = \sup(N^*_\beta \cap j_\varrho)$. Note that if $\iota_\beta^- = j_\varrho$, then 
		by Definition $\ref{levag}$ $(\lim_{\varp < \varrho} p_0^{[j_\varp]})(\beta) = p_0^{[j_\varrho]}(\beta)$, which is equal to $p_0(\beta)$ by \ref{c3}, and similarly for $p_1$. This implies (recalling  $q_\varp \geq p_0^{[j_\varp]}, p_1^{[j_\varp]}$ by \ref{2.-}, \ref{2.}) that $q_\varrho(\beta) \supseteq p_0(\beta) \cup p_1(\beta)$, as desired.
		
		Now suppose that $\iota_\beta^- < j_\varrho$.
		 Then  clearly $p^{[\delta]}(\beta) = p^{[\iota^-_\beta]}(\beta)$ for every condition $p$ and $\delta \in [\iota^-_\beta, j_\varrho)$. Observe that (by Definition $\ref{levag}$ \ref{levd}, and by the fact that  $[f_\beta(i) = j] \in p(\beta)$ implies $j \in N^*_\beta$) the set
		\begin{enumerate}[label = $(\oint)$, ref = $(\oint)$]
			\item \label{vk} $p_0(\beta) \cup p_1(\beta) \setminus (p_0^{[\iota_\beta^-]}(\beta) \cup p_1^{[\iota_\beta^-]}(\beta))$ consists of only symbols of the form $[j \notin \ran(f_\beta)]$, except maybe $[f_\beta(i) = \iota_\beta^-]$ for a unique $i$ (and then necessarily $\iota_\beta^- = j_\varp$ for some $\varp < \varrho$).
		\end{enumerate}
		
		Then recalling  \ref{2.} whenever $\varp < \varrho$ is such that $j_\varp > \iota^-_\beta$, then
		$$q_\varp(\beta) \supseteq \{[f_\beta(i) = j] \in p_0(\beta) \cup p_1(\beta): \ j \leq \iota^-_\beta\},$$
		similarly
		 	$$q_\varp(\beta) \supseteq \{[j \notin \ran(f_\beta)] \in p_0(\beta) \cup p_1(\beta): \ j \leq \iota^-_\beta\}.$$
		 	
		 	This together with \ref{vk} mean that we only have to add $\{ [j \notin \ran(f_\beta)] \in p_0(\beta) \cup p_1(\beta): \ j > \iota^-_\beta\}$ to $q_\varrho$, which is possible, since
			 	$$ q_\varrho(\beta) \subseteq q_\varrho^{[\iota^-_\beta]}(\beta) \cup \{ [\iota^-_\beta \notin \ran(f_\beta)], [f_\beta(i) = \iota^-_\beta]: \ i \in \delta \}.$$
		 	(In fact using $(\iota_\beta^-, j_\varrho) \cap N_\beta^* = \emptyset$ we could have argued that on coordinate $\beta$ $j$'s not belonging to $N^*_\beta$ are irrelevant in terms of the generic embedding $f_\beta$ and the generic filter.)
		
	\end{PROOF}

	\bcl \label{qsor} There exists a sequence $\langle q_\varp: \ \varp < \varrho \rangle$ satisfying \ref{2.-}, \ref{2.}.
	\ecl
	\begin{PROOF}{Claim \ref{qsor}}
	
	We define $q_{0}$ to be the  upper bound of $p_0^{[j_0]}$ and $p_1^{[j_0]}$ to satisfy \ref{a3}, \ref{b3}: For $\beta \in S^*$ if  $p_0(\beta) = (w_\beta, \name A_{0,\beta})$, $p_1(\beta) = (w_\beta, \name A_{1,\beta})$ then we let $s_0(\beta) = (w, \name B_\beta)$ (where $\name B_\beta$ is the $\bbP^1_\beta$-name satisfying $\Vdash_{\bbP^1_\beta} \name B_\beta = \name A_{0,\beta} \cap \name A_{1,\beta}$). Because of $q_0 = q_0^{[j_\varp]}$ (by \ref{2.-}) and recalling \ref{form}/\ref{ggr} for $\gamma=0,1$   $q_0(1)$ can only be the empty condition. Furthermore, for $\beta \in \dom(p_0) \cup \dom(p_1) \setminus S^*$, $\beta >1$ we let 
		\newcounter{bicou} \setcounter{bicou}{0}
	\begin{enumerate}[label = $(\bigtriangleup)_{\arabic*}$, ref= $(\bigtriangleup)_{\arabic*}$]
		\setcounter{enumi}{\value{bicou}}
		\item $q_0(\beta) = \{ [j^* \notin \ran(f_\beta)] \in p_0(\beta) \cup p_1(\beta): \ j < j_0 \ \wedge \ j \leq \sup(N^*_\beta \cap j_0)  \}$.
		\stepcounter{bicou}
	\end{enumerate}
	So	$q_0$, $q_0^+ \in D^0_\alpha$ in fact belong to $D^*_\alpha$, and  we obviously have \ref{2.-}, \ref{2.}.
	
	Now suppose that $q_\theta$'s are already defined for $\theta < \varp$, and we shall construct $q_\varp$, but we need to deal with limit and successor $\varp$'s differently.
	\begin{enumerate}[label = \underline{Case A}:, ref = \underline{Case A}:]
		\item $\varp$ is limit.
	\end{enumerate}
	  Let $s_\varp = \lim_{\theta < \varp} q_\theta \in D^*_\alpha$, we argue that we can choose a suitable extension of $s_\varp$ to be $q_\varp$.
	For $q_\varp$ we only extend $s_\varp$ on coordinates $\beta \in \dom(p_0) \cup \dom(p_1) \setminus (\{1\} \cup S^*)$. So fix such a $\beta$. First, if $j_\varp \notin N^*_\beta$ (hence $N_\beta^*$ is bounded in $j_\varp$) then we let $q_\varp(\beta) = s_\varp(\beta)$. Second, if $j_\varp \in N^*_\beta$, and it is an accumulation point of $N_\beta^*$, then again we do nothing, just let $q_\varp(\beta) = s_\varp(\beta)$. But if $j_\varp$ is a successor of $(j^\beta_\varp)^- = \max(N^*_\beta \cap j_\varp)$ in $N_\beta^*$, then
	first note that
	\begin{enumerate}[label = $(\bigtriangleup)_{\arabic*}$, ref= $(\bigtriangleup)_{\arabic*}$]
		\setcounter{enumi}{\value{bicou}}
		\item $p_0^{[j_\varp]}(\beta) \cup p_1^{[j_\varp]}(\beta) \subseteq  p_0^{[(j^\beta_\varp)^-]}(\beta) \cup p_1^{[(j^\beta_\varp)^-]}(\beta) \cup \{ [j \notin \ran(f_\beta)]: \ j \geq (j^\beta_\varp)^-\} \cup \{ [f_\beta(i) = (j^\beta_\varp)^-]: \ i <(j^\beta_\varp)^- \}$
		\stepcounter{bicou}
	\end{enumerate}
	(in fact $j$'s between two consecutive element of $N^*_\beta$ are irrelevant in terms of the forcing and the embedding $f_\beta$).
	Moreover, as $\varp$ is limit (and $\langle j_\theta: \ \theta < \varrho \rangle$ is closed by \ref{clo1}) there is $\theta \in \varp$ with $j_\theta \in ((j^\beta_\varp)^-, j_\varp)$, and by \ref{2.-}, \ref{2.} we have
	\begin{enumerate}[label = $(\bigtriangleup)_{\arabic*}$, ref= $(\bigtriangleup)_{\arabic*}$]
		\setcounter{enumi}{\value{bicou}}
		\item $q_\theta(\beta) \subseteq r_\varp(\beta) \subseteq r_\varp^{[(j^\beta_\varp)^-]}(\beta) \cup \{ [(j^\beta_\varp)^- \notin \ran(f_\beta)], [f_\beta(i) = (j^\beta_\varp)^-]: \ i <(j^\beta_\varp)^- \}$.
		\stepcounter{bicou}
	\end{enumerate}
	Again 
	\begin{enumerate}[label = $(\bigtriangleup)_{\arabic*}$, ref= $(\bigtriangleup)_{\arabic*}$]
		\setcounter{enumi}{\value{bicou}}
		\item $r_\varp(\beta) \supseteq p_0^{[(j^\beta_\varp)^-]}(\beta) \cup p_1^{[(j^\beta_\varp)^-]}(\beta)$, and
		\item $r_\varp(\beta) \supseteq \left(p_0(\beta) \cup p_1(\beta)\right) \cap \left\{ [(j^\beta_\varp)^- \notin \ran(f_\beta)], [f_\beta(i) = (j^\beta_\varp)^-]: \ i <(j^\beta_\varp)^- \right\}$.
		\stepcounter{bicou} \stepcounter{bicou}
	\end{enumerate}
	so there is no problem adding $\{ [j \notin \ran(f_\beta)] \in p_0(\beta) \cup p_1(\beta): \ (j^\beta_\varp)^- < j < j_\varp\}$ to $s_\varp(\beta)$ obtaining $q_\varp(\beta)$.
	In each of the cases it is also easy to check \ref{2.}.
		\begin{enumerate}[label = \underline{Case B}:, ref = \underline{Case A}:]
		\item $\varp = \theta +1$.
	\end{enumerate}
	We summarize first which symbols would the $q_\varp(\beta)$'s ($\beta \in \dom(p_0) \cup \dom(p_1)$) have to include in order for $q_\varp$ to satisfy $q_\varp^{[j_\varp]} \geq p_0^{[j_\varp]}, p_1^{[j_\varp]}$, and \ref{2.}. Of course only the case  $\beta \notin S^*$ is relevant.
	\begin{enumerate}[label = $(\bigtriangleup)_{\arabic*}$, ref= $(\bigtriangleup)_{\arabic*}$]
		\setcounter{enumi}{\value{bicou}}
	
		\item for  $\beta = 1$ the set to cover is
		\beeq \label{szine} p_0^{[j_\varp]}(1) \cup p_1^{[j_\varp]}(1) \setminus q_\theta(1) = \{ [j_\theta \ R_\tau \ j] \in p_0(0) \cup p_1(0): \ j < j_\theta, \ \tau < \kappa\}. \eeq
		\stepcounter{bicou} 
	\end{enumerate}
	By \ref{d3}
	\begin{enumerate}[label = $(\bigtriangleup)_{\arabic*}$, ref= $(\bigtriangleup)_{\arabic*}$]
			\setcounter{enumi}{\value{bicou}}
		\item for $1 < \beta \in \dom(p_0) \cup \dom(p_1) \setminus S^*$ the set $q_\varp(\beta)$ has to include the set 
		\beeq \label{jthe} \{ [f_\beta(i) = j_\theta] \in p_0(\beta) \cup p_1(\beta): \ i \in \lambda \} \eeq
		(which is actually either a singleton, or the empty set)
		and
		\beeq \label{nemjthe}  \{ [j \notin \ran(f_\beta)]\in p_0(\beta) \cup p_1(\beta): \ j\in \left((\delta_\theta^\beta)^-,\delta_\varp^\beta)^-\right] \cup \{j_\theta\} \setminus \{j_\varp \} \eeq
		(where $(\delta_\theta^\beta)^- = \sup(N_\beta^* \cap (j_\theta+1))$, $(\delta_\varp^\beta)^- = \sup(N_\beta^* \cap (j_\varp+1))$, possibly $(\delta_\theta^\beta)^- =(\delta_\varp^\beta)^-$). Recall that if $[f_\beta(i) = j_\theta] \in p_0(\beta) \cup p_1(\beta)$ for some $i$, then necessarily $j_\theta \in N^*_\beta$, hence $(\delta_\theta^\beta)^- = j_\theta$.
		\stepcounter{bicou}
	\end{enumerate}
	First we extend $q_\theta$ to a condition $q_\theta^+$ with $q_\theta^+(1)$ including the set in $\eqref{szine}$, and for $\beta \in \dom(p_0) \cup \dom(p_1) \setminus S^*$ $q_\theta^+(\beta)$ including the set in $\eqref{jthe}$.
	\bscl \label{qplu}There exists $q_\theta^+ \geq q_\theta$ in $D_\alpha^*$ with
	\begin{enumerate}[label = $(*)_\alph*$, ref = $(*)_\alph*$]
		\item \label{*a}
	$q_\theta^+(1) \supseteq \{ [j_\theta \ R_\tau \ j] \in p_0(0) \cup p_1(0): \ j < j_\theta, \ \tau < \kappa\},$
\item \label{*b}	for each $0< \beta \notin S^*$ 
$$q^+_\varp(\beta) \supseteq \{ [j \notin \ran(f_\beta)] \in p_0(\beta) \cup p_1(\beta): \ j = (j^\beta_\theta)^-\},$$
$$q^+_\varp(\beta) \supseteq \{ [f_\beta(i) = (j^\beta_\theta)^-] \in p_0(\beta) \cup p_1(\beta):  \ i <(j^\beta_\varp)^- \},$$
\end{enumerate}
while
	\begin{enumerate}[label = $(*)_\alph*$, ref = $(*)_\alph*$]
		\setcounter{enumi}{2}
	\item $q^+_\varp(1) \subseteq q^{+[j_\theta]}_\varp(1) \cup \{[j \ R_\nu \ j_\theta]: \ j < j_\theta, \ \nu < \kappa\},$
	\item and for each $1 < \beta \notin S^*$
	$$q^+_\varp(\beta) \subseteq q^{+[(j^\beta_\theta)^-]}_\varp(\beta) \cup \{[f_\beta(i) = j_\theta]: \ i < j_\theta\} \cup \{[j_\theta \notin \ran(f_\beta)]\}.$$
	\end{enumerate}	
	\escl
	Assuming the subclaim (which guarantees that $q_\theta^+$ satisfies \ref{2.}) we only have to add symbols of the form $[j \notin \ran(f_\beta)]$ (sets in $\eqref{nemjthe}$) to the $q_\theta^+(\beta)$'s to obtain the condition $q_{\theta+1} = q_\varp$ satisfying \ref{2.-} and \ref{2.}, therefore Subclaim $\ref{qplu}$ will finish the proof of Claim $\ref{qsor}$
	\begin{PROOF}{Subclaim \ref{qplu}}(Subclaim \ref{qplu})

	\newcounter{bicouu}
	 \begin{enumerate}[label = $(\blacktriangle)_{\arabic*}$, ref= $(\blacktriangle)_{\arabic*}$]
	 	\setcounter{enumi}{\value{bicouu}}
	 	\item For each fixed $\beta$ where $\beta \in \dom(p_0) \cup \dom(p_1)$ with $[f_\beta(i) = (\delta^\beta_\theta)^-] \in p_0(\beta) \cup p_1(\beta)$ for some $i$ let $i^\beta_\theta$ denote this unique  $i$.
	 		\stepcounter{bicouu}
	\end{enumerate}
	Now observe that
	\begin{enumerate}[label = $(\blacktriangle)_{\arabic*}$, ref= $(\blacktriangle)_{\arabic*}$]
		\setcounter{enumi}{\value{bicouu}}
	 	\item \label{elobb}
	 	for each $\beta$ with $i^\beta_\theta$ defined, for each  $j' < j_\varp$ with $[f_\beta(i') = j'] \in q_\varp(\beta)$ for some $i'$ note that $i' < j' \leq (\delta^\beta_\theta)^{-} \leq j_\varp$ and $i^\beta_\varp <(\delta^\beta_\theta)^{-} \leq j_\varp \in N^*_\beta$, so we can apply Lemma $\ref{tech}$, and thus each condition $q$ in $D^*_\alpha$ can be extended to $q' \in D^*_\alpha$ with $q'^{[j_\varp]}$ deciding the color $c_{\name M_\beta}(i', i^\beta_\varp)$.
	 		 	
		\stepcounter{bicouu}
	\end{enumerate}
	So enumerating all possible pairs $(\beta,i')$ (that are  as in \ref{elobb}) and recalling \ref{*lim} we infer that 
	\begin{enumerate}[label = $(\blacktriangle)_{\arabic*}$, ref= $(\\blacktriangle)_{\arabic*}$]
		\setcounter{enumi}{\value{bicouu}}
		\item 	for some $q^* \geq q_\theta$ the condition $q^{*[j_\varp]} \um \beta \in D^*_\alpha$ decides the color  $c_{\name M_\beta}(i', i^\beta_\varp)$ for all such pairs from $\{ (\beta,i'): \ \beta \in \dom(p_0) \cup \dom(p_1), \  \exists j \: \ [f_\beta(i') = j] \in q_\theta\}$,
		\item repeat this for pairs in $\{ (\beta,i'): \ \exists j \: \ [f_\beta(i') = j] \in q^{*[j_\varp]}\}$, and let $q^{**} \in D^*$ be the condition obtained after countable many such steps,
		\stepcounter{bicouu} \stepcounter{bicouu}
	\end{enumerate}
	so
		\begin{enumerate}[label = $(\blacktriangle)_{\arabic*}$, ref= $(\blacktriangle)_{\arabic*}$]
		\setcounter{enumi}{\value{bicouu}}
			\item \label{szindont}  the condition $q^{**} \in D^*_\alpha$, $q^{**} \geq q_\theta$ with $q^{**[j_\theta]} \um \beta$ deciding the color  $c_{\name M_\beta}(i', i^\beta_\varp)$ for all $(\beta,i') \in \{ (\beta,i'): \ \beta \in \dom(p_0) \cup \dom(p_1), \ \exists j \: \ [f_\beta(i') = j] \in q^{**[j_\varp]}(\beta)$,
			\stepcounter{bicouu}
		\end{enumerate}
	
	Finally recall that  by \ref{2.} $q_\theta(1) = q_\theta^{[j_\theta]}(1)$, and for each $\beta \in \dom(q_\theta) \setminus S^*$ if $ q_\theta(\beta) = q_\theta^{[j_\theta]}(\beta)$ can only be non-empty if $\beta \in \dom(p_0) \cup \dom(p_1)$ (and if it is indeed non-empty then it is a singleton $[(j_\theta^\beta)^- \notin \ran(f_\beta)]$ or $[f_\beta(i) = (j_\theta^\beta)^-]$). 
	\begin{enumerate}[label = $(\blacktriangle)_{\arabic*}$, ref= $(\blacktriangle)_{\arabic*}$]
		\setcounter{enumi}{\value{bicouu}}
		\item  
		 This means that after possibly replacing  $q^{**}_\theta(\beta)$ by $q^{**[j_\theta]}(\beta) \cup q_\theta(\beta)$ using \ref{szindont} it is easy to see that we get a condition  $q^{**} \in D^*_\alpha$ (which still satisfies both \ref{2.} and \ref{szindont}).
		\stepcounter{bicouu}
	\end{enumerate}
	Now we are at the position to construct the desired $q^+_\theta$ as an extension of $q^{**}$. (In order to include the symbols listed in \ref{*a}, and \ref{*b} for $\beta$'s with $(j_\theta^\beta)^- = j_\theta$, but constructing a proper condition in $D^*_\alpha$), our task is to determine the color $\nu(j^*,j_\theta) = c_{M_*}(j^*,j_\theta)$ (i.e. add $[j^* \ R_{\nu(j^*,j_\theta)} \ j_\theta]$ to $q^{**}(1)$) for each $j^*$ and $\beta$ such that
\begin{itemize}
	\item  $[f_\beta(i_\theta^\beta) = j_\theta] \in p_0(\beta) \cup p_1(\beta)$,
	\item and for some $i^*$ $[f_\beta(i^*) = j^*] \in q^{**[j_\theta]}(\beta)$,
\end{itemize} 
	so that $\nu(j^*,j_\theta) = c_{\name M_\beta}(i^*,i_\theta^\beta)$ (this latter value is the color forced by $q^{**[j_\theta]} \um \beta$ by \ref{szindont}). Then adding also the symbols $[f_\beta(i_\theta^\beta) = j_\theta] \in p_0(\beta) \cup p_1(\beta)$ will work.
	
	So fix a pair $j^*,j_\theta$ as above.
	Now we will make use of the preparations \ref{form} and \ref{2..} and show that there are no contradicting demands concerning the value of $\nu(j^*,j_\theta)$.
	We distinguish the following cases.
	\begin{enumerate}[label = Case ($1$):, ref = Case ($1$)]
		\item \label{css1} for some $\nu^* < \kappa$ we have $[j^* \ R_{\nu^*} \ j_\theta] \in p_0(1) \cup p_1(1)$.	 		
	\end{enumerate}  
	Then necessarily $j^* = j_\eta$ for some $\eta < \theta$  (and $j_\eta, j_\theta \in L_0$), and the only option is to   
	\beeq \label{nu*} \text{put } [j_\eta \ R_{\nu^*} \ j_\theta] \in q^+_\varp(1), \eeq 
	i.e.\ define $\nu (j_\eta,j_\theta) = \nu^*$. 
	W.l.o.g.\ we can assume that $[j_\eta \ R_{\nu^*} \ j_\theta] \in p_0(1)$.
	Pick an arbitrary $\beta \in \dom(p_0) \cup \dom(p_1)$ satisfying $[f_\beta(i_\theta^\beta) = j_\theta] \in p_0(\beta) \cup p_1(\beta)$ and for some $i^*$  $[f_\beta(i^*) = j_\eta] \in q^{**}(\beta)$.
	
	If $\beta \in \dom(p_0)$, then by \ref{form}/\ref{fo3} we have fact $j_\eta, j_\theta \in L_0$, which implies that both $[f_\beta(i_\theta^\beta) = j_\theta],[f_\beta(i^*) = j_\eta]  \in p_0(\beta)$, so by \ref{form}/\ref{fo0} $p_0^{[j_\theta]} \um \beta$ forces a value for $c_{\name M_\beta}(i^*,i^\beta_\theta)$. Hence,  $q^{**[j_\theta]} \um \beta \geq q_\theta^{[j_\theta]} \um \beta \geq p_0^{[j_\theta]} \um \beta$ forces the same value for $c_{\name M_\beta}(i^*,i^\beta_\theta)$ (by our hypothesis on $q_\theta$ \ref{2.-}), which is $\nu^*$.
	
	Otherwise, assume that $\beta \notin \dom(p_0)$ (so necessarily $\beta \in \dom(p_1)$ and  $[f_\beta(i_\theta^\beta) = j_\theta] \in p_1(\beta)$, and $j_\theta \in L_1$). Then again by \ref{2..}/\ref{a3} the only way that $[f_\beta(i^*) = j_\eta] \in q_\theta$ can happen for some $i^*$ is when $[f_\beta(i^*) = j_\eta] \in p_1(\beta)$, but then \ref{form}/ \ref{ggr} implies that $j_\eta \in L_1$, so $[j_\eta \ R_{\nu^*} \ j_\theta] \in p_1(\beta)$ is a member of $p_1(\beta)$, too, and then we can proceed as in the case above (i.e.\ arguing that  $p_1^{[j_\theta]} \um \beta \Vdash c_{\name M_\beta}(i^*,i^\beta_\theta) = \nu^*$).
		
		\begin{enumerate}[label = Case ($2$):, ref = Case ($2$)]
		\item \label{css2} for no $\nu^* < \kappa$ have we $[j^* \ R_{\nu^*} \ j_\varp] \in p_0(1) \cup p_1(1)$.	 		
	\end{enumerate}  
	\begin{enumerate}[label = Case ($2A$):, ref = Case ($2A$)]
	\item \label{css2a} $j^* = j_\eta$ for some $\eta < \theta$ (so by \ref{fo1} necessarily $|\{ j_\eta, j_\theta \} \cap (L_0 \setminus L_1)| = |\{ j_\eta, j_\theta \} \cap (L_1 \setminus L_0)| = 1$).
\end{enumerate}  
	We can assume, that $j_\eta \in L_0 \setminus L_1$, $j_\theta \in L_1 \setminus L_0$.
	This means that
		\begin{enumerate}[label = $(\blacktriangle)_{\arabic*}$, ref= $(\blacktriangle)_{\arabic*}$]
		\setcounter{enumi}{\value{bicouu}}
		\item \label{eszr} for no $\beta$ there exists $i$ such that $[f_\beta(i) = j_\eta] \in p_1(\beta)$, and similarly, $[f_\beta(i) = j_\theta] \in p_0(\beta)$ is impossible
		\stepcounter{bicouu}
	\end{enumerate}
	 by our assumption \ref{form}/\ref{ggr} on $p_0$ and $p_1$. So by \ref{2..}/\ref{a3} $[f_\beta(i) = j_\eta] \in q_\theta(\beta)$ is only possible for any $\beta \in \dom(p_0) \cup \dom(p_1)$ 
	if $[f_\beta(i) = j_\eta] \in p_0(\beta) \cup p_1(\beta)$, so this case necessarily $[f_\beta(i) = j_\eta] \in p_0(\beta)$. Summing up, for each $\beta$ with the prospective $q_\theta^+$ forcing $j_\eta \in L_0 \setminus L_1$, $j_\theta \in L_1 \setminus L_0$ to be in the range of $f_\beta$ the only possibility is that 
	\beeq \label{e11}  [f_\beta(i^\beta_\theta) = j_\theta]  \in p_1(\beta) \text{, and}  \eeq
		\beeq \label{e22} \text{for some }i^* \ [f_\beta(i^*) = j_\eta] \in p_0(\beta).   \eeq
		Now we argue that  at most one such $\beta \in \dom(p_0) \cup \dom(p_1)$ may exist (then by \ref{szindont} we can put $[j^* \ R_{\nu^*} \ j_\varp] \in q^+_\theta(\beta)$ with $\nu^* < \kappa$ defined by $q^{**[j_\theta]} \um \beta \Vdash c_{\name M_\beta}(i^*, i^\beta_\theta) = \nu^*$, and we are done).
	
	So assume on the contrary, let $ \beta' \neq \beta''$ be such that $\eqref{e11}$ $\eqref{e22}$ hold.
	Then clearly $\beta', \beta'' \in \dom(p_0) \cap \dom(p_1)$, and $j_\theta, j_\eta \in \{\xi_{\beta'}(\rho+1) : \ \rho < \lambda \} \cap \{\xi_{\beta''}(\rho+1): \ \rho < \lambda \}$, then by our assumption (on all the $p_\gamma$'s) \ref{form}/\ref{fo2} contradicts  \ref{eszr}.
	
	\begin{enumerate}[label = Case ($2B$):, ref = Case ($2B$)]
	\item \label{css2b} $j^* $ is not of the form $ j_\theta$ for any $\theta < \varp$.
\end{enumerate}
	This case we argue that at most one $\beta \in \dom(p_0) \cup \dom(p_1)$ could exist with $[f_\beta(i^\beta_\theta) = j_\theta] \in p_0(\beta) \cup p_1(\beta)$ satisfying that for some $i^*$ $[f_\beta(i^*) =j^*] \in q^{**}(\beta)$. (Then again by \ref{szindont} we can put $[j^* \ R_{\nu^*} \ j_\theta] \in q^+_\theta(\beta)$ with $\nu^* < \kappa$, $q^{**[j_\theta]} \um \beta \Vdash c_{\name M_\beta}(i^*, i^\beta_\theta) = \nu^*$.)
	
	So if there are $\beta' \neq \beta'' \in \dom(p_0) \cup \dom(p_1)$ with 
	\begin{itemize}
	\item $[f_{\beta'}(i^*) =j^*]\in q_\varp(\beta')$  for some $i^*$,
	\item  $[f_{\beta''}(i^{**}) =j^*]  \in q_\varp(\beta'')$ for some $i^{**}$,
	\item $[f_{\beta'}(i^{\beta'}_\varp) =j_\varp] \in p_0(\beta') \cup p_1(\beta')$,
	\item $[f_{\beta''}(i^{\beta''}_\varp) =j_\varp] \in p_0(\beta'') \cup p_1(\beta'')$,
	\end{itemize}
then again as in \ref{css2a} we can get to an easy contradiction (i.e.\ $\beta', \beta'' \in \dom(p_0) \cup \dom(p_1)$, and $j^* \in\{\xi_{\beta'}(\rho+1) : \ \rho < \lambda \} \cap \{\xi_{\beta''}(\rho+1): \ \rho < \lambda \}$,  hence \ref{form}/\ref{ggr} and \ref{2..}/\ref{a3} imply $[j^* \notin \ran(f_\beta)] \in p_0(\beta') \cap p_1(\beta')$, similarly for $\beta''$. Now recall $q^{**} \geq q_\theta$ and \ref{2.}).

\end{PROOF}
\end{PROOF}
\end{PROOF}
\end{PROOF}

Having proven that $\bbP^1_\chi$ (and each $\bbP^1_\alpha$, $\alpha \leq \chi$) is the composition of a $\lambda^+$-cc and a $\kappa^+$-cc forcing, so itself $\lambda^+$-cc, we have \ref{e3}. Moreover, recall Claim $\ref{lambdapcc}$ and that $\bbQ^1_0 = Q(\lambda,\chi,\kappa)$, so $\bbQ^1_0$ does not collapse any cardinal, while $\bbP^1_\chi / \bfG^1_1$ is $\kappa^+$-cc,  $<\kappa$-closed, so $\bbP^1_\chi$ being the composition of the forcings preserving cardinals itself does not collapse any cardinal, we get
\ref{e4}. An easy calculation yields the following.
\bcl \label{meret} For each $\alpha < \chi$ we have $\bfV_1^{\bbP^1_\alpha} \models |\bbQ_\alpha^1| \leq \chi$. Therefore, up to equivalence $\bbP^1_\chi$ is of power $\chi$.
\ecl
\begin{PROOF}{Lemma \ref{meret}}
For $\bbP^1_1 = \bbQ^1_0$ we already know $|\bbQ_1^1|$ by Observation $\eqref{|Q_0|}$. We have to prove the two statements simultaneously by induction on $\alpha$. As $\bbP^1_\chi$ is a $<\kappa$-support iteration, and $\chi^{<\kappa} \leq \chi^\lambda = \chi$, by our premises it is enough to prove for the successor case. So for each $\alpha < \chi$ it is enough to show that $\bfV_1^{\bbP^1_\alpha} \models |\bbQ^1_\alpha| \leq \chi$.
For $\alpha =1$ as $\bbQ^1_1$ is a forcing of a $\kappa$-colored random graph on $\lambda$ with conditions of size $<\kappa$ we get that $|\bbQ^1_1| = \lambda^{<\kappa} \leq \chi$ (in fact $|\bbQ^1_1| = \lambda$).

For $\alpha$ with $1< \alpha \notin S^*$ (so Definition $\ref{itdf}$ \ref{case2}). Again, each condition in $\bbQ^1_{\alpha}$ can be coded by a partial function of size $<\kappa$ on $\lambda$ to $\lambda+1$, so $|\bbQ^1_\alpha| = \lambda^{<\kappa} \leq \chi$.

Finally, for $\alpha \in S^*$ (Definition $\ref{itdf}$ \ref{case3}), $\bbQ^1_\alpha = \bbQ_{D_\alpha}$ is the Mathias type forcing from Definition $\ref{c7}$, where 
$D_\alpha$ is a system of subsets of $V_\kappa^{\bfV_1}$ generating a $\kappa$-complete filter, so $|\bbQ^1_\alpha| \leq (2^{|V_\kappa|})^{\bfV_1^{\bbP^1_\alpha}} = (2^{\kappa})^{\bfV_1^{\bbP^1_\alpha}} \leq \chi$ (because $|\bbP^1_\alpha| = \chi$, $\bbP^1_\alpha$ is $\lambda^+$-cc, and we assumed $(\chi^\lambda)^{\bfV_1} = \chi$).

\end{PROOF}

So now we are ready to complete the definition of $\bbP^1_\chi$ by prescribing the names $\name D_\delta$ ($\delta \in S^*$) and $\name M_\delta$ ($1< \delta \notin S^*$), which are standard easy bookkeeping arguments (using $|\bbP^1_\chi| = \chi$ and the $\lambda^+$-cc), but for the sake of completeness we elaborate. This will prove \ref{e5} and \ref{eu}, so complete the proof of Conclusion $\ref{h28}$.
\bcl \label{cll} The system of $\name D_\delta$'s can be chosen so that for every $\bbP^1_\chi$-name
$\name D$ with $\bfV_1 \ \Vdash_{\bbP^1_\chi} \name D \in [\cP(V_\kappa)]^{\leq \lambda}$ there exists a $\delta \in S^*$, such that for the $\bbP^1_\delta$-name $\name D_\delta$ we have $\Vdash_{\bbP^1_\chi} \name D = \name D_\delta$
\ecl
\begin{PROOF}{Claim \ref{cll}}
	It is obvious that by $\chi^\lambda = \chi$ (so $\cf(\chi) > \lambda$) and the $\lambda^+$-cc for every such $\name D$ there is a nice $\bbP^1_\delta$-name for some $\delta < \chi$. 
	As forcing with the $<\kappa$-closed $\bbP^1_\chi$ does not add new elements to $V_\kappa$ we get that for each $\delta$ there are $\chi^{\kappa\cdot \lambda} = \chi$-many such nice names.
	Also, as $|S^*|= \chi$ we can partition $S^* = \bigcup_{\alpha < \chi} S^*_\alpha$ with $S^*_\alpha \cap \alpha = \emptyset$, $|S^*_\alpha| = \chi$, we can let $\langle \name D_\delta: \ \delta \in S^*_\alpha \rangle$ list the nice names for subsets of $\cP(V_\kappa)$.
\end{PROOF}
A similar calculation yields the following.
\bcl \label{cll2} The system of $\name M_\delta$'s can be chosen so that for every $\bbP^1_\chi$-name for a $\kappa$-colored graph $\name M$ on $\lambda$ there exists a $1< \delta \notin S^*$, such that for the $\bbP^1_\delta$-name $\name M_\delta$ we have $\Vdash_{\bbP^1_\chi} \name M = \name M_\delta$. 
\ecl
\begin{PROOF}{Claim \ref{cll2}}
	Easy.
\end{PROOF}

\end{PROOF}

\bibliographystyle{amsalpha}
\bibliography{shlhetal,1185}


\end{document}